\documentclass {article}
\usepackage[utf8]{inputenc}
\usepackage[T1]{fontenc}
\usepackage {amsmath,amssymb,amscd}
\usepackage {dsfont}
\usepackage {bbm}
\usepackage [a4paper] {geometry}

\usepackage {color}
\numberwithin{equation}{section}
\usepackage[french, english]{babel}

\title{On weak Mellin transforms, second degree characters and the Riemann hypothesis}
\author {Bruno Sauvalle \\ MINES ParisTech, PSL - Research University \\60 Bd Saint-Michel, 75006 Paris, France }
\date{February 2015}

\def \cO {\mathcal {O}}
\def \cN {\mathcal {N}}
\def \cH {\mathcal {H}}

\def\AA{\mathbb{A}}

\def\CC{\mathbb{C}}

\def \xd {\mathfrak {d}}
\def \xf {\mathfrak {f}}

\def\QQ{\mathbb{Q}}
\def\RR{\mathbb{R}}

\def\TT {\mathbb{T}}

\def \cf {\mathbbm{1}} 

\def\ZZ{\mathbb{Z}}

\def  \xd {\mathfrak {d}}

\def \xp {\mathfrak {p}}

\providecommand{\abs}[1]{\lvert#1\rvert}
\providecommand{\norm}[1]{\lVert#1\rVert}

\DeclareMathOperator{\bb}{\quad \Box}
\DeclareMathOperator{\Mell}{Mell}

\DeclareMathOperator{\sign}{sgn}

\DeclareMathOperator {\End}{End}
\DeclareMathOperator{\four}{\mathfrak{F}}

\newcommand {\inv}[1]{\frac{1}{#1}}

\newtheorem {prop} {Proposition}
\newtheorem {definition}{Definition}
\newtheorem {theorem}{Theorem}

\AtEndDocument{\bigskip{\footnotesize%
  \textsc{Bruno Sauvalle, MINES ParisTech, PSL - Research University, 60 Bd Saint-Michel, 75006 Paris, France } \par  
  \textit{E-mail address} \texttt{bruno.sauvalle@mines-paristech.fr}  
}}

\begin{document}
\pagenumbering{arabic}

\maketitle
\begin {abstract}

 We say that a function $f$ defined on $\RR$ or $\QQ_p$ has a well defined weak Mellin transform (or weak zeta integral) if there exists some function $M_f(s)$ so that we have $\Mell(\phi \star f,s) = \Mell(\phi,s)M_f(s)$ for all test functions $\phi$  in  $C_c^\infty(\RR^*)$ or $C_c^\infty(\QQ_p^*)$. We show that if $f$ is a non degenerate second degree character on $\RR$ or $\QQ_p$, as defined by Weil, then the weak Mellin transform of $f$ satisfies a functional equation and cancels only for $\Re(s) = \inv 2$. We then show that if  $f$ is a non degenerate second degree character defined on the adele ring $\AA_\QQ$, the same statement is equivalent to the Riemann hypothesis. Various generalizations are provided.

\end {abstract}

\tableofcontents
\section {Introduction}

It has been a standard practice, since Riemann's 1859 paper \cite {Riemann}, when considering the zeta function
$$ \zeta(s) = \sum_{n \ge 1} \inv {n^s} = \prod_p \inv {1 - \inv {p^s}}$$

to add a  ``gamma factor''$\frac {\Gamma(\frac s2)}{\pi^\frac s2}$ in order to build a ``completed zeta function''
$$ \Xi(s) = \frac {\Gamma(\frac s2)}{\pi^{\frac s2}} \zeta(s) $$
which satisfies the functional equation $\Xi(s) = \Xi(1-s)$

However the link between this gamma factor associated to the ``infinite prime'' and the usual factors $\inv {1 - p^{-s}}$ still remains a mistery.

The ``modern'' approach to the functional equation of the the zeta function, developped by John Tate in his thesis \cite {Tate}, is to introduce the ring of adele $\AA_\QQ$,  consider that the function $\Xi(s)$ is the zeta integral of the function 

$$ \phi(x) = e^{-\pi x_\infty^2} \otimes_p \cf_{\ZZ_p}(x_p) $$
and show that the functional equation of $\Xi(s)$ is a consequence of the equality $\four(\phi) = \phi$ using Poisson summation formula.

This approach, however, does not manage in the same way the ``infinite prime'' and the finite primes : the function $e^{-\pi x^2}$ seems to have nothing to do with the functions $\cf_{\ZZ_p}$ and this ``symmetry breakdown'' is not consistent with the modern idea that all the places of a number field should be put on the same footing.  

In this paper, we show that  the gamma factor at the infinite prime and  the factors $\inv {1 - p^{-s}}$ at finite primes are very close to be the Mellin transforms of the ``same'' function , from an algebraic point of view, on $\RR$ and $\QQ_p$. Finding an algebraic equivalent to $e^{-\pi x^2}$ for the finite places does not seem to be possible, but we observe that that finding an algebraic equivalent to $e^{-\pi i x^2}$ is straightforward : this function can be written as $\psi_\RR(\frac {x^2}2)$ where $\psi_\RR$ is the standard additive character $e^{-2\pi i x}$ on $\RR$, so that its ``algebraic equivalent'' on the finite places could be $\psi_p( \frac {x^2}2)$ where $\psi_p$ is the standard additive character on $\QQ_p$.

The Mellin transforms (or zeta integral) of these functions are not well defined in the usual sense, but it is possible to extend the definition of the Mellin transform using the fact that if $f \star g$ is the multiplicative convolution product, we have $\Mell( f \star g) = \Mell(f) \Mell(g)$ : We say that a function $f$ on $\RR$ or $\QQ_p$ has a well defined ``weak Mellin transform'' at the character $\abs x^s$ if there exists a  function $M_f(s)$ so that for any smooth test function $\phi$ with compact support in the multiplicative group $\RR^*$ or $\QQ_p^*$, we have the equality
$ \Mell( \phi \star f ,s) = \Mell( \phi,s)M_f(s) $.
Using this definition, it is possible to prove that  the weak Mellin transforms of $\psi_\infty(\frac {x^2}2)$ and $\psi_p(\frac {x^2}2)$ are well defined for $\Re(s) >0$ and that we have 
\begin {itemize}
\item $ \Mell ( \psi_\infty (\frac {x^2}2)) =  e^{-s \frac {\pi i}4} \frac {\Gamma( \frac s2)}{\pi^{\frac s2}}$
\item For $p = 2$, $\Mell( \psi_2 ( \frac {x^2}2),s)  = \inv {1 - 2^{-s}}(2^{1-s}(1-2^{s-1}) + e^{\frac {\pi i}4} 2^s (1 - 2^{-s}) ) $
\item For $p \neq 2$, $\Mell( \psi_p( \frac {x^2}2),s) = \inv {1 - p^{-s}} $
\end {itemize}

We can in the similar way define on the ring of adeles $\AA_\QQ$ the weak Mellin transform ( or weak zeta integral) of the global function $f(x) = \psi(\frac {x^2}2) = \psi_\infty(\frac {x_\infty^2}2)\otimes \psi_2 ( \frac {x_p^2}2) .. \otimes  \psi_p(\frac {x_p^2}2)..$, show that it is well defined for $\Re(s) >1$, and equal to the product of the local Mellin transforms. We then get the global formula :
$$ \Mell ( \psi( \frac {x^2}2),s)  =  e^{-s \frac {\pi i}4}(2^{1-s}(1-2^{s-1}) + e^{\frac {\pi i}4} 2^s (1 - 2^{-s}) )\Xi(s)$$

If we note $\Xi_f(s) $ this function,  $\Xi_f$ has an analytic continuation and satisfies a functional equation 

$$ \Xi_f(s) = \overline{ \Xi_f(1 - \bar s) } $$
 It is possible to adapt Tate's Thesis in order to show that this functional equation is a consequence of the fact that the Fourier transform of $\psi(\frac {x^2}2)$ ( considered as a distribution) is equal to its complex conjugate $\bar \psi( \frac {x^2}2)$.

$\Xi_f$ appears to be the product of the completed zeta function $\Xi$ and a non trivial entire function $ e^{-s \frac {\pi i}4}(2^{1-s}(1-2^{s-1}) + e^{\frac {\pi i}4} 2^s (1 - 2^{-s}) )$.
It is natural to investigate the zeroes of this  entire function, and it appears that all these zeroes lie on the line $\Re(s) = \inv 2$.

This result can be generalized, and the main objective of this paper is to study how far this generalization is possible.  Let's  first remark that  the function $\psi(\frac {x^2}2)$ belongs to a class of function called ``second degree character'' by Weil in his celebrated 1964 ``Acta'' paper $\cite {Weil64}$. A continuous function $f$ defined on a locally compact abelian group G with values in the torus $\TT$ is called a second degree character if  the function 
$ f(x+y) f(x)^{-1} f(y)^{-1}$ is a bicharacter, i.e. is a group character as a function of $x$ and  as a function of $y$. On $\AA_\QQ$, the second degree characters are of the form $\psi( \frac a2x^2 + bx)$, whith $a$ and $b$ in $\AA$. Weil showed in \cite{Weil64} that if such a second degree character $f$  is non degenerate (  on $\AA_\QQ$, this is equivalent to  $a$ being an idele), its Fourier transform can be written as $\frac {\gamma_f}{\sqrt \abs a}\bar f(\frac xa)$ where $\gamma_f$ is some scalar ( now called the Weil index) satisfying  $\abs {\gamma_f} = 1$.
Combining this result with Tate's Thesis, we show  that if $f$ is a non degenerate second degree character on $\AA_\QQ$, and $\chi$ a Hecke character on the idele group $\AA_\QQ^\times$, the weak Mellin transform of $f$ at the character $\abs x^s \chi(x)$ is well defined for $\Re(s) >1$ and has an analytic continuation with possible poles at 0 and 1. If we keep the notation $\Xi_f(s,\chi)$ for the analytic continuation of $\Xi_f$,  we have a functional equation 

$$  \Xi_f(s,\chi) = {\gamma_{f}} \abs a^{\inv 2-s}\bar  \chi(a) \overline {  \Xi_{ f }(1-\bar s,  \chi)}$$
For $\Re(s)>1$, $\Xi_f$ can be expressed a an Euler product $ \Xi_f(s,\chi) = \prod_v \zeta_{f_v}(s,\chi)$ where each local function $\zeta_{f_v}= \Mell( f_v,s,\chi_v) $ satisfies a functional equation 
$$  \zeta_{f_v}(s,\chi_v) = \rho(s,\chi_v){\gamma_{f_v}} \abs a_v^{\inv 2-s}\bar  \chi_v(a) \overline {  \zeta_{ f_v }(1-\bar s,  \chi_v)}$$
where $\rho(s,\chi_v)$ are the local constants defined by Tate in his Thesis.
The functions $\zeta_{f_v}$ can be explicitly computed. For example, on $\RR$,  the weak Mellin transform of the second degree character $f_\infty(x) = e^{-2\pi i (\frac a2x^2+bx)}$ can be described using the confluent hypergeometric function $_1F_1$ for $a>0$ as 
$$ \zeta_{f_\infty}(s) =  \frac {e^{-s\frac {\pi i }4}}{\sqrt a^{s}}\frac {\Gamma( \frac {s}2)}{\pi^{\frac {s}2} } { _1 F_1}( \frac s2, \inv 2, \frac {\pi ib^2}{a}) $$
and the functional equation of $\zeta_{f_\infty}$ is  equivalent to Kummer's formula 
 $$ e^x {_1F_1}(a,b,-x) =  {_1F_1}(b-a,b,x) $$

The set of zeroes of $\Xi_f(s,\chi)$ can be split in two classes : the ``local '' zeroes, which are zeroes of one of the local functions $\zeta_{f_v}(s,\chi_v) $ at some place $v$,  and the ``global'' zeroes which are the non trivial zeroes of the Hecke L-function $L(s,\chi)$ associated to the Hecke character $\chi$. 

We prove that all the ``local'' zeroes lie on the line $\Re(s) = \inv 2$, so that the Riemann hypothesis for $L(s,\chi)$ is equivalent to the fact that all the zeroes of $\Xi_f(s,\chi)$ lie on the axis $\Re(s) = \inv 2$.  This result is valid for any number field  and any Hecke L-function with the following restrictions on the choice of the second degree character $f$ : 
\begin {itemize}

\item If $f$ is a second degree character on $\AA_F$ with $F \neq \QQ$, $f$ may not be factorizable ( for example $f(x) = \psi(x \sigma(x))$ where $\sigma$ is the automorphism of $\AA_F$ associated to an element of the Galois group of F over $\QQ$). We have, however to assume that this is the case, i.e. that $f(\alpha)$  can be written as $\prod_v f_v(\alpha_v)$  in order to prove that the weak Mellin transform of $f$ is well defined.
\item If $f$ is factorizable, $f_v$ is not necessarily of the form $\psi_v (\frac a2x^2+bx)$ : It can be written as $\psi_v( \inv 2 x \alpha(x) + bx)$ where $\alpha$ is any  continuous additive function on $F_v$, and we have to put some conditions on $\alpha$ in order to get a functional equation.
\item On $\CC$, the zeroes of the Mellin transform of $f(z) = \psi_\CC(\frac a2z^2+bz)$ do not all lie on the axis $\Re(s) = \inv 2$. However, if we take the second degree  character  $f(z) = \psi_\CC( \frac a2 \abs z^2 + bz)$ , then all the zeroes of $\zeta_f$ lie on the axis $\Re(s) = \inv 2$.
\end {itemize}
 The local part of this theory can be  generalized to second degree characters defined over finite dimensional vector spaces. More precisely, if $f$ is any continuous  function defined on $L^n$, where $L$ is some non discrete locally compact field, and if $\phi$ is a function in $C_c^\infty (GL_n(L))$, we define the operator $f \mapsto  \lambda(\phi)f$ by the formula 
$$ \lambda(\phi)f(v) = \int_{ GL_n(L)} \phi(g)f(g^{-1}v) d^\times g $$
 We show that if $f$ is a non degenerate second degree character on $L^n$ and $\phi \in C_c^\infty (GL_n(L))$, then $\lambda(\phi)f$ 
is a Schwartz function on $L^n$. 
We then consider the maximal compact subgroups $K = GL_n(\cO_L)$ if L is a local field and $K = O(n)$ or $U(n)$ if L is $\RR$ or $\CC$, and the invariant norms on $L^n$ associated to K, i.e. $\norm v = \max ( \abs {v_i})$ or $\norm v = \sqrt { \sum_i \abs {v_i}^2}$. We then define the  Mellin transform of a function in $L^n$ y the formula 
$$ M(f,s) = \int_{L^n} f(v) \norm v^s \frac {dv} {\norm v^n} $$
which is well defined for Schwartz functions when $\Re(s) >0$.
 If $\phi$ is  invariant for the left and right action of $K$, i.e. $\phi(k_1gk_2) = \phi(g)$ for all $k$ in $K$, the function $\lambda(\phi) \norm v^s $ is equal to $\norm v^s$  up to a scalar factor, which we note $\xi_s(\phi)$. It is then immediate that if $\phi$ is a function in $C_c^\infty(GL_n(L))$ satisfying the same condition, we have $$ M( \lambda(\phi)f,s) = \xi_{s-n}(\phi^*)M(f,s) $$
where $\phi^*(g) = \inv {\abs {\det g}} \phi(g^{-1}) $.
Using this formula, it is possible to give a definition for the weak Mellin transform of a non degenerate second degree character on $L^n$ similar to the definition given for $n = 1$.
 If we note $\zeta_f(s)$ the weak  Mellin transform of a non degenerate second degree character $f$ defined on $L^n$ and if the endomorphism associated to $f$ is a dilation, then $\zeta_f(s)$ and $\zeta_f(n-s)$ are related by a functional equation and all the zeroes of $\zeta_{f_v}$ lie on the axis $\Re(s) = \frac n2$ if $L$ is a local field or $L = \RR$.

We then show how the concept of weak Mellin transform on a vector space can be generalized, replacing the function $\norm u^s$, which is naturally associated to the trivial representation of $K$, by similar functions $\nu_{\pi,s}(x)$ naturally associated with spherical representations $(\pi,V_\pi)$ of $K$.  If  we note $\zeta_f(s,\pi)$ the weak  Mellin transform of a non degenerate second degree character $f$ defined on $\RR^n$ in this way and if the endomorphism associated to $f$ is a dilation, then we prove again that the zeroes of $\zeta_{f_v}(s,\pi)$ lie on the axis $\Re(s) = \frac n2$.

\section*{Notations}

We will usually call $L$ a  locally compact field,which we always assume to be non discrete and have characteristic zero, $F$ a number field and $\AA_F$ the ring of adeles associated to $F$.

The standard additive character on  $L$ will be noted $\psi_L$, or $\psi$ if no ambiguity is possible. Explicitly, we have 
$\psi_\RR(x) = e^{-2\pi i x}$, $\psi_\CC(x) = e^{-2\pi i(z+\bar z)}$, $\psi_{\QQ_p}(x) = e^{2\pi i \lambda(x)}$ where $\lambda(x) \in \QQ$ is any  rational number of the form $\frac n{p^k}$ satisfying $\lambda(x) - x \in \ZZ_p$. If $\QQ_\xp$ is a local field of residual characteristic p, we have $\psi_{\QQ_\xp} = \psi_{\QQ_p}\circ Tr_{\QQ_\xp / \QQ_p}$. 
We note $dx$ the Haar measure of $L$ considered as an additive group, and $d^\times x$ the Haar measure of $L^*$ considered as a multiplicative group. These Haar measures are normalized following Tate's Thesis  (cf \cite {Tate}): 
 $dx$ is normalized so that the Fourier inversion formula is valid.
The Fourier transform is defined as 
\begin{equation*} \four(f)(y) = \int_L f(x) \psi(xy) dx \end{equation*}
 $d^\times x$  is normalized in the following way : 
on $\RR$, we write $d^\times x = \frac {dx}{\abs x}$; on $\CC$, we write $d^\times z = \frac {dz}{{\abs z}_\CC} = \frac {dz}{{\abs z}^2} $.

On $\QQ_p$, we write $\ZZ_p$ the ring of integers and $\ZZ_p^\times$ the group of units. The additive measure of $\ZZ_p$ is equal to 1, and we normalize the multiplicative measure so that the measure of $\ZZ_p^\times$ is equal to one. We then have  
$d^\times x = \inv {1 - \inv p} \frac {dx}{\abs x} $

If S is any set, we note $\cf_S$ the characteristic function of this set. For example, $\cf_{\ZZ_p}$ is the characteristic function of $\ZZ_p$.

If $\QQ_\xp$ is a general local field of residual characteristic p, with group of units $\cO_\xp^\times$, we normalize $d^\times x$ so that the measure of $\cO_\xp^\times$ is equal to $(\cN \xd)^{-\inv 2}$ ( where $\xd$ is the different of $\QQ_\xp$). The additive measure of the ring of integers $\cO_{\xp}$ is also set to $(\cN \xd)^{-\inv 2}$ and we have 
\begin {equation} d^\times x = \inv {1 - \inv {\cN \xp}} \frac {dx}{\abs x} \end{equation}
On a $L$-vector space $L^n$, we note $K$ the maximal compact subgroup of $GL_n(L)$, i.e. $K = GL_n(\cO_L)$ if $L$ is local, $K = O(n)$ if $L$ is real, $K = U(n)$ if $L$ is complex.

We note  $\AA^\times_F$ the group of ideles of $F$.
The Mellin transforms of a function $\phi$ defined  on $L^\times$ or $\AA_F^\times$ ( called zeta integral in Tate's Thesis)  are defined for $s \in \CC$ and a multiplicative unitary character $\chi$ as 
$$ \Mell(f,s,\chi) = \int_{L^\times} f(x) \abs {x}^s \chi(x)  d^\times x  $$
or 
\begin{equation*} \Mell(f,s,\chi) = \int_{\AA_F^\times} f(x) \abs x^s \chi(x)  d^\times x \end{equation*}
Note that on $\RR$, the definition we use is different from the usual definition of the Mellin transform, wich considers an  integral from 0 to $\infty$ only.
 We say that a Mellin transform is well defined at $(s,\chi)$ if the associated integral converges absolutely.
The functional equations proved by Tate in his thesis can then be written in the following way : 
\begin {prop} (Tate local functional equation) 
If $\phi$ is a Schwartz function on a locally compact field L, we have for $0<\Re(s)<1$ the equality 
\begin{equation} \Mell(\phi,s,\chi) = \rho(s,\chi) \Mell(\four(f),1-s, \bar \chi) \end{equation}
where $\rho(s,\chi)$ is a function of s and $\chi$ but does not depend on $\phi$
\end {prop}
\begin {prop} (Tate global functional equation)
If $\phi$ is a Schwartz function on $\AA_F$ and $\chi$ a Hecke character, then $\Mell(\phi,s,\chi)$ is well defined for $\Re(s)>1$, and has an analytic continuation to $\CC$ with possible poles at 0 and 1. If we keep the notation $\Mell(\phi,s,\chi)$ for the analytic continuations, we have the equality
\begin{equation} \Mell(\phi,s,\chi) = \Mell(\four(\phi), 1-s, \bar \chi) \end{equation}
\end {prop}
These results are proved in \cite {Tate}

If $G$ is a locally compact abelian group, we note $S'(G)$ the space of tempered distributions on $G$, i.e. the space of continuous linear functional on the Schwartz-Bruhat space $S(G)$. If $\mu$ is an element of $S'(L)$, the weak Fourier transform of $\mu$ is defined, following Schwartz, by the usual formula
\begin{equation} < \four(\mu), \phi> = < \mu, \four(\phi)> \end{equation}

\section{ A connection between Tate's Thesis and Weil 1964 'Acta' paper}

\subsection {second degree characters }
Let's now recall Weil's definition ( cf \cite {Weil64}) of a second degree character : 
a continuous function $f$ defined on a locally compact abelian group G with values in the torus $\TT$ is called a second degree character if  the function 
$ f(x+y) f(x)^{-1} f(y)^{-1}$ is a bicharacter, i.e. is a group character as a function of $x$ and as a function of $y$. For example, the function $e^{-2\pi i ( \frac a2x^2+bx)}$ with $a$ and $b$ in $\RR$ is a second degree character on $\RR$

To any such function, we can associate a continuous morphism $\varrho$ from $G$ to $G^*$
by the formula 
\begin{equation}  f(x+y) f(x)^{-1} f(y)^{-1} = < \varrho (y),x > \end{equation}
and it is clear that this morphism has to be symmetic ( ie $< \varrho ( y),x > = < \varrho( x),y>$ ). 
A second degree character is called non degenerate if the associated morphism $\varrho$ is an isomorphism.
We will always assume in this paper that the second degree characters considered are non degenerate and continuous. 

In \cite {Weil64}, Weil gave two formulae describing the weak Fourier transform of a non degenerate second degree character $f$. Theses formulae will be often used in the following sections.
Proofs of these results are avalaible in \cite{Weil64} or  \cite{Cartier64}. The alternative presentations and proofs we propose below are given in order to show the striking 
connection between these formulae and Tate's Thesis. Indeed, both results can be proved using exactly the same methods.

\subsection { The local functional equation}

\begin {prop} (Weil local functional equation )
\label{WLFE}
Let's consider a non degenerate second degree character f on a locally compact abelian group G and note $\varrho$ the morphism associated to f.  Then there exists a complex number $\gamma_f$ satisfying $\abs {\gamma_f} = 1$ so that the weak Fourier transform of f is equal to $ \frac {\gamma_f}{\sqrt {\abs \varrho}} \bar f ( \varrho^{-1}(x))$
\end {prop}

Remark : 
 $\gamma_f$ is now usually called the Weil index associated to the second degree character $f$.

Proof : 
This proposition can also be written in the following form : 
for any Schwartz function $\phi$ in $S(G)$, we have the formula
\begin{equation} \int_G f(x) \four(\phi)(x) dx = \frac {\gamma_f}{\sqrt {\abs \varrho}} \int_{G^*} \bar f(\varrho^{-1}(x)) \phi(x) dx \end{equation}

Let's first remind that a  proof of the Tate local functional equation has been proposed by Weil in $\cite {Weil66}$  using the concept of eigendistribution by using the following proposition : 
\begin {prop}(Weil)
Let's consider a locally compact field $L$ and a continuous multiplicative character $\chi$ on $L^\times$. Then there exists, up to a scalar factor,  one and only one distribution $\Delta_\chi$ satisfying for all Schwartz function $\phi$ in $S(L)$ the formula
\begin{equation} < \Delta_\chi, \phi(ux)> = \chi(u)< \Delta_\chi, \phi> \end{equation}
\end {prop}
Using this proposition, the proof the Tate local functional equation is a straightforward consequence of the fact that $\four( f(ux))(x) = \inv {\abs u} \four(f)(\frac xu) $.

It appears that the local functional equation for second degree characters can be proved in the same way : 

Let's consider a second degree character $f$ with associated symmetric morphism $\varrho$ so that we have 
\begin{equation} f(x+y) \bar f(x) \bar f(y) = <\varrho(y),x> \end{equation}
We can write this expression as
\begin{equation} f(x+y) <-\varrho(y),x> = f(x)f(y)\end{equation}
Let's introduce for $t$ in $\TT$, $u$ in G and $u^*$ in $G^*$ the operator $tU_{u,u^*}$ acting on functions  defined on G  by the formula $tU_{u,u^*}(f)(x) = tf(x+u)<u^*,x>$. 

we can then write the definition of a second degree character $f$ as 
\begin{equation} t U_{u, -\varrho(u)} (f) = t f(u) f \end{equation}
showing that $f$ is an eigendistribution for the action of the operators $t U_{u, -\varrho(u)}$ for all $u$ in G and $t \in \TT$ and that the associated eigenvalue is $ tf(u)$.

Let's recall that the Heisenberg group associated to G can be described as the set $\TT \times G \times G^*$ equipped with the group law $(t,u,u^*)(t',v,v^*) = (tt'<v^*,u>, u+v, u^*+v^*)$. It is immediate that the map $(t,u,u^*) \mapsto tU_{u,u^*}$ is a representation of this Heisenberg group, which is usually called the Schr\"odinger representation.

We now remark that for $\varrho$ fixed, the set of operators of the form $tU_{u, -\varrho(u)}$ for $t \in \TT$ and $u \in G$ is a commutative group ( because $\varrho$ is symmetric), and it is not difficult to see that this set is the image of a maximal commutative subgroup of the Heisenberg group associated to G. The map  which sends the operator $tU_{u,-\varrho(u)}$ to the scalar $t f(u)$ in $\TT$ is a character of this commutative group.

This character restricts to the identity on the center $(\TT,0,0)$ of the Heisenberg group ( because $f(0) = f(0+0) = f(0)^2 $, so that $f(0) = 1$).
We can then use the following proposition, attributed to Cartier, which appears in \cite {Howe}, and is a consequence of the Stone-Von Neumann theorem:  
\begin {prop}
\label{eigenheisenberg}
Let's consider a maximal commutative subgroup A of the Heisenberg group, and a character $\chi$ of A restricting to the identity on its center. Let's note $\rho$ the Schr\"odinger representation of the Heisenberg group. Then there exists, up to a scalar factor, one and only one distribution $\Delta$ satisfying the formula
\begin{equation} <  \Delta, \rho(a)(f)> = \chi(a) < \Delta, f> \end{equation}
for all a in A and all f in the Schwartz space
\end {prop}
An immediate consequence of this proposition is that any distribution $\Delta$ in $S'(G)$ satisfying for all $u$ in $G$ the functional equation 

\begin{equation}  U_{u, -\varrho(u)} (\Delta) = f(u) \Delta\end{equation}
is equal to the second degree character $f$ up to a scalar factor.

The Weil local functional equation for second degree characters is then a straightforward consequence of the commutation relations between the Fourier transform and the  operators $U_{u,u^*}$ :  

We know that $f$ is a second degree character so that we have 
\begin{equation}f(0) = 1 =  f(u-u) = f(u)f(-u) <-\varrho(u),u> \end{equation}
and  $f(u) <-\varrho(u),u>  = \bar f(-u)$
 
The commutation relation bewteen the Fourier transform and the operators $U_{u,u^*}$ can be described as  : 
\begin{equation} \four \circ U_{u,u^*} = <- u^*,u> U_{u^*,-u} \circ \four  \end{equation}
If  we take the Fourier transform of the formula $U_{u, -\varrho(u)} (f) = f(u) f $, we then get the formula 
\begin{equation} U_{-\varrho(u),-u}( \four(f)) = < u, \varrho(u)>^{-1} f(u) \four(f) =  \bar f(-u) (\four(f)) \end{equation}
or, writing $-\varrho(u) = z$, which is possible because $\varrho$ is assumed to be an isomorphism, 
\begin{equation} U_{z, \varrho^{-1}(z)}( \four(f)) = \bar f(\varrho^{-1}(z) ) (\four(f)) \end{equation}
We then remark that  the second degree character $\bar f( \varrho^{-1}(z))$ satisfies the same functional equation, so that it is equal to $\four(f)$ up to a scalar factor using the proposition \ref{eigenheisenberg}.

In order to show that this scalar factor has norm $\inv {\sqrt {\abs \varrho}}$, we use the Fourier inversion formula and the elementary equality  $\four( \bar \phi)(x) = \overline {\four(\phi)(-x)} $ : Let's suppose that $\four(f) = \lambda \bar f \circ  \varrho^{-1}$ for some $\lambda \in \CC^*$. Then we have 
\begin{equation}f(-x) = \four ( \four (f))x = \lambda  \four ( \bar f \circ \varrho^{-1})(x) =\lambda  \abs \varrho \overline{\four( f)( \varrho(-x))}  = \abs {\lambda}^2 \abs \varrho f(-x)\end{equation}
so that $\abs \lambda = \inv {\sqrt \abs \varrho}$  $\bb$

It should be noted that the local Weil formula can be generalized, replacing the Fourier transform $\four $  with any operator M which normalizes the Heisenberg group. Indeed, the induced automorphism of the Heisenberg group $U \mapsto M U M^{-1}$ will map maximal commutative subgroups to maximal commutative subgroups, so that M will map eigendistributions to eigendistributions.
 
\subsection { The global functional equation}

	If the second degree character is constant on some subgroup ( for example, if a second degree character defined on an adele ring $\AA_F$ is trivial on $F$), Weil also  proved the following result :

\begin {prop} ( Weil global functional equation)

Let's consider a non degenerate second degree character f on $G$, suppose that f is equal to 1 on a closed subgroup $\Gamma$ of  $F$, and assume that the symmetric morphism $\varrho$ associated to f is an isomorphism from $(G,\Gamma)$ to $(G^*, \Gamma^*)$. Then $\gamma(f) = 1$
\end {prop}

 Let's recall that when applied to an adele ring $\AA_F$, with $\Gamma$ equal to the field $F$ embedded in $\AA_F$, this theorem gives a proof of the quadratic reciprocity law on this field.
We can write this proposition in the following way : 
  for any Schwartz function $\phi$ in $S(G)$, we  have the formula
\begin{equation} \int_{G} f(x) \four(\phi)(x) dx = \inv {\abs \varrho} \int_{G^*} \bar f(\varrho^{-1}(x)) \phi(x) dx \end{equation}

In order to prove his own global functional equation, Tate in his Thesis considers the integral 
\begin{equation*} \int_{\AA_F^\times} \phi(x) \chi(x) \abs x^s d^\times x \end{equation*}
He then  splits this integral using a fundamental domain for the action of $F^\times$ on $\AA_F^\times$ and applies the Poisson summation formula to the sum on $F^\times$.

It appears that the same method can be used to prove the Weil global formula if $G$ is an adele ring : Let's suppose that $G$ is an adele ring $\AA_F$,that $\varrho$ is a bijection from $F$ to $F$, take $\Gamma = F$, and identify $G$ with $G^*$ and $\Gamma$ with $\Gamma^*$ using the standard additive character $\psi$ on $\AA_F$. 
Considering that $\abs \varrho = 1$ because $\AA/F$ is compact, we have to prove that if $f$ is our second degree character and $\phi$ any Schwartz function,we have
\begin{equation} \int_{\AA_F} f(x) \four(\phi)(x) dx =  \int_{\AA_F} \bar f(\varrho^{-1}(x)) \phi(x) dx \end{equation}

Let's then split the integral 
\begin{equation*}  \int_{\AA_F} f(x) \phi(x) d x \end{equation*}

 using a fundamental domain  D for the (additive) action of $F$ on $\AA_F$, and apply the Poisson summations formula  : 
\begin{equation} \sum_{x \in F} f(x) = \sum_{x \in F} \four(f)(x)  \end{equation}

Any element of $\AA_F$ can be described in a unique way as $r  = x + \delta$ with $x \in F$ and $\delta \in D$ 
Note that we have, considering that $f$ is trivial on $F$:  
\begin{equation} f(x+\delta) =  f(x) f(\delta)\psi(x \varrho(\delta)) = f(\delta)\psi(x \varrho(\delta))\end{equation}
We now decompose the integral 
\begin{equation} \int_{G} f(x) \phi(x) dx  =   \int_{ \delta \in D}  \sum_{ x \in F}f(x+\delta) \phi(x+\delta)   d\delta\end{equation}

\begin{equation} = \int_D f(\delta)  \lbrace \sum_{x \in F} \psi( x \varrho(\delta))  \phi(x+\delta)  \rbrace  d\delta\end{equation}
We then apply the Poisson summation formula to the inner sum: we consider the Schwartz function $\varphi$
\begin{equation}\varphi(y) =  \psi(y \varrho(\delta))\phi(y+ \delta)  \end{equation}
and compute its Fourier transform :

\begin{equation} \four( \varphi)(y) = \psi(-\delta ( \varrho(\delta)+y) )  \four(\phi)(y+\varrho(\delta)) \end{equation}
The integral is then equal to 
\begin{equation}  \int_{\delta \in D}  \sum_{x \in F}  f(\delta)  \psi( -\delta(  \varrho(\delta) +x)) \four(\phi)(x+ \varrho(\delta))   d\delta\end{equation}
The definition of the second degree character $f$ allows us to write, considering that $-\varrho^{-1}(x) = \delta +( -\delta-\varrho^{-1}(x)) $ 
\begin{equation}  \psi( -\delta(  \varrho(\delta) +x)) = \psi( \delta(  \varrho(-\delta  -\varrho^{-1}(x))) = \bar f (\delta)\bar f ( -\delta - \varrho^{-1}(x) )f( -\varrho^{-1}(x))\end{equation}
Considering that $f( -\varrho^{-1}(x)) = 1$ for $x \in F$, the integral becomes
\begin{equation}  \int_{\delta \in D}  \sum_{x \in F} \bar f( -\delta - \varrho^{-1}(x) ) \four(\phi)(x+ \varrho(\delta))   d\delta\end{equation}
Let's write $\varrho(\delta) = \delta'$. 
It is  immediate that $\varrho$ maps a fundamental domain of $\AA_F$ for the action of $F$ to another fundamental domain D' \begin{equation}  \int_{\delta' \in D'}  \sum_{x \in F} \bar f( -\varrho^{-1}(\delta' + x) ) \four(\phi)(x+ \delta' )   d\delta'\end{equation}

\begin{equation} =   \int_{\AA_F} \bar f ( \varrho^{-1}(x)) \four(\phi)(-x) dx \bb\end{equation}

\section {The weak Mellin transform of second degree characters defined over locally compact fields}

\subsection { definition of the weak Mellin transform}

The definition of the weak Mellin transform that we will use is  different than the one given for the weak Fourier transform and uses the properties of the convolution product.
If $L$ is a non discrete locally compact field, we note $C_c^\infty(L^\times)$  the space of smooth functions with compact support  on $L^\times$ considered as a multiplicative group . The expression ``smooth'' means as usual $C^\infty$ if $L$ is equal to $\RR$ or $\CC$, and locally constant if $L$ is a local field. 

It is clear that $C_c^\infty(L^\times)$ is a  commutative algebras for the convolution product and that if $g$ and $h$ lies in $C_c^\infty (L^\times)$,  we have 
\begin{equation} \Mell( g \star h, s, \chi) = \Mell( g,s, \chi) \Mell(h, s, \chi) \end{equation}

This motivates the following definition : 

\begin {definition}
Lets' consider a function f  defined on $L^\times$ and assume that  $ \phi \star f$ is well defined for all $\phi \in C_c^\infty(L^\times)$.   We say that the function  $M_f(s,\chi)$ is a weak Mellin transform of $f$ at $(s,\chi)$ if for any function $\phi$ in $C_c(L^\times)$ , the Mellin transform of  $  \phi \star f  $  is  well defined at $(s,\chi)$ and satisfies the formula 
\begin{equation} \Mell( \phi \star f,s,\chi) = M_f(s,\chi) \Mell(\phi,s,\chi) \end{equation}
\end {definition}
It is immediate that this definition extends the usual definition of the Mellin transform, so that we will use the same notation for Mellin transforms and weak Mellin transforms.

\subsection {description of  second degree characters on a locally compact field}

Let's now describe more explicitly the non degenerate  second degree characters when $G$ is a non discrete locally compact field $L$ in characteristic zero.  $L$ is then a finite extension of $\QQ_p$ or $\RR$, which we call the base field $L_0$ of $L$.
The group characters of $G$ are simply the  functions of the form $\psi(ax)$ with $a$ in $L$.
Any function of the form $\psi(\inv 2\alpha(x)x)$ where $\alpha$ is a continuous homomorphism of $L$ considered as an additive group  is clearly a second degree character,  and any function of the form $< x, \varrho (x)> $ can be written in this form using the isomorphism between $L$ and $L^*$ given by $a \mapsto \chi_a : \chi_a(x) = \psi(ax)$.
All  second degree characters $f$  can then be written in the form 
\begin{equation} f(x) = \psi( \inv 2 \alpha(x)x+ bx) \label{AB} \end{equation}
where $\alpha$ is any continuous $\ZZ-$module  homomorphism from $L$ to $L$ satisfying $\psi(\alpha(x)y) = \psi(x\alpha(y))$. 
Since we assume $f$ non degenerate,  $\alpha$ is then  also $\QQ$-linear. Using continuity and the fact that the closure of $\QQ$ in $L$ is equal to $L_0$,  we then see  that $\alpha$ has to be $L_0$-linear.  
 For example, if $\sigma$ is an element of the Galois group of $L$, it is immediate that $\psi(\inv 2 a\sigma(x)x + bx)$ is a second degree character, with $\alpha(x) = a\sigma(x)$ and $\abs \alpha = \abs a$.

\subsection { The existence of the weak Mellin transform of a seconde degree character}

Let's  now consider a second degree character $f$ defined on $L$ and suppose that $\alpha$ and $b$ are defined as in equation \ref{AB}. We have a natural left action $\lambda$ of the multiplicative group $L^\times$ on $f$ by the formula  
\begin{equation} \lambda(x) f(y) = f(x^{-1}y) \end{equation}
The integrated form of this action can be written, for $\phi \in C_c^\infty(L)$ as 
\begin{equation} \lambda(\phi) f(y)   = \int_{L^\times} \phi(x) f(x^{-1}y) d^\times x \end{equation} 
We do not use the notation $\phi \star f$ because the domain of $\lambda(\phi)f$ is $L$, not $L^*$

\begin {prop}
\label {regul}
If f is a non degenerate second degree character and $\phi \in C_c^\infty(L^\times)$, then $\lambda(\phi)f$ is a Schwartz function on L

\end {prop}

Proof :

Let's first suppose that $L$ is a local field, $L = \QQ_\xp$. We have to show that $\lambda(\phi)f$ is continuous  ( i.e. locally constant) and has compact support. The continuity is immediate since $\phi$ has compact support in $L^\times$. Let's now show that the support of $\lambda(\phi)f$ is compact.

We get, assuming $y \neq 0$ and using the commutativity of the convolution product
\begin{equation} \lambda(\phi)(f)(y) =  \int_{\QQ_\xp^\times} f(x) \phi(x^{-1} y)  {d^\times x} \end{equation}

We  observe that the integral on $\QQ_\xp^\times$ can be written as an integral on $\QQ_\xp$, using the relation 
$d^\times x = \inv {1 - \inv {\cN\xp}}\frac {dx}{\abs x} $.
Let's note $\phi^*$ the function $\inv {\abs x} \phi(\inv x)$, which is also in $C_c^\infty(\QQ_\xp^\times)$, and write  $\phi^*(0) = 0$, so that $\phi^*$ can also be considered as a Schwartz function on $\QQ_\xp$.
We get 
\begin{equation} \lambda(\phi)(f)(y)  = \inv {1 - \inv {\cN\xp}} \inv {\abs y} \int_{\QQ_\xp} f( x) \phi^*(y^{-1}x) dx  \end{equation}
We now  use the local Weil functional equation ( proposition \ref{WLFE})
\begin{equation} = \inv {1 - \inv {\cN\xp}} \frac {\gamma_f}{\sqrt {\abs \alpha}}\int_{\QQ_\xp} \bar f(\alpha^{-1}(z))\four( \phi^*)(-zy) dz  \end{equation}

$\four(\phi^*)$is Schwartz, so that it has compact support on $\QQ_p$.We can then suppose that its support is included in a ball of radius R. We also know that $f\circ \alpha^{-1}$ is continuous and equal to 1 near zero, so that there exists some $\epsilon$ so that if $\abs z < \epsilon$, then $\bar f \circ \alpha^{-1}(z) = 1$. It is then immediate that if $\abs y > \frac R \epsilon$, then the integral becomes
\begin{equation} \inv {1 - \inv {N\xp}} \frac {\gamma_f}{\sqrt {\abs \alpha}}\int_{\abs{yz} < R} \bar f(\alpha^{-1}(z))\four( \phi^*)(-zy) dz    \end{equation}
\begin{equation}= \inv {1 - \inv {N\xp}} \frac {\gamma_f}{\sqrt {\abs \alpha}}\int_{\abs {z} < \frac R{\abs y}} \bar f(\alpha^{-1}(z))\four( \phi^*)(-zy) dz \end{equation}

considering that  $\frac R{\abs y} < \epsilon$, we get  
\begin{equation} = \inv {1 - \inv {N\xp}} \frac {\gamma_f}{\sqrt {\abs \alpha}}\int_{\QQ_\xp} \four( \phi^*)(-zy) dz \end{equation}
 \begin{equation} = \inv {1 - \inv {N\xp}} \frac {\gamma_f}{\sqrt {\abs \alpha}}\inv{\abs y}  \phi^*( 0) = 0 \end{equation}
 which shows that $\lambda(\phi)f$ has compact support.

Let's now consider the case $L = \RR$. The proposition can be considered as a simple application of the method of stationary phase, but can also be proved directly in the following way, which has the advantage of being fully similar to the local field case : 

We consider the integral 
\begin{equation} \lambda(\phi)(f)(y) = \int_{\RR^*} \phi(x) f(x^{-1}y)d^\times x \end{equation}
Considering that the support of  $\phi$ is compact, we can exchange the derivation and integrations signs and conclude that $\lambda(\phi)(f)$ is $C^\infty$.
In order to show that $\lambda(\phi)(f)(y)$ is  $O(\inv {y^n})$ in $\infty$ for all $n >0$, we write   again $\phi^*(x) = \inv {\abs x} \phi(\inv x)$, $\phi^*(0) = 0$ so that $\phi^*$ is a Schwartz function and replace x with $\inv x $ in the integral  so that it becomes
\begin{equation}  \int_{\RR} f(xy) \phi^*(x)d x \end{equation}
We then use the local Weil formula

\begin{equation} =  \frac {\gamma_{f}}{\abs y \sqrt  {\abs \alpha}}  \int_{\RR} \bar f \circ \alpha^{-1}(\frac xy) \four(\phi^*) (-x ) dx  \end{equation}
We then remark that 
\begin{equation} \int_\RR \four(\phi^*)(-x) dx = \phi^*(0) = 0 \end{equation}
and more generally, using $\four(f')(x) = 2\pi i x \four(f)(x)$,  that 
\begin{equation}  \int_\RR x^n \four(\phi^*)(-x) dx  =(  \inv {2\pi i })^n  \int_\RR  \four((\phi^*)^{(n)})(-x) dx = ( \inv {2\pi i })^n  (\phi^*)^{(n)}(0) = 0 \end{equation}
 so that if $P$ is any Polynomial, the expression is  equal to 
\begin{equation} \frac {\gamma_{f}}{\sqrt {\abs a}\abs y}  \int_{\RR} (\bar f \circ \alpha^{-1}(\frac xy) -P(\frac xy)) \four(\phi^*) (-x ) dx  \end{equation}

Let's choose $P$ to be the polynomial of order n associated to the Taylor expansion of $\bar f \circ \alpha^{-1}$ in zero so that $\bar f \circ \alpha^{-1} (x) -P(x) = O(x^n)$ near zero, and write 
\begin{equation}\delta(x) = \bar f \circ \alpha^{-1} (x) -P( x )\end{equation}
The integral becomes
\begin{equation} \frac {\gamma_{f}}{\sqrt {\abs {a}}\abs y  } \int_{\RR} \delta(\frac xy)  \four(\phi^*) (-x ) dx  \end{equation}
Let's suppose for example  $y >1$ and write $x' = \frac x{\sqrt { y}}$, we get 
\begin{equation}  \frac {\gamma_{f}}{\sqrt {\abs {a y}} } \int_{\RR} \delta(\frac {x'}{\sqrt y})  \four(\phi^*) (-x' \sqrt y ) dx'  \end{equation}
and split the integral in according to the condition $\abs {x'} <1$ and $\abs {x'}>1$ : 
\begin{equation} = \frac {\gamma_{f}}{\sqrt {\abs {a y}} }\lbrace  \int_{\abs {x'} <1} \delta(\frac {x'}{\sqrt y})  \four(g) (-x' \sqrt y ) dx'  +  \int_{\abs {x'} >1} \delta(\frac {x'}{\sqrt y})  \four(g) (-x' \sqrt y ) dx' \rbrace \end{equation}
In the first integral, we remark that $\four(\phi^*)$ is bounded ( it is a Schwartz function) and that $\delta(\frac x{\sqrt y}) < K (\frac x{\sqrt y})^n$ for $y$ large enough, so that the expression is bounded by $\frac {K' }{y^{\frac n2}}$ for some constant $K'$. In the second integral, we use the fact that $\delta(\frac x{\sqrt y})$ is bounded by $K(1+ (\frac x{\sqrt y})^n)$ for some K ( $\delta(x)$ is the sum of a polynomial and a function of module 1) that $\four(g)(-x\sqrt y)$ is bounded by $(\inv {x\sqrt y})^{2n}$ for y large enough to get the result

In order to show that the derivatives of $\lambda(\phi)f$ are also fast decreasing, we remark that for $y \neq 0$, we can exchange the integration and the derivation signs in the expression 
\begin{equation} \frac d{dy}( \lambda(\phi))(f)(y) = \frac {d}{dy} \int_{\RR} f(x)\phi(\frac yx) \frac {dx}{\abs x} \end{equation}
so that we get 
\begin{equation} =  \int_{\RR} f(x) \inv x \phi'(\frac yx) \frac {dx}{\abs x} = \inv {y} \lambda( \inv x \phi')(f)(y)\end{equation}
since $\inv x \phi'(x)$ is also in $C_c^\infty(\RR^\times)$, we see that the expression is fast decreasing when $y \rightarrow \infty$

The proof is the same for  the case $L = \CC$ : We consider $\CC$ as a $\RR$ vector space, $z = x+iy$ and observe that the Fourier transform maps multiplication with $x$ or $y$ to differential operators. As a consequence, If a Schwartz function $\phi^*(z) = \phi^*(x+iy)$ has all its derivatives in zero equal to zero, then we have 
\begin{equation} \int_{\CC} x^ny^m \four(\phi^*)(x+iy) dz = 0 \end{equation}
and the proof can be carried in the same way using the real Taylor expansion in zero of $\bar f \circ \alpha^{-1}(z)$ considered as a function of x and y.
 $\bb$\\
 
Remark : this proposition can be extended without difficulty to division rings but not to split simple algebras (i.e. $GL_n(D)$ where D is a division ring and $n \ge 2$).

\begin {prop}
If f is a non degenerate second degree character defined on a locally compact field $L$, then the weak Mellin transform of f is well defined for $\Re(s) >0$.

\end {prop}

Proof :  In order to prove that the weak Mellin transform of  a second degree character $f$ is well defined, we  have to prove that for each pair $(\chi,s)$ with $\Re(s)>0$, there exists a constant M which does not depend of the choice of $\phi$ so that $\Mell( f \star \phi, \chi,s) = M \Mell(\phi,\chi,s) $
 
Using the fact that $\lambda(\phi)(f)$ is a Schwartz function, it is immediate that the Mellin transform of $f \star \phi$ is well defined for $\Re(s) >0$

Let's now show that if we have  two functions $\phi$ and $\mu$ in $C_c^\infty(L^\times)$, we have for $\Re(s) >0$  the equality 
\begin{equation} \Mell( f \star \phi ,\chi, s)\Mell(\mu,\chi,s) = \Mell( \phi,\chi,s) \Mell( f \star \mu ,\chi,s)\end{equation}

Since all the Mellin transforms appearing in this equality are well defined for $\Re(s) >0$,  it is enough to prove that  we have the equality 
\begin{equation}(f \star \phi) \star \mu = \phi \star ( f \star \mu) \end{equation}
The associated double integral is absolutely convergent since both $\phi$ and $\mu$ are in $C_c^\infty (L^\times)$  and $f$ is bounded, so that we can change the order of the integrals and get the result. 

Let's now take any function $\mu$ in $C_c^\infty(L^\times)$ satisfying $\Mell(\mu, s) \neq 0$ and $\Mell( f \star \mu,s) \neq 0$. If no such function exist, then we can say that the weak Mellin transform of $f$ at $(\chi,s)$ is equal to zero and there is nothing else to prove. If we can find such a $\mu$, we then have for any $\phi$ in $C_c^\infty(L^\times)$ the equality 
\begin{equation} \Mell( f\star \phi,\chi,s) = \frac {\Mell(f\star \mu,\chi,s)}{\Mell(\mu,\chi,s)} \Mell(\phi,\chi,s) \end{equation}
which shows that the weak Mellin transform of f is well defined and equal to $\frac {\Mell(f\star \mu,\chi,s)}{\Mell(\mu,\chi,s)}$ $\bb$

If $f$ is a non degenerate second degree character, we will note $\zeta_f(s,\chi)$ the weak Mellin transform of $f$ at the multiplicative character $\abs x^s \chi(x)$.

\subsection {The functional equation of $\zeta_{f}$}

 Let's first give an elementary equality  for $\zeta_{f}$ : 
Using the formula $\Mell(\phi(ax),s) = \abs a^{-s} \bar \chi(a)\Mell(\phi,s,\chi)$ which is valid for $\phi \in C_c^\infty(L^\times)$, It is immediate that we have 

\begin{equation} \zeta_{f(ax)}(s,\chi) =    \abs a ^{-s} \bar \chi(a) \zeta_{f} (s) \end{equation}

Let's  now show that $\zeta_f$ has an analytic continuation :   

\begin {prop}
\label {prop fourier lambda phi}
for $\phi \in C_c^\infty(L^\times)$ and f a second character defined on $L$, The Fourier transform of the Schwartz function $\lambda(\phi)f (y)$ is equal to the Schwartz function $\frac { \gamma_{f}}{\sqrt {\abs \alpha}} \lambda(  \phi^*)  (\bar f \circ \alpha^{-1})(  y)$ 
where the function $\phi^*(x)$ is defined by the formula $\phi^*(x) = \inv {\abs x}\phi( \inv x)$

\end {prop}
Proof : Since we know that $\lambda(\phi)(f)$ is Schwartz, It is enough to show that this is true in the weak sense, i.e. we have to show that for any Schwartz function $\varphi$, we have the equality
\begin{equation}  \int_{L} \lambda(\phi)(f)(y) \four(\varphi)(y) dy = \frac { \gamma_f}{\sqrt {\abs \alpha}}  \int_{L}  \lambda (\phi^* )( \bar f \circ \alpha^{-1})( y) \varphi (y) dy\end{equation}

 The first integral  is equal to

\begin{equation}  \int_{L} \int_{L^\times} \phi(x)f(x^{-1} y)  d^\times x \four( \varphi)(y) dy \end{equation}

The double sum is  absolutely convergent since $\four(\varphi)$ and $\phi$ are summable on $L$ and $L^\times$ , so that we can  exchange the order of the integrals and use the Weil local functional equation : 
\begin{equation}   \int_{L^\times}\lbrace  \int_{L}f(x^{-1} y)   \four( \varphi)(y) dy \rbrace  \phi(x)d^\times x\end{equation}

\begin{equation} = \frac { \gamma_{f}}{\sqrt {\abs \alpha}}\int_{L^\times} (\int_{L} \abs x \bar f( \alpha^{-1}(yx))    \varphi(y ) dy) \phi(x)d^\times x  \end{equation}

writing $t = \inv x $
\begin{equation} = \frac { \gamma_f}{\sqrt {\abs \alpha}}  \int_{L} (  \int_{L^\times}  \inv {\abs t}   \phi(\inv t ) \bar f(\alpha^{-1} ( yt^{-1}))    d^\times t ) \varphi(y ) dy  \end{equation}

\begin{equation} = \frac { \gamma_f}{\sqrt {\abs \alpha}}  \int_{L} \lambda(\phi^*)( \bar f \circ \alpha^{-1})( y) \varphi (y) dy  \bb \end{equation}

This leads to the following formula :  
\begin {prop}
If f is a non degenerate second degree character on $L$, we have for $0<\Re(s)<1$ the formula
\begin{equation} \zeta_f(s,\chi) =   \frac {\gamma_{f}}{\sqrt {\abs \alpha }} \rho(s, \chi) \zeta_{\bar f \circ \alpha^{-1}}(1-s, \bar \chi)\end{equation}
Where $\rho(s,\chi)$ is the local factor appearing in Tate's local functional equation.
\end {prop}
Proof : 
We know by Tate's Thesis that if $\varphi$ is a Schwartz function, then  $\Mell(\varphi,s, \chi) = \rho(s, \chi)\Mell( \four(\varphi), 1-s, \bar \chi) $  for $0<\Re(s)<1$ so that we get 

\begin{equation} \Mell( \lambda(\phi)(f), s, \chi) = \rho(s, \chi) \Mell( \frac {\gamma_{f}}{\sqrt {\abs \alpha }}   { \lambda(\phi^*)(\bar f \circ \alpha^{-1})(  y)}, 1-s,\bar  \chi)\end{equation}
using the definition of weak Mellin transforms, we get that for any function $\phi$ in $C_c^\times(L^\times)$, we have 
\begin{equation}  \Mell(\phi,s,\chi) \zeta_{f}(s, \chi)=  \frac {\gamma_{f}}{\sqrt {\abs \alpha }} \rho(s, \chi)  \Mell( \phi, s, \chi)\zeta_{\bar f \circ \alpha^{-1}}(1-s, \bar \chi)  \end{equation}
and we get the result by choosing any function $\phi$ so that  $\Mell( \phi, s, \chi)  \neq 0$ $\bb$

It is then immediate that $\zeta_{f}$ has an analytic continuation for $\Re(s) \le 0$ with possible poles at the poles of $\rho(s, \chi)$.
The previous formula is not really a functional equation since $\zeta_f$ and $\zeta_{\bar f \circ \alpha^{-1}}$ are not the same functions. We can however get true functional equations from this under some additional hypothesis on $f$. 

Let's for example suppose that  $f$ is of the form 
\begin{equation} f(x) = \psi(\frac a2x^2+bx) \end{equation} 
We then have $\alpha(x) = ax$ and 
\begin{equation} \bar f \circ \alpha^{-1}(x) = \psi( - \inv {2a} x^2 - \frac ba x) = \bar f(\frac xa)\end{equation}
so that we have 
\begin{equation} \zeta_{\bar f \circ \alpha^{-1}}(s,\chi) = \abs a^{s} \chi(a) \bar  \zeta_f(\bar s, \bar \chi) \end{equation}
which leads to the functional equation 
\begin{equation}  \zeta_f(s,\chi) = {\gamma_{f}} \rho(s, \chi)\abs a^{\inv 2-s}\bar  \chi(a) \bar  \zeta_{ f }(1-\bar s,  \chi)\end{equation}

One can also suppose consider an element $\sigma$ of the Galois group of $L$ and a function $f$ of the form 
 
\begin{equation} f(x) = \psi(\frac a2 \sigma(x)x + bx) \end{equation}
with $\sigma(a) = a$ and $\sigma(b) = b$

 We then have $\alpha(x) = a\sigma(x)$ and using $\psi(\sigma(y)) = \psi(y)$  ( because the trace of $\sigma(y)$ is equal to the trace of $y$) and  $\alpha^{-1}(x) = \sigma^{-1}(\frac xa)$ we get  
\begin{equation} \bar f \circ \alpha^{-1}(x) = \bar \psi(\frac a2 \frac xa {\sigma^{-1}(\frac xa)}+ b \sigma^{-1} ( \frac xa) ))  = \bar \psi(\inv {2a} \sigma(x) x + \sigma(b) \frac xa)  \end{equation}

\begin{equation} = \bar \psi( \frac a2 \sigma ( \frac xa) \frac xa + b \frac xa) = \bar f(\frac xa) \end{equation}
which leads to the same functional equation.

\subsection {The weak Mellin transform of $\psi_p(\frac {x^2}2)$ on $\QQ_p$}

Let's now prove the results given in the introduction for the Mellin transform of $\psi_p(\frac {x^2}2)$ on  $\QQ_p$ for $p$  a rational prime : 

\begin {prop}
\label{qp}
\begin {itemize}

\item On $\QQ_p$ with $p \neq 2$, the weak Mellin transform of $\psi_p( \frac {x^2}2)$ at the character $\abs x^s$ is equal to $\inv {1 - p^{-s}}$.

\item On $\QQ_2$, The weak Mellin transform of $\psi_2( \frac {x^2}2)$ at the character $\abs x^s$ is equal to $\inv {1 - 2^{-s}}(2^{1-s}(1-2^{s-1}) + e^{\frac {\pi i}4} 2^s (1 - 2^{-s}) ) $
\end{itemize}
\end {prop}
Proof : 
We take the function $\cf_{\ZZ_p^\times}$ as test function $\phi$, note that the Mellin transform of this test function at any unramified character is equal to 1 for any value of $s$,  compute explicitly $\lambda( \cf_{\ZZ_p^\times})(f)$ and take the Mellin transform of the result:  

For $p \neq 2$, Let' s consider the integral 
\begin{equation} \lambda(\cf_{\ZZ_p^\times})f(y) = \int_{\QQ_p^\times}\cf_{\ZZ_p^\times} (x) \psi(\inv 2 \frac {y^2}{x^2}) d^\times x  = \int_{\QQ_p^\times}\cf_{\ZZ_p^\times} (x) \psi(\inv 2 {y^2}{x^2}) d^\times x\end{equation}
If the valuation of $y$ is positive or zero, it is immediate that the result is equal to 1 because $\psi $ is equal to 1 on $\ZZ_p$. If the valuation of $y$ is negative strict, we replace the integral on $\QQ_p^\times$ by an integral on $\QQ_p$ and  use the local Weil formula (proposition \ref{WLFE})  : 
\begin{equation} = \inv {1 - \inv p} \int_{\QQ_p} \psi( \inv 2  {y^2}{x^2}) (\cf_{\ZZ_p}(x) - \cf_{p\ZZ_p}(x) dx \end{equation}
\begin{equation} = \inv {1 - \inv p} \frac {\gamma_f}{\abs y} \int_{\QQ_p} \psi( -\inv 2  \frac {x^2}{y^2}) ( \cf_{\ZZ_p}(x) - \inv p \cf_{\inv p \ZZ_p}(x)) dx \end{equation}
and we remark that the restriction of the function $\psi( -\inv 2 \frac {x^2}{y^2})$ to $\ZZ_p$ and $\inv p \ZZ_p$ is equal to 1 because the valuation of $y$ is negative strict. The integral is then equal to zero. 
We then have 
\begin{equation} \lambda( \cf_{\ZZ_p^\times})f  = \cf_{\ZZ_p} \end{equation} so that, taking Mellin transform, we get  
\begin{equation} \zeta_f(s) = \inv {1 - p^{-s}}\end{equation}
Lets' now consider the case $p = 2$. We have to compute the integral 
\begin{equation} \lambda(\cf_{\ZZ_2^\times})f (y)  = \int_{\QQ_2^\times} \psi(\inv 2 {y^2}{x^2})\cf_{\ZZ_2^\times} (x) d^\times x\end{equation}
If the valuation of $y$ is equal or greater than 1, the result is 1 because the function $\psi(\inv 2 y^2x^2)$ remains equal to 1. 
Let's now suppose that the valuation of $y$ is zero. We can suppose that $y$ is equal to 1 because $\lambda(\cf_{\ZZ_p^\times})(f)$ is clearly unramified (i.e. invariant under the action of $\ZZ_p^\times$). We then have to compute the integral 
\begin{equation} \int_{\ZZ_2^\times} \psi(\frac {x^2}2)d^\times x \end{equation} 

For $x \in \ZZ_2^\times$, we can write $x = 1 + 2z$ with  $z \in \ZZ_2$ and it is immediate that $ \psi( \frac {x^2}2) = \psi(\inv 2)  = e^{\pi i } = -1$, so that the integral is equal to $-1$.
Let's now suppose that the valuation of $y$ is equal to $-1$, for example $y = \inv 2$.
We then have to compute the integral 
\begin{equation} \int_{\ZZ_2^\times} \psi(\frac {x^2}8)d^\times x \end{equation} 
 any element of $\ZZ_2^\times$ can be written as $x = k + 4z$ with $z \in \ZZ_2$
and $k$ equal to 1 or 3. We  have 
\begin{equation} \psi(\frac {x^2}8) =  \psi(\frac {k^2}8) = \psi(\inv 8) = e^{\frac {\pi i}4} \end{equation}
If the valuation of $y$ is equal to $-2$ or lower, we use the local Weil formula in the same way as for the case $p \neq 2$ and find that the result is zero.
We can then write 
\begin{equation} \lambda(\cf_{\ZZ_2^\times})f = \cf_{2\ZZ_2} - \cf_{\ZZ_2^\times} + e^{\frac {\pi i}4} \cf_{\inv 2 \ZZ_2^\times} \end{equation}
taking Mellin transform, we get 
\begin{equation} \zeta_f(s) = \frac {2^{-s}}{1 - 2^{-s}} -1 + e^{\frac {\pi i }4}2^s = \inv {1 - 2^{-s}}(2^{1-s}(1-2^{s-1}) + e^{\frac {\pi i}4} 2^s (1 - 2^{-s}) ) \bb\end{equation}
\subsection { the value of $\zeta_f(1)$}

\begin {prop} Let's consider some non degenerate second degree character f defined on a locally compact field f, and note $\alpha$ the associated endomorphism.
We have on $\QQ_\xp$ the equalities 
\begin{equation} \zeta_f(1) = \inv {1 - \inv {\cN \xp}} \frac {\gamma_f}{\sqrt {\abs \alpha}} \end{equation}
and on $\RR$ and $\CC$ 
\begin{equation}  \zeta_f(1) =  \frac {\gamma_f}{\sqrt {\abs \alpha}} \end{equation}
\end {prop}
Remark :  1  is then never a zero of $\zeta_f(s)$.\\
Proof : Let's consider the case $L = \QQ_\xp$,  some test function $\phi$ and compute
\begin{equation} \Mell(\lambda(\phi)f,1) = \int_{\QQ_\xp} \lambda(\phi)f(x) \abs x d^\times x \end{equation}

\begin{equation} = \inv {1 - \inv {\cN \xp}}\int_{\QQ_\xp} \lambda(\phi)f(x) dx 
 =  \inv {1 - \inv {\cN \xp}}\four ( \lambda(\phi)f)(0) \end{equation}
using proposition \ref{prop fourier lambda phi}
\begin{equation}= \inv {1 - \inv {\cN \xp}} \frac {\gamma_f}{\sqrt {\abs \alpha}} \lambda ( \phi^*)\bar f \circ \alpha^{-1}(0) \end{equation}
we have seen that $\bar f \circ \alpha^{-1}(0) = 1$
\begin{equation} =  \inv {1 - \inv {\cN \xp}}\frac {\gamma_f}{\sqrt {\abs \alpha}} \int_{\QQ_\xp^\times} \inv {\abs x} \phi(\inv x) d^\times x \end{equation}
\begin{equation} = \inv {1 - \inv {\cN \xp}}\frac {\gamma_f}{\sqrt \alpha} \Mell(\phi,1) \end{equation}

The proof is the same for $\RR$ and $\CC$

$\bb $

The Weil indices associated to second degree characters of the form $\psi(ax^2)$ are explicitly described in \cite{Perrin} for all locally compact fields. For example, the Weil index of the function $\psi_\RR(\frac {x^2}2) = e^{-\pi ix^2}$ is equal to $e^{-\frac {\pi i }4}$. Let's note $\gamma_a$ the Weil index of the second degree character $\psi( \frac a2x^2)$, and consider more general second degree characters :  

\begin {prop}
\label{Weilindex}
The Weil index of the second degree character $f(x) = \psi(\frac a2x^2+bx)$ is equal to $\gamma_f = \gamma_a \psi( - \frac {b^2}{2a})$
\end {prop}
proof : 
We consider the integral 
\begin{equation} \int_L f(x) \phi(x) dx = \int_L \psi(\frac a2x^2+bx) \phi(x) dx\end{equation}
Let's write $\phi_b(x)$ the function $\phi(x)\psi(bx)$. This is again a Schwartz function. We use the Weil local functional equation associated to the second degree character $\psi(\frac a2x^2)$ : 
\begin{equation}= \int_L \psi(\frac a2x^2)\phi_b(x) dx\end{equation}
\begin{equation} = \frac {\gamma_a}{\sqrt {\abs a}}\int_L \psi( - \inv {2a}x^2) \four(\phi_b)(x) dx \end{equation}
\begin{equation} = \frac {\gamma_a}{\sqrt {\abs a}}\int_L \psi( - \inv {2a}x^2) \four(\phi)(x+b) dx \end{equation}
writing $y = x+b$
\begin{equation} = \frac {\gamma_a}{\sqrt {\abs a}}\int_L \psi( - \inv {2a}(y-b)^2) \four(\phi)(y) dy \end{equation}
\begin{equation} = \frac {\gamma_a}{\sqrt {\abs a}}\int_L \psi( - \inv {2a}(y^2-2yb+b^2) \four(\phi)(y) dy \end{equation}
\begin{equation} = \frac {\gamma_a}{\sqrt {\abs a}}\psi( -\frac {b^2}{2a}) \int_L \psi( - ( \frac a2(\frac ya)^2+b \frac ya) ) \four(\phi)(y) dy \bb\end{equation} 
\subsection {The function $\zeta_{a,b}(s)$ considered as a function of $b$}

Let's consider a family of second degree character of the form $f_{a,b}(x) = \psi(\inv2 ax^2+bx)$ and note $\zeta_{a,b}(s)$ the weak Mellin transform at $s$ of $f_{a,b}$. We consider in this subsection the function $\zeta_{a,b}(s)$ as a function of $b$.
Let's note $D_s$ the distribution on $S(L)$ defined for $\Re(s) >0$ by the formula
\begin{equation} < D_s, \phi> = \int_{\QQ_p^\times} \phi(x) \abs x^s d^\times x \end{equation}

We have the following alternative definition of $\zeta_{a,b}(s)$ :

\begin {prop}
\label{fourier}
The function $\zeta_{a,b}(s)$ considered as a function of b ( or, more precisely, as a distribution on the variable b)  is equal to the weak Fourier transform of the distribution $\psi(\frac a2x^2)D_s$.
\end {prop}
Remark : The Fourier transform and the map $\varphi \mapsto \psi(\frac a2x^2)\varphi$ both belongs to a group of unitary operators called the metaplectic group (cf \cite{Weil64}). The function $\zeta_{a,b}(s)$ considered as a function of $b$  is then the image of $D_s$ under the action of a metaplectic operator. Considering that the distribution $D_s$ can be defined, up to a scalar factor, by the fact that it is an eigendistribution for the dilation group ( cf \cite{Weil66}), which is also a subgroup of the metaplectic group, we see that the function $\zeta_{a,b}(s)$, considered as a function of $b$, can also be defined, up to a scalar factor, as an eigendistribution for a subgroup of the metaplectic group conjugate to the dilation group.

Note also  that a consequence of this proposition is that $\zeta_{a,b}(s)$ considered as a function of $b$ is never a square integrable function, and is never the zero function.

Proof : we have to prove that if $\varphi$ is a Schwartz function on F, then 
\begin{equation} \int_{b \in F} \zeta_{a,b}(s)\varphi(b) db = \int_{y \in F^\times} \psi(\frac a2y^2)\four(\varphi)(y)\abs y ^s  d^\times y \end{equation}
Considering that $\four(\phi)$ is a Schwartz function, the right hand side is simply the Mellin transform of the Schwartz function $\psi(\frac a2y^2)\four(\varphi)(y)$ at s.

Let's choose a function $\phi \in C_c^\infty(F^\times) $  so that $\Mell( \phi,s) \neq 0$  and consider the product 
\begin{equation}  (\int_{b \in F} \zeta_{a,b}(s) \varphi(b) db)\Mell( \phi,s) = \int_{b \in F} \zeta_{a,b}(s)\Mell( \phi,s) \varphi(b) db \end{equation}
using the definition of the weak Mellin transform
\begin{equation} = \int_{b \in F} \Mell(\lambda( \phi)f_{a,b},s) \varphi(b) db \end{equation}
\begin{equation} = \int_{b \in F } (\int_{ x \in F^\times} \lambda(\phi)f_{a,b}(x)\abs x^s d^\times x)  \varphi(b)db \end{equation}
We remark that this double integral is absolutely convergent  for $0< \Re(s) <1$ since we have 
\begin{equation} \abs {\lambda(\phi)f_{a,b}(y)} \le  K \end{equation}
and, using Weil formula, 
\begin{equation} \abs {\lambda(\phi)f_{a,b}(y)} \le \frac {K'}{\abs y} \end{equation}
where $K$ and $K'$ do not depend of $b$ . 
As a consequence, we can exchange the integration signs and get 
\begin{equation} \int_{x \in F^\times} ( \int_{b \in F } \lambda(\phi)f_{a,b}(x) \varphi(b) db)\abs x^s d^\times x \end{equation}
We now remark that the inner integral can be described as 
\begin{equation}  \int_{b \in F } \lambda(\phi)f_{a,b}(x) \varphi(b) db =  \int_{b \in F }  \int_{ y \in F^\times}\phi( y) \psi(\inv2 a( y^{-1} x)^2+b y^{-1}x)  f (b)d^\times y db \end{equation}
and that the double integral is again absolutely convergent so that we can again exchange the summation signs and get 
\begin{equation} \int_{y \in F^\times} \phi(y) \psi(\frac a2( y^{-1} x)^2) \four(f)(y^{-1}x) d^\times y \end{equation}
let's reintroduce the operator $\lambda(\phi)$
\begin{equation} = \lambda(\phi) (\psi( \frac a2 x^2 \four(f)(x) )(x) \end{equation}
If we insert this in the former expression, we get 
\begin{equation} \Mell ( \lambda(\phi) \psi( \frac a2 x^2 \four(f)(x) ,s) \end{equation}
\begin{equation} = \Mell(\phi,s) \Mell ( \psi( \frac a2 x^2 \four(f)(x),s) \end{equation}

If $\Re(s) $ is not in the interval $]0,1[$, we use the unicity of the analytic continuation since both expressions are analytic in s $\bb$

 \subsection {The zeroes of $\zeta_f(s,\chi) $ on a local field}

In this section, we will study the zeroes of $\zeta_f(s,\chi)$ on a local field $L = \QQ_\xp$
We split the study in two parts : $\chi = 1$ on the unit  group ( unramified character) and $\chi \neq 1$ on the unit group.

\subsubsection {  unramified character}

Let's consider a local field $\QQ_\xp$ and a unramified character on $\QQ_\xp$, i.e. a character of the form $\abs x^s$

 Let's note $\varpi$ an uniformizer, $\cO$ the ring of integer, $\cO^\times$ the group of units, and $q = \cN \xp = \abs {\varpi}^{-1}$
 We also note $\xd$ the different ideal and define $d$ by the formula $\xd = \xp^d$
 
\begin {theorem}
\label{zerounramified}
Let's consider a non degenerate second degree character on $\QQ_\xp$ of the form $f(x) = \psi(\inv 2ax^2 + bx)$. Then all the zeroes of $\zeta_f(s)$ lie on the line $\Re(s) = \inv 2$
\end {theorem}

Proof : 
We remark that the zeroes of $\zeta_f$ do not change if we replace $f(x)$ with $f(cx)$ with $c \neq 0$. We can then suppose that  $\abs {a}$ is equal to $q^{d}$ or $q^{d-1}$.

It is immediate that $\Mell( \cf_{\cO^\times},s,\chi)$ is equal to the measure of $\cO^\times$, i.e. $(\cN \xd)^{-\inv 2}$
so that we have the formula 
\begin{equation} \zeta_f(s,\chi) = (\cN \xd)^{\inv 2}\Mell ( \lambda(\cf_{\cO^\times})f,s) \end{equation}
Let's then compute 
\begin{equation} \lambda(\cf_{\cO^\times})f(y) = \int_{\QQ_{\xp}^\times}\cf_{\cO^\times}(x)f(\frac yx)  d^\times x \end{equation}
writing $z = \inv x$
\begin{equation} = \inv {1 - \inv q}\int_{\QQ_{\xp}}   \cf_{\cO^\times}(z)f(yz)dz \end{equation}
an element of $\QQ_{\xp}$ is in $\cO^\times$ if and only if it is in $\cO$ but not in $\varpi \cO$
\begin{equation} = \inv {1 - \inv  q}\int_{\QQ_{\xp}} f(yz)  (\cf_{\cO}(z)- \cf_{\cO}(\frac z \varpi) )dz \end{equation}
\begin{equation} =  \inv {1 - \inv  q}\int_{\QQ_{\xp}} ( f(yz) - \inv q f(\varpi yz)  (\cf_{\cO}(z))dz \end{equation}
Let's define the function $\theta_f$ as 
\begin{equation} \theta_f(y) = \int_{\QQ_{\xp}} f(yx)  \cf_{\cO}(x)dx  \end{equation}
we then have 
\begin{equation} \lambda(\cf_{\cO^\times})f(y) = \inv {1 - \inv { q}}(\theta_f(y) - \inv q \theta_f( y\varpi ))\end{equation} so that if the Mellin transform of $\theta$ is well defined at some s with $\Re(s) >0$, we have 
\begin{equation} \Mell( \lambda(\cf_{\cO^\times})f,s) = \inv {1 - \inv { q}}( 1 - q^{s-1}) \Mell( \theta_f,s) \end{equation} 

Considering that 1 is never a zero of $\zeta_f$, we then see that the zeroes of $\zeta_f$ are the same as the zeroes of $\Mell(\theta_f,s)$
 Let's now give an explicit description of $\theta_f$ and $\Mell( \theta_f,s)$.

We know by the definition of the local different and our choice of the haar measure that the Fourier transform of $ \cf_{\cO }$ is $  (\cN \xd)^{-\inv 2} \cf_{\xd^{-1}} $  
so that we get using the local Weil formula ( proposition $\ref{WLFE}$)
\begin{equation} \theta_f(y) = \frac {\gamma_{f}}{\sqrt {\abs a}} \inv {\abs y} \int_{\QQ_{\xp}} \bar f (\frac x{ay}))(\cN\xd)^{-\inv 2} \cf_{\xd^{-1}}(x) )dx \end{equation}

write $z = x \varpi^{d} $ so that $ dz = (\cN \xd)^{-1} dx  $
\begin{equation} =  (\cN \xd)^{\inv 2}\frac {\gamma_f}{\sqrt {\abs a}} \inv {\abs y} \int_{\QQ_{\xp}} \bar f(\frac z{ya\varpi^d}))  \cf_{\cO}(z)   dz \end{equation}
when then have the functional equation 
\begin{equation} \theta_f(y)  = (\cN \xd)^{\inv 2}\frac {\gamma_f}{\sqrt {\abs a}} \inv {\abs y} \bar \theta_f(\inv {y a \varpi^d})\end{equation}

Let's now first suppose that  the valuation of $a$ is equal to $-d$. 
The functional equation becomes, considering that $\theta_f$ is unramified ( i.e. invariant under the action of $\cO^\times$)
\begin{equation} \theta_f(y) = \gamma_{f} \inv {\abs y} \bar \theta_f( \inv y) \end{equation}
In order to compute $\theta_f$, we can then suppose that the valuation of $y$ is $\ge 0$.

The function $g(x) =f(yx)$, if restricted to $x \in \cO$, it is an additive character : If $x$ and $z$ are in $\cO$, we have
\begin{equation} g(x+z)\bar g(x)\bar g(z) = f(yx+yz)\bar f(yx)\bar f(yz) =\psi(ay^2xz) = 1 \end{equation}
Indeed, we have $\abs {a y^2  x z} \le  \abs a = q^d $ which shows that $ay^2xz \in \xd^{-1}$.
Considering the definition of $\theta_f$ and the fact that $\cO$ is an additive group,  we then have 
$ \theta_f(y) = (\cN \xd )^{-\inv 2}$ if and only if $f(xy) = 1$ for all $x \in \cO$, because the integral of an additive character on a compact abelian group is equal to zero if the character is not trivial on this group, or the measure of this  group if the character is trivial.
It is then immediate that if $\theta_f(y) = (\cN \xd )^{-\inv 2}$, then $\theta_f(z) = (\cN \xd )^{-\inv 2}$ for all  $z$ having a higher valuation than $y$. 
We can then write the restriction of $\theta$ to $\cO$ as $(\cN \xd)^{-\inv 2} \cf_{\varpi^k \cO}(y)$ for some $k \ge 0$.
Using the functional equation, we get that if the valuation of $y$ is $\le 0$, we have 
\begin{equation} \theta_f(y) = \gamma_{f} \inv {\abs y} (\cN \xd)^{-\inv 2} \cf_{\varpi^k \cO}(\inv y)\end{equation} 
Let's first suppose that $k = 0$. Note that it implies that $\gamma_{f} = 1$ (using the functional equation with $y = 1$)
We then have, avoiding double counting for $\abs y = 1$ , the following description of $\theta$ :  
\begin{equation} \theta_f(y) = (\cN \xd)^{-\inv 2} ( \cf_{\cO}(y) + (1 - \cf_{\cO}(y))\inv {\abs y} )  \end{equation}
We then see that the Mellin transform of $\theta_f$ is well defined for $0< \Re(s) <1$ and compute that
 \begin{equation}\Mell( \theta_f,s) = (\cN \xd)^{-1} (\inv {1 - q ^{-s}} +  \frac {q ^{s-1}}{1 - q^{s-1}} )\end{equation}
 \begin{equation}= \frac {(\cN \xd)^{-1}}{(1 - q^{-s})(1 - q^{s-1})}(1 - q^{s-1} + q^{s-1}(1 - q^{-s}) )\end{equation}
 \begin{equation} = \frac {(\cN \xd)^{-1}}{(1 - q^{-s})(1 - q^{s-1})} (1 - \inv q) \end{equation}
 We then have 
 \begin{equation}  \zeta_f(s,\chi) = (\cN \xd)^{\inv 2} \Mell( \lambda(\cf_{\cO^\times})f,s) =(\cN \xd)^{\inv 2} \inv {1 - \inv { q}}( 1 - q^{s-1}) \Mell(\theta_f,s) = \end{equation}
 \begin{equation} = (\cN \xd)^{-\inv 2} \inv {1 - q^{-s}} \end{equation}
 and there is no zero.
 
Let's now suppose that $k \ge 1$ We then have 
\begin{equation} \theta_f(y) = ( \cN \xd) ^{-\inv 2} (\cf_{\varpi ^k \cO}(y)+ \gamma_{f} \inv {\abs y} \cf_{\varpi^k \cO}(\inv y) )\end{equation}
and the Mellin transform is 
\begin{equation} \Mell( \theta_f,s) =( \cN \xd) ^{-1}( \frac {q^{-ks}} {1 - q^{-s}}+ \gamma_{f} \frac {q^{k(s-1)}}{1 - q^{s-1}}) \end{equation}
The roots of $\zeta_f(s)$ are then the roots of the equation 
\begin{equation} q^{-ks}(1 - q^{s-1}) + \gamma_{f} q^{k(s-1)}(1 - q^{-s}) = 0 \end{equation}
Let's write $X = q^{s-\inv 2} $ so that $\abs X = 1$ if and only if $\Re(s) = \inv 2$.
The equation becomes
\begin{equation} X^{-k}( 1 - \frac X {\sqrt {q}}) + \gamma_f X^k(1 - \inv X \inv {\sqrt {q}}) = 0 \end{equation}
\begin{equation} \gamma_f X^{2k} -\frac {\gamma_f}{\sqrt q} X^{2k-1} - \frac X{\sqrt q} +1 = 0 \end{equation}
The number of solution of this polynomial equation cannot be greater than $2k$. Let's then show that we have exactly $2k$ solutions on the unit circle : we write $X =  e^{ i \phi}$ so that if we choose some square root of $\gamma_f$,  the equation becomes, after multiplication with $\inv {\sqrt {\gamma_f}}e^{-i k \phi}$
\begin{equation} \sqrt \gamma_f e^{i k\phi} - \frac {\sqrt \gamma_f}{\sqrt q} e^{ i (k-1)\phi} - \frac { \overline {\sqrt {\gamma_f}}}{\sqrt q} e^{ i (1-k)\phi}
 + \overline {\sqrt \gamma_f} e^{-i k\theta }= 0 \end{equation}
 This expression is twice the real part of $ \sqrt \gamma_f e^{i k\phi} - \frac {\sqrt \gamma_f}{\sqrt q} e^{ i (k-1)\phi}$, and it is clear ( because q>1) that when $\phi$ moves from zero to $2 \pi$, this expression has 2k sign changes, so that it has 2k distinct zeroes.

Lets' now suppose that $\abs a = q^{d-1}$. The functional equation of $\theta$ becomes 
\begin{equation} \theta(y) = {\gamma_f}{\sqrt q} \inv {\abs y} \bar \theta( \inv {\varpi y}) \end{equation}
the same argument gives that there exists some $k \ge 0$ so that for $y \in \cO$, we have 
\begin{equation} \theta(y) = (\cN \xd)^{-\inv 2} \delta_{ \varpi^k \cO}(y) \end{equation}
using the functional equation, we then get the following expression for $\theta$, valid for any value of $k \ge 0$ : 
\begin{equation} \theta(y) =  (\cN \xd)^{-\inv 2} (\cf_{ \varpi^k \cO}(y) + \gamma_{f}\sqrt q  \inv {\abs y}\cf_{\varpi^k}( \inv {\varpi y})) \end{equation}
The Mellin transform of $\theta$ is then well defined for $0<\Re(s)<1$ and we have
\begin{equation} \Mell(  \theta,s) =  (\cN \xd)^{-1} (\frac {q^{ -ks}}{1 -q^{-s}} + \frac {\gamma_f}{ \sqrt q} q^s \frac {q^{ k(1-s)}}{1 -q^{s-1}})\end{equation}
And a similar argument allows to prove that all the roots of this equation satisfy $\Re(s) = \inv 2$ $\bb$

\subsubsection { ramified characters}

We consider a ramified character $\chi$ on the unit group of $\QQ_\xp^\times$ and extend it to a character of $\QQ_\xp$ writing $\chi(\varpi) = 1$ ( cf Tate's Thesis).

We then have the following  theorem : 

\begin {theorem}
Let's consider a non degenerate second degree character on $\QQ_\xp$ of the form $\psi( \frac a2x^2+bx)$  and let's assume that $\zeta_f(s,\chi)$ is not identically zero as a function of s. Then all the zeroes of $\zeta_f(s, \chi)$ lie on the line $\Re(s) = \inv 2$
\end {theorem}
Remark : its is clear that if $\chi$ is odd, the weak Mellin transform of $\psi( \frac a2x^2)$ at $(s,\chi)$ is always equal to zero.

Proof : In order to compute $\zeta_f(s,\chi)$  we consider the test function $\phi_{ \bar \chi}(x) =  \chi(x) \cf_{\cO^\times}(x)$
 which satisfy $\Mell( \phi_{\bar \chi}, s, \chi) = (\cN \xd)^{-\inv 2}$ for all values of s.
We note $\xf$ the conductor of $\chi$ and define $n$ so that $\xf = \xp^n$

We recall  from Tate's Thesis  (\cite{Tate}, p 322 )that if $\chi$ is ramified, the local factor $\rho(\chi,s)$ appearing in Tate's local functional equation is  described by the formulae
\begin{equation} \rho(s,\chi) = (\cN (\xf \xd))^{s-\inv 2} \rho_0(\chi) \end{equation}
where the term $\rho_0(\chi)$ satisfies $\abs {\rho_0(\chi)} = 1$ and is described by the formula 
\begin{equation} \rho_0(\chi) = (\cN \xf)^{-\inv 2} \sum_\epsilon \chi(\epsilon) \psi( \frac {\epsilon}{\varpi^{d+n}}) \end{equation}
where $\{\epsilon\}$ is a set of representatives of the cosets of $1 + \xf$ in $\cO^\times$ 
\begin {prop}
\label{fourchi}
The Fourier transform of $\phi_\chi(x) = \chi(x) \cf_{\cO^\times}(x)$ is equal to $ \chi(-1) \rho_0( \chi) q^{-\frac {n+d}2}  \phi_{\bar \chi} ({ \varpi^{n+d}}x) $
\end {prop}

Proof : 
Considering that $\cO^\times$ is a compact multiplicative group, we remark that the Mellin transform of $\phi_{\chi}$ at the character $\chi'$ is equal to zero if $ \chi' \neq \bar \chi $ on $\cO^\times$ and $(\cN \xd)^{-\inv 2}$ if $\chi' = \bar \chi$ on $\cO^\times$. Using Tate's local functional equation, we then see that  the Mellin transform of $\four(\phi_\chi)$at   $\chi'(x) \abs x^s$ is equal to zero for $\chi' \neq  \chi$. For $\chi'= \chi$, we get that 
\begin{equation} \Mell( \four(\phi_\chi), \chi,s) = (\cN \xd)^{-\inv 2}\chi(-1)\rho(s,\chi) = (\cN \xd)^{-\inv 2}(\cN (\xf \xd))^{s-\inv 2} \chi(-1) \rho_0(\chi)\end{equation}
\begin{equation} = q^{-\frac d2} q^{(s-\inv 2)(n+d)}\chi(-1)\rho_0(\chi) =  \chi(-1) \rho_0(\chi)q^{-\frac d2-\frac {n+d}2 +s(n+d)}\end{equation} .

We then observe that the function $  \chi(-1)\rho_0(\chi) q^{-\frac {n+d}2}\phi_{\bar \chi}( \varpi^{n+d}x) $ has the same Mellin transform as $\four(\phi_\chi)$ for all multiplicative characters $\chi$, so that these two functions have to be equal on $\QQ_\xp^\times$. The equality for $x = 0$ is immediate $\bb$

Let's now consider a second degree character of the form 
\begin{equation} f(x) = \psi(\frac a2 x^2+bx) \end{equation}
with $a \neq 0$
We have 
\begin{equation} \lambda(\phi_{\bar \chi})f(y) = \int_{\QQ_\xp^\times}f(\frac yx) \phi_{\bar \chi}(x) d^\times x \end{equation}
\begin{equation} = \inv {1 - \inv q} \int_{\QQ_\xp}f(yz) \phi_{ \chi}(z) dz \end{equation}
using the local Weil formula and the proposition \ref{fourchi}, we get 
\begin{equation} =  \inv {1 - \inv q}\frac {\gamma_f}{\sqrt {\abs a}}\inv {\abs y} q^{-\frac {n+d}2} \chi(-1) \rho_0( \chi) \int_{\QQ_\xp}\bar f(\frac z{ya})\phi_{\bar \chi}(-\varpi^{n+d}z) 
dz \end{equation}
writing $x = \varpi^{n+d}z$, the expression becomes
\begin{equation} =  \inv {1 - \inv q}\frac {\gamma_f}{\sqrt {\abs a}}  \inv {\abs y}q^{\frac {n+d}2}  \rho_0( \chi) \int_{\QQ_\xp}\bar f(\frac { x}{ya\varpi^{n+d}})\phi_{\bar \chi}(x) dx \end{equation}
so that if we note $\theta_{f, \bar \chi}(y) = \lambda(\phi_{\bar \chi})f(y)$, we have the functional equation 
\begin{equation} \theta_{f, \bar \chi}(y) =  {\gamma_f  \rho_0(\chi)}\frac {q^{\frac {n+d}2}}{\sqrt {\abs a}}\inv {\abs y}  \bar \theta_{f, \bar \chi} ( \inv {y a \varpi^{n+d}}) \end{equation}
We can again suppose without loss of generality that the valuation of $a$ is equal to $-n-d$ or $-n-d+1$. 

We write $a = u \varpi^{-n-d +\delta}$ with $\abs u = 1$ and $\delta$ is equal to 0 or 1, so that the functional equation becomes
\begin{equation} \label {FEthetaf} \theta_{f, \bar \chi}(y) = \gamma_\chi  \rho_0(\chi) \frac {q^{\frac \delta 2} }{\abs y} \bar \theta_{f, \bar \chi } (\inv  {u y \varpi^\delta}) \end{equation}
In order to compute $\theta_{f, \bar\chi}(y)$, we can then suppose that the valuation of $y$ is higher or equal to zero. Considering the integral 
\begin{equation} \theta_{f, \bar \chi}(y) = \inv {1 - \inv q} \int_{\cO_\xp^\times}f(yz) \chi(z) dz \end{equation}
We split $\cO_\xp^\times$ in  cosets modulo  the subgroup $1+\xf$, so that $\chi$ is constant on each coset, and choose a representative $z_i$ of each coset. 

\begin{equation} = \inv {1 - \inv q}\sum_i \int_{\xf}f(yz_i(1+x)) \chi(z_i) dx \end{equation}
We then remark, considering that $f$ is a second degree character, that we have for $x \in \xf$
\begin{equation} f(yz_i + yz_ix)\bar f(yz_i) \bar f(yz_ix)  = \psi(ax y^2z_i^2)  = 1\end{equation}
Indeed, if $x$ is in $\xf$, the valuation of $x$ is higher or equal to $n$, so that the valuation of $ax$ is higher or equal to $-d$, and the valuations of $y$ and $z_i$ are positive or zero.
The integral can then be written as 
\begin{equation} \inv {1 - \inv q}\sum_i  \lbrace f(yz_i)  \chi(z_i)\int_{\xf}f(yz_ix) dx  \rbrace \end{equation}
The change of variable $x' = z_ix$ in the last integral shows that it is a constant as a function of $i$, so that we get

\begin{equation} =  \inv {1 - \inv q}( \sum_i   \chi(z_i) f(z_iy) )(\int_{ \xf} f(yx) dx) \end{equation}

We then remark that if $x_1$ and $x_2$ are in $\xf$, then the valuations of $x_1$ and $x_2$ are higher than n, so that the valuation of  $ax_1$ is higher or equal to $- d$ and we have again
\begin{equation} f(x_1 +x_2)\bar f(x_1)\bar f(x_2) = \psi( ax_1x_2) = 1\end{equation}

Considering that the restriction of $f(x)$ to $\xf$ is then an additive character, the restriction of $f(yx)$ to $\xf$ is also an additive character if the valuation of $y$ is positive or zero. As a consequence, there exists a unique $k \in \ZZ$  so that $f(yx)$ is  equal to 1 for $x \in \xf$ if the valuation of $y$ is equal or higher than $k$, and not constant on $\xf$ if the valuation of $y$ is lower than $k$. We then have $\theta_{f, \bar \chi}(y) = 0$ if the valuation of $y$ is lower  (strict) than k.  
 If $k < 0$, then $\theta_{f, \bar \chi}(y) = 0$ thanks to the functional equation and there is nothing else to prove. 

We then suppose $k\ge 0$. We now show that $\theta_{f, \bar \chi}(y) $ is also equal to zero  if the valuation of y is equal or higher than $k+1$

Let's  consider the sum $\sum_i   \chi(z_i) f(z_iy) $. If $f(zy)$ is constant as a function of $z$  for  $z \in \cO^\times$, it is immediate that this sum is zero because $\chi$ is assumed to be ramified. More generally,   we can compare $f(z_iy)$ and $f(z_jy)$ using the formula
\begin{equation} f(yz_j) = f(y(z_j-z_i))f(yz_i)\psi(ay^2z_i(z_j-z_i)) \end{equation}
Let's now suppose that the valuation of $y$ is equal to $k+1$ or higher. Then the definition of $k$ shows that  $f(y(z_i-z_j))$ is equal to 1 if $z_i - z_j \in \inv {\varpi} \xf $.  We also have under the same conditions $\psi(ay^2z_i(z_j-z_i)) = 1$ :  the valuation of $y^2$ is at least equal to 2, the valuation of $a$ is at least equal to $-n-d$ the valuation of $z_i$ is 0 and the valuation of $z_j - z_i$ is assumed to be higher than $n-1$.

We then see that if $z_i$ and $z_j$ are in the same cosets modulo $1+\inv {\varpi}\xf$, then $f(z_iy) = f(z_jy)$. 
Let's then renumber the $z_i$, writing $z_i = z_{j,k}$ where j indicates to which coset of $1+\inv {\varpi}\xf$ it belongs.
We then have 
\begin{equation} \sum_i   \chi(z_i) f(z_iy) =  \sum_{j,k}   \chi(z_{j,k}) f(z_{j,k}y) = \sum_{j}  \alpha_j \sum_k  \chi(z_{j,k}) \end{equation}
with $\alpha_j = f(z_{j,k}y)$ for any choice of $k$\\

However the character $ \chi$ is  constant on the subgroup $1+\xf$ but not on the subgroup $1 + \inv {\varpi} \xf$ ( this is the definition of $\xf$), and for $j$ fixed, the set of $z_{j,k}$ is a full set of representatives of cosets modulo $1+\xf$ inside a coset modulo $1+\frac {\xf} \varpi$.
The sum over $k$ is then equal, up to a scalar, to an integral of a non constant multiplicative character on a subgroup, so that it is equal to zero.

 We then see that if the valuation of $y$ is equal to $k+1$ or higher, then $\theta_{f, \bar \chi}(y) = 0$, and that if it is lower than $k$ and positive, then we also have $\theta_{f, \bar \chi}(y) = 0$.
 
 Considering that the restriction of $\theta_{f, \bar \chi}$ to $y \in \varpi^k \cO^\times$ has to be of the form $C \bar \chi(y)$ for some constant $C = \theta_{f ,\bar \chi} (\varpi^k)$ ( because it is immediate from the definition of $\theta_{f, \bar \chi} $ that if $w$ is a unit, we have $\theta_{f, \bar \chi}(wx) = \bar \chi (w) \theta_{f, \bar \chi}(x))$, we can then write that if the valuation of y is positive or zero we have 
 \begin{equation} \theta_{f, \bar \chi} (y) = C \bar \chi(y) \cf_{ \varpi^k \cO^\times}(y) \end{equation}

Using the functional equation \ref{FEthetaf}, we can then write that for $y \in \QQ_{\xp}$, we have 
\begin{equation} \theta(y) =   C\bar \chi(y) \cf_{\varpi^k \cO^\times}(y) + \gamma_f \bar \rho_0(\chi) \frac {q^{\frac \delta 2}} {\abs y} \chi(\inv {uy})\bar C   \cf_{\varpi^k\cO^\times}(\inv  {u y \varpi ^\delta })\end{equation}
 If $\inv{u y \varpi^\delta}$ is in $\varpi^k\cO^\times$, then $\abs y = q^{k+\delta}$, so that the expression can also be written as

\begin{equation} = \bar \chi(y)C  \lbrace \cf_{\varpi^k \cO^\times}(y) + \gamma_f \bar \rho_0(\chi) \bar  \chi(u)\frac {\bar C}C q^{-k-\frac {\delta}2} \cf_{\varpi^{-k-\delta}\cO^\times}(y)\rbrace\end{equation}
Let's note $\omega$ the espression $\gamma_f \bar \rho_0(\chi) \bar  \chi(u)\frac {\bar C}C $. It is clear that $\abs \omega = 1$
We get 
\begin{equation} \Mell ( \theta,s,\chi) = (\cN \xd )^{-\inv 2}C( q^{-ks} + \omega q^{-k - \frac {\delta}2} q^{(k+\delta)s}) \end{equation}
and it is immediate that the zeroes of this function lie on the axis $\Re(s) = \inv 2$
$\bb$

\subsection {The function $\zeta_f(s,\chi) $ on $\RR$}

 If a second degree character $f$ is defined on $\RR$, then the automorphism $\alpha$  associated to $f$ has to be $\RR$-linear, so that it can be written as $\alpha(x) = ax$. Any non degenerate second degree character $f$ on $\RR$  is then of the form 
\begin{equation} f(x) = \psi(\frac a2x^2+bx)  = e^{-2\pi i (\frac a2x^2+bx)}\end{equation}
with $a \in \RR^\times$ and $b \in \RR$.  
and the functional equation 
\begin{equation}  \zeta_f(s,\chi) = {\gamma_{f}} \rho(s, \chi)\abs a^{\inv 2-s}\bar  \chi(a) \bar  \zeta_{ f }(1-\bar s,  \chi)\end{equation}
is then always valid.
With this description of $f$, we introduce the notation $\zeta_f(s,\chi) = \zeta_{a,b}(s,\chi)$.
We have only two unitary characters $\chi$ on the unit group : the identity and the sign function $\sign(x)$, which we note $\pm(x)$. Note that since $e^{-2\pi i (\frac a2x^2)}$ is even, we have $\zeta_{a,0}(s,\pm) = 0$ for all s.

We have the following description of the weak Mellin transform of a second degree character : 

\begin {prop}
\label {zetareal}
The weak Mellin  transform of the second degree character $\psi(\frac a2x^2+bx)$ at the character $\abs x^s$ for $a>0$ is 
\begin{equation} \zeta_{a,b}(s) =   \frac {e^{-s\frac {\pi i }4}}{\sqrt a^{s}}\frac {\Gamma( \frac {s}2)}{\pi^{\frac {s}2} } { _1 F_1}( \frac s2, \inv 2, \frac {\pi ib^2}{a}) \end{equation}

\end {prop}
Remark : 
 the functional equation of $\zeta_{a,b}$ is  then consistent with Kummer's formula 
\begin{equation} e^x {_1F_1}(a,b,-x) =  {_1F_1}(b-a,b,x) \end{equation}
Indeed, the functional equation can be written

\begin{equation} \frac {e^{-s\frac {\pi i }4}}{\sqrt a^{s}}{ _1 F_1}( \frac s2, \inv 2, \frac {\pi ib^2}{a})  = \gamma_f  \abs a^{\inv 2 -s} \frac {e^{(1-s)\frac {\pi i }4}}{\sqrt a^{1-s}}{ _1 F_1}( \frac {1-s}2, \inv 2, -\frac {\pi ib^2}{a})  \end{equation}
which simplifies to 
\begin{equation} { _1 F_1}( \frac s2, \inv 2, \frac {\pi ib^2}{a})  = \gamma_f   {e^{\frac  {\pi i }4}}{ _1 F_1}( \frac {1-s}2, \inv 2, -\frac {\pi ib^2}{a})  \end{equation}
using Kummer formula, we get 
\begin{equation} {e^{\frac  {-\pi i }4}}e^{\frac {\pi i b^2}a} = \gamma_f   \end{equation}
and this formula is a consequence of proposition \ref{Weilindex}.\\

Proof of proposition \ref{zetareal} : 
We remark that the formula 
\begin{equation} \four( e^{-2\pi i ( \inv 2 ax^2+bx)})  = \frac {e^{\pi i \frac {b^2}a}}{\sqrt{ ai}} e^{2\pi i ( \inv 2 a ( \frac xa)^2+b \frac xa)} \end{equation}
is valid not only for $a \in \RR_+^*$ and $b \in \RR$, but also for $a \in \CC$ with $\Im(a)<0$ and $\Re(a)>0$ and $b \in \CC$.
We can then define $\zeta_{a,b}(s)$ using the same method  for any $a \in \CC$ with $\Im(a)<0$, $\Re(a)>0$ and $b \in \CC$, and looking at the proof of proposition \ref {regul}, it is not difficult to see that $\zeta_{a,b}(s)$ is a continuous function of $a$ provided $a$ does not cross the lines $\Re(a) = 0$ or $\Im(a) = 0$. It is however not difficult to compute $\zeta_{a,0}(s)$ when $\Im(a)<0$ and $\Re(a)>0$ : this is a regular Mellin transform, and we get  
\begin{equation} \zeta_{a,0}(s) =   (\inv {\sqrt a^s})e^{-s \frac {\pi i }4} \frac {\Gamma( \frac s2)} { \pi  ^{\frac s2}} \end{equation}
Thanks to the continuity of $ \zeta_{a,0}(s)$ as a function of $a$ , this formula is also valid for $a \in \RR_+^*$

Let's now suppose that $a$ is fixed in $\RR_+^*$ and consider $\zeta_{a,b}$ as a function of $b$ .  We observe that $\zeta_{a,b}$ considered as a function of $b$, is the limit of the functions $\zeta_{a',b}$, which are analytic in $b \in \CC$, when $a' \in \CC$, $\Im(a')<0$  converges to $a \in \RR_+^*$, and that the convergence is uniform if $b$ stays in some compact set. 

 As a consequence, $\zeta_{a,b}$ is an analytic function of the variable $b$ ( as a uniform limit of complex analytic functions of $b$), so that we can describe it using its Taylor expansion in zero. 

\begin{equation} \zeta_{a,b}(s) = \sum_{k \ge 0} \frac {\partial ^k}{\partial b^k} \zeta_{a,b}(s)(b=0) \frac {b^k}{k!} \end{equation}

We remark that $\zeta_{a,b}(s)$ is an even function of b, so that all the odd derivatives in zero are equal to zero. 
In order to evaluate the even derivatives, we use the following proposition : 

\begin {prop}
\label{exchange}
For s  fixed, $\zeta_{a,b}(s)$ satisfies the  equations 
\begin{equation}  \frac {\partial}{\partial a}  \zeta_f(s) = -\pi i \zeta_f(s+2) \end{equation}

\begin{equation} \frac {\partial^2}{\partial b^2} \zeta_f(s) = (-2\pi i)^2\zeta_f(s+2) \end{equation}
and 
\begin{equation} \frac {\partial}{\partial a} \zeta_{a,b}(s) +\inv {4 \pi i } \frac {\partial^2 }{\partial b^2} \zeta_{a,b}(s) = 0 \end{equation}
\end {prop}
Proof : 
The function $e^{-2\pi i (\frac a2x^2+bx)}$ satisfies the  equations
\begin{equation} \frac {\partial}{\partial a}e^{-2\pi i (\frac a2x^2+bx)} = (-\pi i )x^2e^{-2\pi i (\frac a2x^2+bx)}\end{equation}
\begin{equation} \frac {\partial^2}{\partial b^2}e^{-2\pi i (\frac a2x^2+bx)} = (-2\pi i )x^2e^{-2\pi i (\frac a2x^2+bx)}\end{equation}

so that we mainly have to prove that we can exchange integration and differentation signs.
Let's for example consider the first equation. 
We have to prove that the Mellin transform of $\frac {\partial}{\partial a} \lambda(\phi)e^{-2\pi i (\frac a2x^2+bx)}$ converges absolutely and find a uniform bound for the associated absolute integral. 
We remark that this expression is equal to 
\begin{equation}  \lambda(\phi){\lbrace \frac {\partial}{\partial a}}e^{-2\pi i (\frac a2x^2+bx)}\rbrace\end{equation}
The function $g(x) = \frac {\partial}{\partial a}e^{-2\pi i (\frac a2x^2+bx)} = -\pi i x^2 e^{-2\pi i (\frac a2x^2+bx)} ) $
is $C^\infty$ and its Fourier transform is also $C^\infty$ and can be computed explicitly  using the commutation relation of the Fourier transform and differential operators.

We can then use the proof of proposition \ref{regul}  in the real case,  and show that the Mellin transform of $\lambda( \phi)g $ is well defined for $\Re(s) >0$ with an absolute bound which remains finite if $a$ stays in some compact in $\RR^*$ which does not contain 0. As a consequence, the Mellin transform of the complete expression is also well defined for $\Re(s)>0$, and we have a uniform bound for the integral defining the Mellin transform which allows to exchange the differentiation and Mellin integration signs.

 $\bb$

We then have for $\Re(a)>0$
\begin{equation} \frac {\partial^{2k}}{\partial b^{2k}} \zeta_{a,b}(s)(b = 0) = (-4\pi i)^k \frac {d^k}{da^k}\zeta_{a,b}(s)(b=0) =  (-4\pi i)^k   {e^{ -s \frac {\pi i }4}}\frac {\Gamma( \frac s2)} { \pi  ^{\frac s2}}  \frac d{da^k}( \inv {\sqrt a^s})\end{equation}
\begin{equation} = (-4\pi i)^k   {e^{-s \frac {\pi i }4}}\frac {\Gamma( \frac s2)} { \pi  ^{\frac s2}}  (-\frac s2)(-\frac s2-1)..(-\frac s2-(k-1)) ( \inv {\sqrt a^{s+2k}})\end{equation}
\begin{equation} = (4\pi i)^k   {e^{-s \frac {\pi i }4}}\frac {\Gamma( \frac s2)} { \pi  ^{\frac s2}}  (\frac s2)_{k}  ( \inv {\sqrt a^{s+2k}})\end{equation}
\begin{equation} = \frac  {e^{ -s \frac {\pi i }4}}{\sqrt a^s} \frac {\Gamma( \frac s2)} { \pi  ^{\frac s2}}  \frac {(4\pi i)^k (\frac s2)_k}{a^{ k}}\end{equation}
The Taylor expansion of $\zeta_{a,b}(s)$ becomes
\begin{equation} \zeta_{a,b}(s) =  \frac  {e^{- s \frac {\pi i }4}}{\sqrt a^s} \frac {\Gamma( \frac s2)} { \pi  ^{\frac s2}} \sum_{k \ge 0}  {(\frac s2)_k} \frac { (4\pi i)^k b^{2k}}{a^{k} (2k)!}\end{equation}
Writing $(2k)! = (2.4.6..2k)(1.3.5..2k-1) = (4^k)k! (\inv 2)(\inv2 +1) ..(\inv 2 + k-1) = 4^k k! (\inv 2)_k$, we  recognize a Kummer confluent hypergeometric function $_1F_1$ 
\begin{equation} =   \frac {e^{-s\frac {\pi i }4}}{\sqrt a^{s}}\frac {\Gamma( \frac {s}2)}{\pi^{\frac {s}2} }\sum_{k = 0}^\infty \frac {(\frac s2)_k}{(k)!(\inv 2)_k } (\frac {\pi i b^2}a)^k =  \frac {e^{-s\frac {\pi i }4}}{\sqrt a^{s}}\frac {\Gamma( \frac {s}2)}{\pi^{\frac {s}2} }{_1F_1}( \frac s2, \inv 2, \frac {\pi i b^2}a) \bb \end{equation}

\begin {prop}
The weak Mellin  transform of $\psi_\RR(\frac a2x^2+bx)$ at the character $\sign(x)\abs x^s$ is  
\begin{equation} \zeta_{a,b}(s, \pm) =  - 2\pi  i b  \frac {e^{-(s+1)\frac {\pi i }4}}{\sqrt a^{s+1}}\frac {\Gamma( \frac {s+1}2)}{\pi^{\frac {s+1}2} } { _1 F_1}( \frac {s+1}2, \frac 3 2, \frac {\pi ib^2}{a})\end{equation}

\end {prop}
Proof : We start from the equality 
\begin{equation} \frac {\partial}{\partial b} e^{-2\pi i( \frac a2x^2+bx)} = -2\pi i x e^{-2\pi i ( \frac a2x^2+bx)}\end{equation}
which leads, after exchanging the integration and derivation signs, to the identity
\begin{equation} \frac {\partial}{\partial b}\zeta_{a,b}(s) = -2\pi i \zeta_{a,b}(s+1; \pm) \end{equation}
so that we have 
\begin{equation} \zeta_{a,b}(s, \pm) = -\inv {2\pi i } \frac {\partial}{\partial b}\zeta_{a,b}(s-1)\end{equation} 
\begin{equation} = -\inv {2\pi i }  \frac {e^{-(s-1)\frac {\pi i }4}}{\sqrt a^{s-1}}\frac {\Gamma( \frac {s-1}2)}{\pi^{\frac {s-1}2} } \frac {\partial}{\partial b}{ _1 F_1}( \frac {s-1}2, \inv 2, \frac {\pi ib^2}{a})\end{equation}
We then use the elementary formula 
\begin{equation} \frac {\partial}{\partial z} {_1 F_1 }(\alpha,\beta,z) = \sum_{k \ge 0} \frac {(\alpha)_{k+1}}{(\beta)_{k+1}}\frac {z^k}{k!} =  \frac \alpha \beta {_1F_1} (\alpha+1,\beta+1,z) \end{equation}
to get  the result $\bb$

Let's now look at the location of the zeroes of $\zeta_{a,b}$. 
\begin {prop}
\label{zero}
Let's suppose that  $_1F_1(u,v,z) = 0$ with  z imaginary and $v \in \RR_+^*$. Then we have $\Re(u) = \frac v2$
\end {prop}
Proof : 
It is well known that the function ${_1 F_1}(u,v,z)$ considered as a function of $z$ is a solution of the Kummer differential equation 
\begin{equation} z f''(z)+(v-z)f'(z)- uf = 0\end{equation}
We can suppose without loss of generality that the zero of ${_1 F_1}(u,v,z)$ has a positive imaginary part, so that we can write it as $z = it_0^2$ with $t_0 \in \RR$.

We consider the function for $t \ge 0$
\begin{equation} \Phi(t) =  t^{\alpha}e^{-\frac {it^2}2}{_1 F_1}(u,v,it^2)\end{equation}
with $\alpha = v -\inv 2 $.

Elementary calculations show that the Kummer differential equation becomes 
\begin{equation} \Phi''(t) = \Phi(t) ( -\frac {\alpha^2}t+\frac {\alpha^2}{t^2} -t^2 +2i (2u -v))  \label {kummer2} \end{equation}
consider for $t>0$ the function 
\begin{equation} W(t) = \Phi(t)  \bar \Phi' (t) - \bar \Phi (t) \Phi'(t) \in i \RR\end{equation}
 
 The derivative of $W(t)$ can be computed as 
\begin{equation} W'(t) = \Phi(t) \bar \Phi''(t) - \bar \Phi  (t)  \Phi''(t) \end{equation}

using equation \ref{kummer2}
\begin{equation} = \abs {\Phi(t)^2}(-2i(2\bar u -v)-2i(2u-v) )\end{equation} 
\begin{equation} = -4i  \abs {\Phi(t)^2}(2\Re(u)  -v)\end{equation} 
 We then have 
\begin{equation} W(t_2) - W(t_1) = \int_{t_1}^{t_2}W'(t)dt = -4i (2\Re(u)  -v) \int_{t_1}^{t_2} \abs {\Phi(t)^2}dt  \end{equation}
Considering that $\Phi(t)$ is square integrable on the interval $[0,t_0]$ ( because $v>0$), it is immediate that $W(t)$ converges to some value when $t$ converges to zero, and it is not difficult to check that this value is zero. Considering that we also have $W(t_0) = 0$, we get 
\begin{equation} 0 = W(t_0) - W(0) = -4i (2\Re(u)  -v) \int_{0}^{t_0} \abs {\Phi(t)^2}dt  \end{equation}
which shows that $2\Re(u) - v = 0 \bb$


\begin {theorem}
Let's consider a non degenerate second degree characters f defined on $\RR$, a unitary character $\chi$ on $\{-1,1\}$ and assume that $\zeta_f(s,\chi)$ is not the zero function as a function of s. Then all the zeroes of $\zeta_f(s,\chi)$ lie on the line $\Re(s) = \inv 2$
\end {theorem}
Proof : 
this is an immediate consequence of the previous propositions.

\subsection {The function $\zeta_f(s,\chi) $ on $\CC$}

On $\CC$ we will study the second degree characters of the form $\psi_\CC(\frac a2z^2+bz)$ and $\psi_\CC(\frac a2 \abs z^2 + bz)$.

 The characters on the unit group of $\CC^*$ are of the form $c_n(z) = ( \frac {z}{\abs z})^n$ with $n \in \ZZ$
 If $f$ is a second degree character, we will then note $\zeta_f(s,n)$ the weak Mellin transform of $f$ at the character $(s,c_n)$ 
 \subsubsection { second degree characters of the form $\psi_\CC(\frac a2\abs z^2+bz)$}
We consider second degree characters which can be written as $f(z) = \psi(\frac a2\abs z^2+bz)$ with $a \in \RR_+^*$ and $b \in \CC$. 

 We note $Z_{a,b}(s,n)$ the weak Mellin transform of $\psi(\frac a2\abs z^2+bz)$ at the character $\abs z_\CC^sc_n(z)$
\begin {prop}
The weak Mellin transform of $\psi_\CC(\frac a2\abs z^2+bz)$ at the character $\abs z_\CC^sc_n(z)$ is for $n = 0$
\begin{equation} Z_{a,b}(s,0) =   \frac {e^{-\frac {\pi i }2s}}{a^s} \frac {\Gamma(s)}{(2\pi)^{s-1}}{_1F_1}(s,1,\frac {2\pi i \abs b^2}a) \end{equation}
and for $n>0$
\begin{equation} Z_{a,b}(s,n) =(-1)^n  \frac {e^{-\frac {\pi i }2(s- \frac n2)}}{a^{s }} \frac {\Gamma(s+\frac n2)}{(2\pi)^{s-1}}\inv{n!} (  \frac {\sqrt {2 \pi} }{\sqrt a}\bar b )^n{_1F_1}(s+\frac n2,1+n,\frac {2\pi i \abs b^2}{ a}) \end{equation}
\end {prop}
Remark : 
As a consequence of proposition \ref{zero}, the zeroes of $Z_{a,b}(s,n)$ lie on the axis $\Re(s) = \inv 2$ \\

Proof : 
We have, writing $z = x+iy$ 
\begin{equation} \psi( \frac a2\abs z^2 + bz) = e^{-2\pi i ( a\abs z^2 + bz+\bar b \bar z)}  = e^{-2\pi i ( a\abs z^2 + 2\Re(b)x - 2\Im(b)y)}\end{equation}
We make the same observation as in the real case : We consider the right expression and remark that if we take  $a \in \CC$, $\Im(a)<0$, and replace $\Re(b)$ and $\Im(b)$ with complex values, then we can still define using the same method its weak Mellin transform, which is continuous as a function of $a$, $\Re(b)$ and $\Im(b)$. However, If we suppose that $\Im(a) <0$,  then it has a well defined regular Mellin transform, which is analytic in $\Re(b)$ and $\Im(b)$.
We can then use the same method  as before :  first compute $Z_{a,b}(s,n)$ for $b = 0$ and then use the Taylor expansion of $Z_{a,b}$ at $b = 0$.

Let's first remark that for $n \neq 0$, we have $Z_{a,0}(s,n) = 0$ : the function $\psi(\frac a2\abs z^2)$ is invariant if we replace $z$ with $uz$ with $\abs u = 1$.
For $n = 0$, the computation of $Z_{a,0}(s,0)$ is straightforward and gives 
\begin{equation}Z_{a,0}(s,0) =  \Mell( e^{-2\pi ia \abs z^2}, \abs z^{2s},c_0 ) = \frac {e^{ -s\frac {\pi i}2}}{a^s} \frac {\Gamma(s)}{ (2\pi)^{s-1}} \end{equation}

Let's now suppose that $b \neq 0$ and $n = 0$
 We consider the equation ( using the Wirtinger operator $\frac { \partial}{\partial b}$)
\begin{equation}\frac { \partial}{\partial b} e^{-2\pi i ( \inv 2 a{\abs z}^2+bz+\inv 2 \bar a {\abs  z}^2 + \bar b \bar z)} = -2\pi i z e^{-2\pi i ( \inv 2 a{\abs z}^2+bz+\inv 2 \bar a { \abs z}^2 + \bar b \bar z)} \end{equation} 
Which leads, using  the same kind of argument as in the real case  , to
\begin{equation} \frac { \partial}{\partial b}Z_{a,b}(s,n) = -2\pi i Z_{a,b}(s+\inv 2,n+1) \end{equation}
we also have 
\begin{equation}  \frac { \partial}{\partial \bar b}Z_{a,b}(s,n) = -2\pi i Z_{a,b}(s+\inv 2,n-1) \end{equation}
Let's now write the Wirtinger Taylor expansion of $Z_{a,b}(s)$ near zero : 
\begin{equation} Z_{a,b}(s) = \sum_{ n \ge 0, p \ge 0} \frac {\partial^{p+n}}{\partial^p b \partial^n \bar b } Z_{a,b}(s)(b=0)\frac {b^p\bar b^n}{p!n!} \end{equation}
The Wirtinger derivatives with $p \neq n$ cancel so that we get 
\begin{equation} = \sum_{ n \ge 0} \frac {\partial^{2n}}{\partial^n b \partial^n \bar b } Z_{a,b}(s)(b=0)\frac {\abs b^{2n}}{(n!)^2} \end{equation}
\begin{equation} = \sum_{ n \ge 0} (-2\pi i )^{2n} Z_{a,b}(s+n)\frac {\abs b^{2n}}{(n!)^2} \end{equation}
\begin{equation} =    \frac {e^{-\frac {\pi i }2s}}{a^s} \frac {\Gamma(s)}{(2\pi)^{s-1}} (\sum_{ n \ge 0} (-2\pi i )^{2n}\inv {(2\pi)^n}(-i)^n \frac {(s)_n}{a^n}  \frac {\abs b^{2n}}{(n!)^2} \end{equation}
\begin{equation} =   \frac {e^{-\frac {\pi i }2s}}{a^s} \frac {\Gamma(s)}{(2\pi)^{s-1}} (\sum_{ n \ge 0}(\frac {2\pi i }a)^n (s)_n \frac {\abs b^{2n}}{(n!)^2} \end{equation}

We recognize again a confluent hypergeometric function
\begin{equation} =    \frac {e^{-\frac {\pi i }2s}}{a^s} \frac {\Gamma(s)}{(2\pi)^{s-1}}{_1F_1}(s,1,\frac {2\pi i \abs b^2}a) \end{equation}

In order to compute $Z_{a,b}(s,n)$ for $n > 0$, we use the formula
\begin{equation} Z_{a,b}(s,n) =   \inv {(-2\pi i )^n}\frac { \partial^n}{\partial b^n}Z_{a,b}(s - \frac n2, 0) \bb \end{equation}
%
%
%

 \subsubsection { second degree characters of the form $\psi(\frac a2 z^2+bz)$}

We note $\zeta_{a,b}(s,n)$ the weak Mellin transform of the second degree character $f_{a,b}(z) = \psi_\CC(\frac a2z^2+bz)$.
In order to compute $\zeta_{a,b}(s,n)$, we use the same method as for the real case or for $Z_{a,b}(s)$ :  Compute $\zeta_{a,b}(s)$ for $b = 0$, then use a Taylor expansion to compute $\zeta_{a,b}(s)$ for all $b$.

Let's first compute $\zeta_{a,b}(s,n)$   for $b = 0$ : 
\begin {prop}
 If n is odd, then 
\begin{equation} \zeta_{a,0}(s,n) = 0\end{equation}
If n is even, then 
\begin{equation} \zeta_{a,0}(s,n) = \abs a^{-s}c_{-\frac n2}(a) (-i)^{\abs {\frac n2}}\pi^{1-s} \frac { \Gamma( \frac s2 + \frac {\abs n}4)}{ \Gamma( (1-\frac s2) + \frac {\abs n}4)} \end{equation}

\end {prop}
Remark : we then see that for $n = 0$, $\zeta_{a,0}(s,0)$ cancels for even positive integer values of $s$, so that the zeroes of $\zeta_{a,0}$ are not all on the line $\Re(s) = \inv 2$.

Proof : The function $\psi(z^2)$ is even, and the function $c_n$ is odd if n is odd, which proves the first formula, because the multiplicative convolution of any function with an even function gives an even function. 
If $n$ is even, we use the fact that $\CC$ is quadratically closed, i.e. that the map $z \mapsto z^2$ is onto.
A consequence of this property is that for any Schwartz function $\phi$ on $\CC$ and n even, we have 
\begin{equation} \Mell ( \phi(x^2),s,n) = \int_{\CC^*} \phi(x^2) \abs x^s ( \frac x{\abs x})^n d^\times x = 2 \int_{\CC^*} \phi(y) \abs y^\frac s2 ( \frac y{\abs y})^{\frac n2} \frac {d^\times y}4  = \inv 2 \Mell( \phi, \frac s2, \frac n2) \end{equation}
and this equality can easily be generalized to weak Mellin transforms. 
We observe, however, that for $n$ even and $0<\Re(s)<1$, the weak Mellin transform of the function $\psi_\CC(z)$ is well defined and equal to the the function $\rho(c_n \abs {}  ^s)$ appearing in Tate's local functional equation  ( cf \cite{Tate}, p 319) : 

\begin{equation}\Mell( \psi_\CC(z), s,n) =  \rho(c_n \abs {} ^s) = (-i)^{\abs n} \frac {(2\pi)^{1-s} \Gamma( s + \frac {\abs n}2)}{(2\pi)^s \Gamma( (1-s) + \frac {\abs n}2)}  \end{equation}

Indeed, Let's consider a function $\phi$  in $C_c^\infty (\CC^*)$. We want to compute the Mellin transform of 
\begin{equation} \lambda(\phi) \psi_\CC (z) = \int_{\CC^*} \phi(x) \psi_\CC(x^{-1}z) d^\times x \end{equation}
We can, however describe this integral as 
\begin{equation} = \int_{\CC} \inv {\abs y_\CC} \phi(\inv y ) \psi(y z) dy  = \four ( \phi^*)(z) \end{equation}
where the function $\phi^*$ is defined to be $\inv {\abs y_\CC} \phi(\inv y )$ and is a Schwartz function on $\CC$.
However we know by Tate's thesis that 
\begin{equation} \Mell( \phi^*,s,n) = \rho(c_n\abs {}^s)\Mell ( \four(\phi^*),1-s,-n) \end{equation}
\begin{equation} = \rho(c_n\abs {}^s)\Mell(\lambda(\phi)\psi ,1-s,1-n) \end{equation}
We then see that the weak Mellin transform of $\psi$ is well defined and that we have 
\begin{equation} \Mell( \psi, 1-s, -n) = \inv { \rho(c_n \abs {}^s)} =  \rho(c_{-n} \abs {}^{1-s}) \end{equation}

which shows that $\rho(c_n \abs {}^s )$ is indeed the weak Mellin transform of $\psi(z)$. It is then immediate that the weak Mellin transform of $\psi(az)$ is equal for $n$ even to 
\begin{equation} \Mell ( \psi(az),s,n) = \abs a_\CC^{-s}c_{-n}(a) \rho(c_n \abs {}^s )\end{equation}
so that we get 
\begin{equation} \Mell ( \psi(\frac a2z^2),s,n) = \inv 2 \Mell ( \psi(\frac a2 z),\frac s2,\frac n2) = \inv 2  \abs {\frac a2}_\CC^{-\frac s2}c_{-\frac n2}(a) \rho(c_{\frac n2} \abs {}^\frac s2 )\end{equation}
 which proves the formula for n even $\bb$
 
 Let's now consider the case $b \neq 0$.

\begin {prop}
The function $\zeta_{a,b}(s,n)$ satisfy the equations 
\begin{equation} \frac {\partial \zeta_{a,b}(s,n)}{\partial b} = -2\pi i \zeta_{a,b}(s+ \inv 2, n+1) \end{equation}
\begin{equation}  \frac {\partial \zeta_{a,b}(s,n)}{\partial \bar b} = -2\pi i \zeta_{a,b}(s+ \inv 2, n-1) \end{equation}
\end {prop}
Proof : these formula can be considered to be the Mellin transform of the formulae 
\begin{equation} \frac {\partial }{\partial b}\psi ( \frac a2 z^2+bz) = -2\pi i z \psi ( \frac a2 z^2+bz)\end{equation}
and 
\begin{equation}  \frac {\partial }{\partial \bar  b}\psi ( \frac a2 z^2+bz) = -2\pi i \bar z \psi ( \frac a2 z^2+bz)\end{equation}
the justification of the exchange of derivation and integration signs is done as in proposition \ref{exchange} $\bb$

Let's now give an explicit description of the weak Mellin transform of $\psi_\CC(\frac a2z^2+bz)$. Considering that we have  $\zeta_{a,b}(s) = \abs a^{-s} c_{-\frac n2}(a) \zeta_{1, \frac {b}{\sqrt a}}(s) $, we can suppose that $a = 1$
\begin {prop}
\label {zetacomplex}
The weak Mellin transform of $\psi_\CC(\frac 12z^2+bz) $ at the character $\abs z_\CC^s $ is equal to 
\begin{multline} \zeta_{1,b} (s) = \pi^{1-s} \lbrace \frac {\Gamma(\frac s2)}{\Gamma(1 - \frac s2)} {_1 F_1}(\frac s2, \inv 2, {i \pi b^2}){_1 F_1}(\frac s2, \inv 2, {i \pi \bar b^2})\\
 -  {4\pi} \abs b^2 \frac {\Gamma(\frac {s+1}2 )}{ \Gamma( 1 - \frac {s+1}2 )}  {_1 F_1}(\frac {s+1}2, \frac 32,  {\pi i b^2}){_1 F_1}(\frac {s+1}2, \frac 32, {\pi i \bar b^2})\rbrace \end{multline}
\end {prop}
Proof : 

We use the Wirtinger Taylor expansion
\begin{equation} \zeta_{1,b}(s) = \sum_{ n \ge 0, p \ge 0} \frac {\partial^{p+n}}{\partial^p b \partial^n \bar b } \zeta_{1,b}(s)(0)\frac {b^p\bar b^n}{p!n!} \end{equation}
\begin{equation}  = \sum_{ n \ge 0, p \ge 0}\frac {(-2\pi i b)^p(-2\pi i \bar b)^n}{p!n!} \zeta_{1,0}(s+\frac {p+n}2,p-n) \end{equation}

\begin{equation} = \pi^{1-s} \sum_{ n \ge 0, p \ge 0, \text {p+n even} }\frac {(2\pi i b)^p(2\pi i \bar b)^n}{p!n!} \frac { (-i)^{ \frac {\abs {p-n}}2}}  {  \pi^{ \frac {p+n}2}} \frac {\Gamma(\frac s2+ \frac {p+n}4 + \frac {\abs {p-n}}4)}{ \Gamma(1-\frac s2 - \frac {p+n}4 + \frac {\abs {p-n}}4)}\end{equation}
\begin{equation} =  \pi^{1-s}\sum_{ n \ge 0, p \ge 0, \text{p+n even} } { (-i)^{ \frac {\abs {p-n}}2}}   \frac {(2\sqrt {\pi} i b)^p(2\sqrt {\pi} i \bar b)^n}{p!n!} \frac {\Gamma(\frac s2+ \frac {\max(p,n)}2)}{ \Gamma(1-\frac s2 - \frac {\min(p,n)}2 )}\end{equation}
We can split the sum in two : the first sum $S_1$ is for n and p even, the second sum $S_2$ for n and p odd. 
\begin{equation} =  \pi^{1-s}( S_1+S_2) \end{equation}
The first sum $S_1$ becomes, writing $p = 2k$ and $n = 2l$ 
\begin{equation}  S_1 = \sum_{ k, l \ge 0 } { (-i)^{ {\abs {k-l}}}}   \frac {(2\sqrt {\pi} i b)^{2k}(2\sqrt {\pi} i \bar b)^{2l}}{(2k)!(2l)!} \frac {\Gamma(\frac s2+  {\max(k,l)})}{ \Gamma(1-\frac s2 - {\min(k,l)} )}\end{equation}
We now use the elementary formulas involving the pochhammer symbol : 
\begin{equation} \Gamma(\frac s2 +\max(k,l)) = (\frac s2+\max(k,l)-1)(\frac s2+\max(k,l)-2)...(\frac s2+1)\frac s2 \Gamma(\frac s2) = (\frac s2)_{\max(k,l)} \Gamma(\frac s2)  \end{equation}

and 

\begin{equation}{\Gamma(1-\frac s2)}=    (1 - \frac s2-1)(1-\frac s2-2)..(1-\frac s2-\min(k,l)) {\Gamma( 1 - \frac s2 - \min(k,l))}\end{equation}
\begin{equation} =  ( - \frac s2)(-\frac s2-1)..(-\frac s2-\min(k,l)+1) {\Gamma( 1 - \frac s2 - \min(k,l))}\end{equation}
\begin{equation} = (-1)^{\min(k,l)} (\frac s2)_{\min(k,l)} {\Gamma( 1 - \frac s2 - \min(k,l))}\end{equation}

We also have the  identities $(\frac s2)_{\max(k,l)} (\frac s2)_{\min(k,l)} = (\frac s2)_k (\frac s2)_l $
and $ (-i)^{\abs {k-l}}(-1)^{\min(k,l)} = (-i)^{\abs {k-l}}(-i)^{k+l -\abs {k-l}} = (-i)^{k+l}$
so that we get 
\begin{equation} S_1 = \frac {\Gamma(\frac s2)}{\Gamma(1 - \frac s2)}\sum_{ k, l \ge 0 } { (-i)^{ k+l}}  \frac {(2\sqrt {\pi} i b)^{2k}(2\sqrt {\pi} i \bar b)^{2l}}{(2k)!(2l)!} (\frac s2)_k ( \frac s2)_l  \end{equation} 

This expression  can be factored

\begin{equation} = \frac {\Gamma(\frac s2)}{\Gamma(1 - \frac s2)}(\sum_{k \ge 0} \frac {(-i)^k (2\sqrt{ \pi} ib)^{2k}}{ (2k)!}( \frac s2)_k)(\sum_{l \ge 0} \frac {(-i)^l (2\sqrt {\pi} i\bar b)^{2l}}{ (2l)!}(\frac s2)_l)\end{equation}
\begin{equation} = \frac {\Gamma(\frac s2)}{\Gamma(1 - \frac s2)}(\sum_{k \ge 0}(\frac  {(4 i \pi b^2)^k} {(2k)!} ( \frac s2)_k )(\sum_{l \ge 0}( \frac {(4i \pi \bar b^2)^l}  {(2l)!} ( \frac s2)_l )\end{equation}

and we recognize the product of two confluent hypergeometric functions, writing again $(2k)! = 4^k k!(\inv 2)_k$ : 
\begin{equation} = \frac {\Gamma(\frac s2)}{\Gamma(1 - \frac s2)}(\sum_{k \ge 0}(  ({i \pi b^2})^k \inv {(k)!} \inv {(\inv 2)_k} ( \frac s2)_k )(\sum_{l \ge 0}( ({i \pi \bar b^2})^l \inv {(l)!}\inv {(\inv 2)_l} ( \frac s2)_l )\end{equation}
\begin{equation} =  \frac {\Gamma(\frac s2)}{\Gamma(1 - \frac s2)} {_1 F_1}(\frac s2, \inv 2,  {i \pi b^2}){_1 F_1}(\frac s2, \inv 2, {i \pi \bar b^2})\end{equation}
The computations for the second sum $S_2$ are similar $\bb$

It should be noted that our proof of proposition \ref{zetacomplex} is not complete, considering that we have not proved that $\zeta_{a,b}(s)$ is analytic as a function of $\Re(b)$ and $\Im(b)$ ( The method used in the real case or for the case $\psi(\frac a2 \abs z^2 +bz)$ is not valid here). Let's sketch how the proof can be completed in this case : we remark that, as a function of $b$, $\zeta_{a,b}(s)$ is defined up to a scalar factor as  an eigendistribution for the action of a commutative subgroup of the metaplectic group ( cf remark following proposition \ref{fourier}). If we look at  infinitesimal generators of this subgroup, we see that $\zeta_{a,b}(s)$ can be defined up to a scalar factor as an eigenvector for the action of these infinitesimal generators, i.e. as a solution of some set of partial differential equations in the variables $\Re(b)$ and $\Im(b)$. In order to get a totally complete proof, it is then enough to write explicitly these generators and check that the function given in proposition \ref{zetacomplex} is indeed a solution of these partial differential equations.

\section {The weak Mellin transform of second degree characters defined on adele rings}
We consider a number field F and the associated adele ring $\AA_F$.
\subsection { Factorizable second degree characters on $\AA_F$}

 We know that the continuous characters of $\AA_F$ are of the form $\psi(bx)$ with $b \in \AA_F$, so that we can write any second degree character on $\AA_F$ as 
\begin{equation} f(x) = \psi( \inv 2 \alpha(x) x +bx) \end{equation}
with $b \in \AA_F$ and $\alpha$ is continuous morphism of additive group from $\AA_F$ to $\AA_F$ so that $\alpha^{-1}$ is also continuous.

We say that a second degree character is factorizable if it can be written as a tensor product ( noting $P_F$ the set of places of F)
\begin{equation} f(x) = \otimes_{v \in P_F} f_v \end{equation}
so that if an element of $\AA_F$ is written as $a = (a_v)$, we have $f(a) = \prod_{v \in P_F} f_v(a_v)$. 
For example, the second degree character $\psi(\frac a2x^2+bx)$ is factorizable, but if $\sigma$ is a Galois automorphism, and if we keep the notation $\sigma$ for its natural action on $\AA_F$, the second degree character $\psi(\frac a2\sigma(x)x+bx)$ is not in general factorizable.

If  $f$ is factorizable,we then also have $\alpha = \otimes_{v \in P_F} \alpha_v$, and the continuity of $\alpha$ and $\alpha^{-1}$ means that there exists a finite set S of valuations so that if $v$ is not in S, we have $\alpha_v(\cO^\times_v) = \cO^\times_v$, so that $\abs {\alpha_v} = 1$. 
It is clear that on $\AA_\QQ$, all second degree characters are factorizable  ( because there is no continuous non trivial additive continuous map from $\QQ_p$ to $\QQ_{p'}$ with $p \neq p'$).

\subsection { The existence of the weak Mellin transform}
On an adele ring, the weak Mellin transform is defined as follows : 
\begin {definition}
We say that a function $f$ defined on $\AA^\times$ has a well defined Mellin transform at the character $\abs x^s \chi(x)$ if exist a function $M_f(s,\chi)$ so that for any test function $\phi \in C_c^\infty (\AA_F^\times)$, we have 
\begin{equation} \Mell( \phi \star f , s,\chi) = \Mell( \phi,s,\chi) M_f(s,\chi)\end{equation}
\end {definition} 

Let's first prove the existence of the weak Mellin transform of a non degenerate second degree character defined on an adele ring.
We consider  a second degree character $f$, a function $\phi$ in $C_c^\infty(\AA^\times_F)$   and the map 
\begin{equation} \lambda(\phi)f(y) = \int_{\AA_F^\times} \phi(x) f(x^{-1} y)  d^\times x \end{equation}
where $y$ is an element of $\AA_F$

We then have the following proposition : 
\begin {prop}
If f is factorizable and $\phi \in C_c^\infty(\AA^\times)$, then $\lambda(\phi)f$ is a Schwartz function on $\AA$  
\end {prop}
Proof :  Considering that $\phi$ is locally constant as a variable of $x_v$ for all the finite places and that the support of $\phi$ is compact, we can write $\phi$ as a finite sum of functions of the form $\phi_\infty \otimes_{ \text {\text{finite} $v$}} \phi_v$. 
where $\phi_\infty \in C_0^\infty( F_{v_1} \times F_{v_2} \times..\times F_{v_n})$ (product of all the archimedean places)  and $\phi_v$ are functions defined on each local fields $F_v$ associated to finite valuations, with nearly all the $\phi_v$ equal to $\cf_{\cO^\times_v}$.
Let's then prove the proposition for such a function.
It is immediate that we have
\begin{equation} \lambda(\phi)(f) = \lambda (\phi_\infty) f _\infty \otimes_{ \text{finite } v } \lambda(\phi_v)f_v\end{equation}
 We have already proved that all the function $\lambda(\phi_v)f_v$ are in $S(F_v)$ for $v$ finite.
The two remaining points to prove is  that  $\lambda (\phi_\infty) f _\infty $ is Schwartz and that nearly all the functions $\lambda(\phi_v)f_v$ are equal to $\cf_{\cO_v}$. 
Let's first prove that $\lambda (\phi_\infty) f _\infty $ is Schwartz. Let's for example suppose  that we have two real places to consider, so that the seconde degree character $f_\infty$ can be decomposed as $f_1 \otimes f_2$ where $f_1$ and $f_2$ are non degenerate seconde degree characters defined over $\RR$
We then have to consider the function 
\begin{equation} \lambda (\phi_\infty) f _\infty (y_1,y_2) = \int_{\RR^{*2}} \phi_\infty(x_1,x_2) f_1(x_1^{-1}y_1)f_2(x_2^{-1}y_2) d^\times x_1 d^\times x_2 \end{equation}
Let's write $\phi_\infty^*(x_1,x_2) = \inv { \abs {x_1x_2}}\phi_\infty(\inv {x_1},\inv {x_2})$. We extend this function to $\RR^2$ by writing $\phi_\infty^*(0,x_2) = \phi_\infty^*(x_1,0) = 0$, so that $\phi^*_\infty$ is a Schwartz function an the integral becomes
\begin{equation} \lambda (\phi_\infty) f _\infty (y_1,y_2) =  \int_{\RR^{2}} \phi_\infty^*(z_1,z_2) f_1(z_1y_1)f_2(z_2y_2) d z_1 d z_2 \end{equation}
Let's for example show that this function is fast decreasing as a function of $y_1$ : we apply the local Weil functional equation to the integral on $z_1$

\begin{equation} = \frac {\gamma_{f_1}}{ \abs {y_1} \sqrt {\abs {a_1}}}\int_{\RR^{2}} \four_{z_1}(\phi_\infty^*)(z_1,z_2) f_1(\frac {z_1}{a_1y_1})f_2(y_2z_2) d z_1 d z_2 \end{equation}
considering that $\phi^*$ and all its derivatives cancel at $(0,z_2)$ for any value of $z_2$, this is equal, for any polynomial P, to 
\begin{equation} = \frac {\gamma_{f_1}}{ \abs {y_1} \sqrt {\abs {a_1}}}\int_{\RR^{2}} \four_{z_1}(\phi_\infty^*)(z_1,z_2)( f_1(\frac {z_1}{a_1y_1})- P( \frac {z_1}{a_1y_1}))f_2(y_2z_2) d z_1 d z_2 \end{equation}
We then choose P to be the Taylor expansion of $f_1( x)$, and the proof can be finished as in the real case.
In order to show that all the partial derivatives are fast decreasing, we use again the fact that the operator $x_1 \frac {\partial}{\partial x_1}$ 
satisfy 
\begin{equation} x_1 \frac {\partial}{\partial x_1} \lambda(\phi) f =   \lambda(x_1\frac {\partial}{\partial x_1}\phi) f \end{equation}
The proof is the same for a function $\phi_\infty$ defined on $\RR^n \times \CC^m$.\\

Lets' now show that nearly all the functions $\lambda(\phi_v)f$ are equal to $\cf_{\cO_v}$.
We remark that 
\begin {itemize}
\item 
 nearly all the $\phi_v$ are equal to $\cf_{\cO_v^\times}$

  \item nearly all the $\alpha_v$ satisfy $\abs {\alpha_v} = 1$ ( because $\alpha$ is a continuous and its inverse is also continuous)
 
 \item nearly all the $b_v$ satisfy $\abs b_v \le 1$ ( $b$ is an adele)
 \item the local different $\xd_v$ is equal to $\cO_v$ for nearly all valuations \end {itemize}
 
It is then  enough to prove that for the valuations $v \in P_F$ satisfying $\abs {\alpha_v}=1$  ,$\abs {b_v} \le 1$, and $\xd_v = \cO_v$, the functions $\lambda(\cf_{\cO_v^\times})(f_v)$ is equal to $\cf_{\cO_v}$. We will also suppose that $\abs 2_v = 1$, which is also true for nearly all valuations. 
The computation is then exactly the same as in proposition $\ref{qp}$ for $p \neq 2$ : 
\begin{equation} \lambda(\cf_{\cO_v^\times})f_v(x) = \int_{F_v^\times} f_v(\frac yx) \cf_{\cO_v\times}(x) d^\times x \end{equation}
\begin{equation} = \int_{F_v^\times} \psi(\inv 2 \alpha(\frac yx)\frac yx + b \frac yx) \cf_{\cO_v^\times}(x) d^\times x \end{equation}
If $y \in \cO_v$, and $x \in \cO_v^\times$, then  $\frac yx \in \cO_v$ so that $\alpha(\frac yx) \in \cO_v$ and  the expression $ \psi(\inv 2 \alpha(\frac yx)\frac yx + b \frac yx)$ is equal ro 1 because $\psi$ is trivial on $\cO_v$. The integral is then equal to the (multiplicative) measure of $\cO_v^\times$ which is  $(\cN \xd_v)^{-\inv 2} = 1$  
Let's now suppose that the valuation of $y$ is negative strict. We replace the integral on $\cO_v^\times$ by an additive integral on $\cO_v$ 
\begin{equation} =  \inv {1 - \inv {\cN \xp}} \int_{F_v} \psi(\inv 2 \alpha(\frac yx)\frac yx + b \frac yx) \cf_{\cO_v^\times}(x) \frac {d x}{\abs x} \end{equation}
We remark that the term $\abs x$ can be removed since the integral is on $\cO_v^\times$ and use the Weil local functional equation :
\begin{equation}  \inv {1 - \inv {\cN \xp}}\frac {\gamma_f}{\sqrt {\abs \alpha}} \int_{F_v} \psi(-\inv 2 \frac yx\alpha^{-1}(\frac yx) - b\alpha^{-1}( \frac yx)) \four(\cf_{\cO_v^\times}(x))  {d x}\end{equation}
Let's note $\varpi$ an element of $F_v$ with valuation 1. We can then write that 
\begin{equation} \cf_{\cO_v^\times} = \cf_{\cO_v} - \cf_{\varpi \cO_v}\end{equation}
Using $\xd_v = \cO_v$, we see that $\four( \cf_{\cO_v}) = \cf_{\cO_v}$ so that 
\begin{equation} \four ( \cf_{\cO_v^\times}) = \cf_{\cO_v} - \inv {\cN \xp} \cf_{\inv \varpi \cO_v} \end{equation}
The support of $\four(\cf_{\cO_v^\times})$ is then included in $\inv {\varpi} \cO_v$ so that we can restrict the integral to the set of elements $x$ having a valuation $\ge -1$. We remark, however, that if the valuation of $y$ is $<0$ and the valuation of $x$ is $\ge -1$, then the valuation of $\frac yx$ is $\ge 0$, so that the function $ \psi(-\inv 2 \frac yx\alpha^{-1}(\frac yx) - b\alpha^{-1}( \frac yx))$ is equal to 1. The integral is then equal to 
\begin{equation}\inv {1 - \inv {\cN \xp}}\frac {\gamma_f}{\sqrt {\abs \alpha}} \int_{F_v} \four(\cf_{\cO_v^\times}(x))  {d x} = \inv {1 - \inv {\cN \xp}}\frac {\gamma_f}{\sqrt {\abs \alpha}} \cf_{\cO_v^\times}(0) = 0 \end{equation}
We have then proved that $\lambda(\cf_{\cO_v^\times})f_v(y)$ is equal to 1 if the valuation of $y$ is positive or zero, and zero if the valuation of $y$ is negative strict, so that $\lambda(\cf_{\cO_v^\times})f_v = \cf_{\cO_v}$

$\bb$
\begin {prop}
If f is a factorizable non degenerate second degree character defined on $\AA_F$, then the weak Mellin transform of f is well defined for $\Re(s) >1$.
\end {prop}
Proof : the proof is  the same as for the local case ( replace $>0$ by $>1$).

We will  use the notation $\Xi_{f}(s,\chi)$ for the weak Mellin transform of a  second degree character $f$ defined on an adele ring  at the character $\abs x^s \chi(x)$

\subsection {The functional equation of $\Xi_{f}$}

\begin {prop}

for $\phi \in S(\AA^\times)$ and f a non degenerate second character defined on $\AA$, The Fourier transform of the Schwartz function $\lambda(\phi)f (y)$ is equal to $\frac { \gamma_{f}}{\sqrt {\abs \alpha}} \lambda (\phi^*) \bar f \circ \alpha^{-1}(  y)$, where $\phi^*$ is defined by the formula $\phi^*(x) = \inv {\abs x} \phi(\inv x)$

\end {prop}
Proof : the proof is exactly the same as for the local case $\bb$

\begin {prop}
Let's assume that $\chi$ is a Hecke character on $\AA_F^\times$. Then the function $\Xi_f(s,\chi)$ considered as a function of s has an analytic continuation to $\CC$, with possible poles at 0 and 1  if $\chi $ is unramified. If we keep the notation $\Xi_f(s,\chi)$ for the analytic continuation, we have the equality 
\begin{equation} \Xi_f(s,\chi) = \frac {\gamma_f}{\sqrt {\abs \alpha}} \Xi_{\bar f \circ \alpha^{-1}}(1-s, \chi) \end{equation}
\end {prop}
Proof : The proof is the same as for the local case, but we have to replace the Tate local functional equation by the global functional equation $\bb$

\subsection { The connection with Hecke L-functions}

We remind that if $\chi = \otimes \chi_v$ is a Hecke character defined on $\AA_F^\times$, the Hecke L-function $L(s,\chi)$ is defined for $\Re(s) >1$ as 
\begin{equation} L(s,\chi) = \prod_{\text{$v$ finite},  \chi_v \text{unramified at }v} \inv {1 - \chi_v (\varpi_v) (\cN \xp_v)^{-s}} \end{equation}

\begin {theorem}
Let's consider a number field F, a unitary Hecke character $\chi$ and note $L(\chi,s)$ the associated Hecke L-function. 
Let's consider a factorizable non degenerate second degree character $f$and note $\Xi_f(s,\chi)$ the Mellin transform of f at $(s,\chi)$. Then $(s,\chi)$ is a zero of $\Xi_f$ if and only if it is a non trivial zero of $L(s,\chi)$ or a zero of one of the local functions $\zeta_{f_v}$

\end {theorem}

Proof : 
 
We have already computed that for nearly all valuations, we have 
\begin{equation} \lambda( \cf_{\cO_v^\times})f_v = \cf_{\cO_v} \end{equation}
so that 
\begin{equation} \zeta_{f_v}(s,\chi) \Mell( \cf_{\cO_v}^\times, s, \chi) = \Mell (\cf_{\cO_v},s,\chi) \end{equation}
Nearly all these valuations satisfy the condition that $\chi_v$ is unramified at $v$. \\
Let's then consider the two following finite sets : The set $S$ is the set of all valuations which  are either infinite, or finite with $\chi_v$ ramified, or satisfy $\lambda( \cf_{\cO_v^\times})f_v \neq \cf_{\cO_v}$ or satisfy $\xd_v \neq \cO_v$. The set T is the set of all valuations which are eiher infinite or finite with $\chi_v$ ramified. It is clear that $T \subset S$. 

Let's suppose that $v$ is not in $S$. Since $\chi$ is unramified and  $\cN \xd_v = 1$, we have $\Mell(\cf_{\cO_v}^\times,s,\chi) = 1$ so that we get 
\begin{equation} \zeta_{f_v}(s,\chi)= \Mell (\cf_{\cO_v},s,\chi)  = \int_{F_v^\times} \cf_{\cO_v}(x) \chi(x) \abs x^s d^\times x \end{equation}
we split the integral according the the valuation of x and get 
\begin{equation} 1 + \chi(\varpi_v)\cN \xp^{-s} + ..\end{equation}
\begin{equation} = \inv {1 - \chi(\varpi_v)\cN \xp^{-s}} \end{equation}
We can then write
\begin{equation} \Xi_f(s,\chi) = \prod_{v \in S} \zeta_{f_v}(s,\chi_v) \prod_{ v \notin S} \inv {1 - \chi(\varpi_v)\cN \xp^{-s}}\end{equation}
let's compare this with 
\begin{equation} L(s,\chi) = \prod_{ v \notin T} \inv {1 - \chi (\varpi_v) (\cN \xp_v)^{-s}} \end{equation}

 We can split $L(s,\chi)$ as 
\begin{equation} L(s,\chi) = \prod_{ v \in S, v \notin T }  \inv {1 - \chi(\varpi_v)\cN \xp^{-s}} \prod_{  v \notin S, v \notin T} \inv {1 - \chi(\varpi_v)\cN \xp^{-s}} \end{equation}
and $\Xi_f$ as 
\begin{equation} \Xi_f(s,\chi) = \prod_{v \in S, v \notin T} \zeta_{f_v}(s,\chi_v) \prod_{v \in S, v \in T} \zeta_{f_v}(s,\chi_v)\prod_{ v \notin S, v \notin T } \inv {1 - \chi(\varpi_v)\cN \xp^{-s}}\end{equation} 
which leads to the equality
\begin{equation}  \Xi_f(s,\chi) = \prod_{v \in S, v \notin T}(1 - \chi(\varpi_v)\cN \xp^{-s}) \zeta_{f_v}(s,\chi_v)( \prod_{v \in S, v \in T} \zeta_{f_v}(s,\chi_v)) L(s,\chi) \end{equation}
Note that both products are finite, because S is finite.

We know that the local functions $\zeta_{f_v}(s,\chi_v)$ do not have any pole for $\Re(s)>0$ and that the only pole of $L(s,\chi)$ is at $\chi(x)\abs x^s = \abs x$. It is then immediate, since S and T are finite sets, that a zero of $L(s,\chi)$ satisfying $\Re(\rho)>0$  is also a zero of $\Xi_f(s,\chi)$ an that any zero of $\zeta_v(s,\chi_v)$  is also a zero of $\Xi_f$ (we have seen that 1 is never a zero of $\zeta_f$ for unramified characters). Conversely, it is immediate that a zero of $\Xi_f$ has to be either a zero of some $\zeta_{f_v}$for $v \in S$ or a zero of the product $\prod_{v \in S,   v \notin T} (1 - \chi(\varpi_v)\cN \xp^{-s})$ or a zero of $L(s,\chi)$. We remark that the zeroes of $\prod_{ v \in S, v \notin T} (1 - \chi(\varpi_v)\cN \xp^{-s})$ are not zeros of $\Xi_f$ because they are canceled by the same terms appearing in the definition of $L(s,\chi)$. The non trivial zeroes of $L(s,\chi)$ at negative integers are also not zeroes of $\Xi_f(s,\chi)$ : if this were the case, we would also find these zeroes at positive integers thanks to the functional equation, and we know that it is not possible $\bb$.

\section {Weak Mellin transforms and second degree characters defined on vector spaces }

Let's now  come back to Riemann's proof of the functional equation of $\zeta$, which is based on the Poisson summation formula on $\ZZ$, i.e. to the fact that the distribution $\delta_{\ZZ}$ is equal to its Fourier transform. We know that this Poisson summation formula can be generalized to distributions of the form $\delta_{\ZZ^n}$ defined on $\RR^n$, leading to the functional equation of Epstein Zeta functions or Eisenstein series. The idea of this section is to perform a similar generalization, replacing a second degree character defined on a field by a second degree character defined on a vector space. The main result of this section is that the natural generalization of the local functions $\zeta_{f_v}$ to vector spaces have their zeroes on the line $\Re(s) = \frac n2$ under reasonable conditions.

\subsection { A local functional equation on vector spaces}
In order to generalize our results to second degree characters defined on vector space, we first need a generalization of Tate's local functional equation to Schwartz functions defined on vector spaces.
We define the following maximal compact subgroups $K_L$ of $GL_n(L)$ : for  $L = \RR$, we write $K_\RR = SO(n)$. For $L = \CC$, we write $K_\CC = U(n)$. If L is a local field,  we write $K_{L} = GL_n(\cO_L)$
We define the norm of an element of the vector space $L^n$ as follows : if L is a local field, $\norm x = \sup ( \abs {x_1}, \abs {x_2}, , \abs {x_n})$. If L is equal to $\RR$,  $\norm x$ is the usual norm. If L is equal to $\CC$, we have to take the square of the usual norm, i.e. $\norm x_\CC = \norm x^2$. Note that $\norm x$ is always invariant under the action of the compact group $K_L$

Let's  consider for any Schwartz function $\varphi$ defined on $L^n$ and $\Re(s) >0$  the integral 
\begin{equation}M( \varphi,s)  = \int_{L^n} \varphi(x) \norm x^s  \frac {dx}{{\norm x}^n} \end{equation}
This integral is well defined : The convergence near zero is a consequence of the fact that $\varphi$ is continuous. The convergence for $\norm x$ large is a consequence of the fact that $\varphi$ is Schwartz. We call this integral the Mellin transform of $\varphi$.

If a function $\varphi$ defined on $L^n$, is invariant under the action of $K$, we say that it is a radial function.

\begin {prop} 
\label {LFT vector space}
For any Schwartz function f ,  and $s \in \CC$ with $0<\Re(s)<n$, the  Mellin transform of f and $\four(f)$ are related by the following formula
\begin{equation} M(f,s) = \rho_n(s) M(\four(f),n-s)\end{equation}
for some scalar $\rho_n(s)$ which does not depend on f 
\end {prop}

Proof :  
Let's first suppose that $L = \RR$. Then this proposition simply states that the Fourier transform of $\norm x^{s-n}$ considered as a distribution is equal to $\norm x^{-s}$ up to a scalar factor. This is a well known result in the theory of homogeneous distributions ( cf \cite {Grafakos}, Th 2.4.6). Considering that a radial homogeneous distribution on $\CC^n$ can also be considered as a radial homogeneous disribution on $\RR^{2n}$, the result is then also true for $L = \CC$.

Let's now consider the case where $L$ is a local field, $L = \QQ_{\xp}$.
Let's first consider the Schwartz functions which are radial (i.e. invariant under the action of $K$). Such a function $\phi$ can be written as a finite sum of the form $\phi = \sum_k a_k\cf_{\varpi^{k}\cO^n}$ ( because $\phi$ has compact support and is continuous in zero), and its Fourier transform is equal to $\four(\phi) = (\cN \xd)^{-\frac n 2} \sum_k a_k \inv {q^{nk}}\cf_{\varpi^{-k}\xd^{-1}\cO^n}$, and it is immediate, writing $\xd = \xp^d$ and $\abs \varpi = q$  that we have 
\begin{equation} M(\phi,s) = q^{d(\frac n2-s)}\frac {1-q^{s-n}}{1 - q^{-s}}M(\four(\phi), n-s) \end{equation}
Let's now suppose that $\phi$ is not radial. We remark that if $\phi_K$ is the radial function obtained by averaging $\phi$ under the action of K, we have $M(\phi,s) = M(\phi_K,s)$. Considering that this averaging action commutes with the Fourier transform, we get the result for all $\phi$ $\bb$

\subsection { From second degree characters on vector spaces to Schwartz functions}

We also need a generalization of proposition \ref{regul} for vector spaces. This is given by the following proposition : 
\begin {theorem}
\label{regulvector}
Let's consider a non degenerate second degree character f on a L-vector space $L^n$ of finite dimension n, where L is a locally compact field. Let's consider a function $\phi$ in $C_c^\infty(GL_n(L))$ and define $\lambda(\phi)(f)$ as
\begin{equation} \lambda(\phi)(f)(v) = \int_{GL_n(L)}\phi(x) f(x^{-1}v) d^\times x\end{equation}
Then $\lambda(\phi)(f)$ is a Schwartz function on $L^n$
\end {theorem}
Proof : Let's first prove an elementary proposition

\begin {prop}
\label {transitive}
Any element $v \in L^n - \{0\}$ can be written as $kv_1$ where k is in K and the only non zero coordinate of $v_1$ is the first one. (i.e $v_1 = (x_1,0,0..)$ for some $x_1$ in L)

\end {prop}
Proof : 
This is clear for $\RR^n$ and $\CC^n$. Let's then suppose that L is local and note $\varpi$ an uniformizer, with $\abs \varpi = \inv q$.It is enough to show that any vector  $x$ in $L^n$ satisfying $\norm x = q^{-m}$ is in the orbit of the  vector $v_1 = (\varpi^{m} ,0,0,0..)$  under the action of $K_L$ for any $m \in \ZZ$. After multiplication by $\varpi^{-m}Id$, which commutes with K, we can suppose that $m = 0$.

Considering that $\norm x = \max {\abs x_i} = 1$, all the coordinates of $x$ are in $\cO_L$ and at least one of the coordinates of $x$ is a unit. 
Since the map exchanging the basis vectors $e_1$ and $e_i$ is in K, for all i, we can suppose that the first coordinate $x_1$ is a unit. We can then write that 
\begin{equation} \begin {pmatrix}  x_1 \\ x_2 \\x_n  \end {pmatrix} =   \begin {pmatrix}  x_1 & 0 & 0 .. \\ x_2 & 1 & 0.. \\x_3 & 0  & 1 ..  \end {pmatrix} \begin {pmatrix}  1 & 0 & 0   \end {pmatrix}  \end{equation}
and the square matrix is in $GL_n(\cO_L)$ $\bb$ 

Let's now prove the theorem.
Suppose that L is local. It is immediate that $ \lambda(\phi) f  $ is locally constant, so that we have to show that it has a compact support in V. 
In order to simplify notations, we suppose that $n = 2$, but the proof remains the same for all $n$. 
We first convert  the integral on $GL_2(L)$
\begin{equation} \lambda(\phi)(f)(v) = \int_{GL_2(L)}\phi(x) f(x^{-1}v) d^\times x\end{equation}
into an integral on $M_2(L)$ : If $dx$ is the standard haar measure on $M_2(L)$,  $\frac {dx}{\abs {\det x}^2}$ is a haar measure on $GL_2(L)$ so that it is equal to the standard haar measure on $GL_2(L)$ up to a constant scalar factor.  $\lambda(\phi)(v)$ is then equal, up to a constant scalar factor, to the integral
\begin{equation}  \int_{M_2(L)}\phi(x) f(x^{-1}v) \frac {d x}{\abs {\det x}^2} \end{equation}
writing $ \phi^\times(x) = \phi(x^{-1})\inv  {\abs {\det x}^2} $  ( not to be confused with $\phi^* = \phi(x^{-1})\inv  {\abs {\det x}} )$ for $x \in GL_2(L)$ and $\phi^\times(x) = 0$ for $x \notin GL_2(L)$, the integral becomes
\begin{equation}  \int_{M_2(L)}\phi^\times(x)f(xv)dx \end{equation}

 where $\phi^\times$ is a Schwartz function on $M_2(L)$. 

Let's first  suppose that the vector $v$ is of the form $(y,0)$
Let's write the matrix $x$ as $x =  \begin {pmatrix}  \alpha & \beta\\ \gamma & \delta  \end {pmatrix}$, so that $xv = (\alpha y, \gamma y) = y(\alpha,\gamma)$ and the integral can be described as 
\begin{equation} \int_{\alpha, \beta, \gamma, \delta \in \QQ_p} f(\alpha y, \gamma y ) \phi^\times( \begin {pmatrix}  \alpha & \beta\\ \gamma & \delta  \end {pmatrix}) d\alpha d\beta d\gamma d\delta \end{equation}
\begin{equation} = \int_{ \beta,  \delta \in \QQ_p}\lbrace \int_{\alpha, \gamma,  \in \QQ_p} f(\alpha y, \gamma y ) \phi^\times (  \begin {pmatrix}  \alpha & \beta\\ \gamma & \delta  \end {pmatrix}) d\alpha d \gamma \rbrace  d\beta d\delta \end{equation}
We have assumed that the  function $f(\alpha, \gamma)$, is a non degenerate second degree character on $L^2$. If we note $\varrho$ the morphism from $V$ to $V^*$ associated to $f$, and identify $V$ with $V^*$ using the standard additive character $\psi$,  then the weak fourier transform of $f(\alpha, \gamma)$ is equal to $\frac {\gamma_f} {\sqrt {\varrho}}\bar f( \varrho^{-1} (\alpha, \gamma))$ ( cf Weil $\cite {Weil64}$). The fourier transform of $f( y (\alpha, \beta))$ is then equal to $\inv {\abs y} \frac {\gamma_f} {\sqrt {\varrho}}\bar f(y^{-1} \varrho^{-1} (\alpha, \gamma))$
We then get the equality, with the notation $\four_{\alpha, \gamma}(\phi^\times)$ for the Fourier transform of $\phi^\times(\alpha, \beta, \gamma, \delta)$ considered as a function of $\alpha$ and $\gamma$ only. 
\begin{equation} = \inv {\abs y} \frac {\gamma_f} {\sqrt {\varrho}}\int_{ \beta,  \delta \in \QQ_p}\lbrace \int_{\alpha, \gamma,  \in \QQ_p} \bar f(-\varrho^{-1}(\frac {\alpha}y , \frac {\gamma} y ) \four_{\alpha, \gamma} (\phi^\times)( \alpha,\beta, \gamma, \delta) d\alpha d \gamma \rbrace  d\beta d\delta \end{equation}
We then remark that the function $\four_{\alpha, \gamma}(\phi^\times) ( \alpha, \beta, \gamma, \delta)$ has compact support ( it is a Schwartz function on $L^4$).As a consequence
we can  suppose that its support is included in a ball of radius R ( using the sup norm $\norm {\alpha,\beta,\gamma,\delta} = \max ( \abs \alpha, \abs \beta, \abs \gamma, \abs \delta$) ) . We also know that $f\circ \varrho^{-1}$ is continuous and equal to 1 near zero, so that there exists some $\epsilon$ so that if $\abs {\alpha} < \epsilon$ and $\abs {\gamma} < \epsilon$, then $\bar f \circ \varrho^{-1}(\alpha, \gamma) = 1$. It is then immediate that if $\abs y > \frac R \epsilon$, then the integral becomes zero :  The expressions $\abs {\frac \alpha y }$ and $\abs {\frac \gamma y}$ are always  lower than $\epsilon$ if $\alpha$ and $\gamma$ are in the support of $\four_{\alpha,\gamma}(\phi^\times)$ so that the integral becomes, considering that the matrix $  \begin {pmatrix}  0 & \beta\\ 0 & \delta  \end {pmatrix}$ is not in $GL_2(L)$  : 
\begin{equation} \inv {\abs y} \frac {\gamma_f} {\sqrt {\varrho}}\int_{ \beta,  \delta \in \QQ_p}\lbrace \int_{\alpha, \gamma,  \in \QQ_p}  \four_{\alpha, \gamma} (\phi^\times)( \alpha,\beta, \gamma, \delta) d\alpha d \gamma \rbrace  d\beta d\delta \end{equation}
\begin{equation}  \inv {\abs y} \frac {\gamma_f} {\sqrt {\varrho}}\int_{ \beta,  \delta \in \QQ_p} \phi^\times( \begin {pmatrix}  0 & \beta\\ 0& \delta  \end {pmatrix} )    d\beta d\delta  = 0\end{equation}
Let's now suppose that the vector $v$ is not of the form $(y,0)$. We have seen that it is always possible to write $v$ as $v = kv'$ with $k \in GL_2(\cO)$ and $v' = (y',0)$ for some $y' \in L$. Note that the sup norm $\norm {v}$ is equal to $\abs {y'}$.
It is immediate that we have 
\begin{equation} \lambda (\phi)(v) = \lambda(\phi) (kv')= \lambda(\phi(kx))(v') \end{equation} 
If the support of $\phi$ is included in a ball of radius R, then the support of $\phi(kx)$ is included in the same ball, since we have the equality of sup norms $\norm{kx} = \norm x$ for all $ k \in GL_2(\cO)$ and all $x$ in $GL_2(L)$. 
The function is then equal to zero if $\abs {y'} = \norm {v} > \frac R{\epsilon} $, which proves the theorem for L local.

Let's now consider the case $L = \RR$ or $\CC$. 

Let's for example take $L = \RR$.
If $\phi \in C_c^\infty(GL_2(\RR))$, we write again $\phi^\times(x) = \phi(x^{-1}) \inv {\abs {\det(x)}^2}$.
 The support of $\phi^\times$ is also compact in $M_2(\RR)$ for the topology of $M_2(\RR)$( because the inclusion map from $GL_2(\RR)$ to $M_2(\RR)$ is continuous) so that if $x$ is any element of $M_2(\RR)$ which is not in $GL_2(\RR)$, then g is zero in a neighborhood of $x$ so that all its derivatives in $x$ cancel.
 
It is immediate that $\phi^\times$ is Schwartz on $M_2(\RR)$ because it is $C^\infty$ with compact support.
We first consider the case $v= (y,0)$ and the same computation as for the case L ultrametric leads to the integral 
\begin{equation} \inv {\abs y} \frac {\gamma_f} {\sqrt {\varrho}}\int_{ \beta,  \delta \in \QQ_p}\lbrace \int_{\alpha, \gamma,  \in \QQ_p} \bar f(\varrho^{-1}(\frac {\alpha}y , \frac {\gamma} y ) \four_{\alpha, \gamma} (\phi^\times)( \begin {pmatrix}  \alpha & \beta\\ \gamma & \delta  \end {pmatrix}) d\alpha d \gamma \rbrace  d\beta d\delta \end{equation}
We then observe, that for $n \ge 0, m \ge 0$, the integral 
\begin{equation} \int_{\alpha, \gamma,  \in \QQ_p}\alpha^n \gamma^m \four_{\alpha, \gamma} (\phi^\times)(  \begin {pmatrix}  \alpha & \beta\\ \gamma & \delta  \end {pmatrix}) d\alpha d \gamma  \end{equation}
is equal, up to a constant, to 
\begin{equation}  \int_{\alpha, \gamma,  \in \QQ_p} \four_{\alpha, \gamma} (\frac {\partial^{n+m}}{\partial \alpha^n \partial \gamma ^m} \phi^\times)(  \begin {pmatrix}  \alpha & \beta\\ \gamma & \delta  \end {pmatrix}) d\alpha d \gamma  = \frac {\partial^{n+m}}{\partial \alpha^n \partial \gamma ^m} \phi^\times  \begin {pmatrix}  0 & \beta\\ 0 & \delta  \end {pmatrix} = 0\end{equation}
If $P(\alpha,\beta)$ is the polynomial associated to the Taylor expansion of degree $n$ of $\bar f(\varrho^{-1}({\alpha},  {\gamma} ) $, the integral is then equal to 
\begin{equation}  \inv {\abs y} \frac {\gamma_f} {\sqrt {\varrho}}\int_{ \beta,  \delta \in \QQ_p}\lbrace \int_{\alpha, \gamma,  \in \QQ_p}( \bar f(\varrho^{-1}(\frac {\alpha}y , \frac {\gamma} y ))  - P(\frac \alpha y, \frac \gamma y))\four_{\alpha, \gamma} (\phi^\times)( \begin {pmatrix}  \alpha & \beta\\ \gamma & \delta  \end {pmatrix}) d\alpha d \gamma \rbrace  d\beta d\delta \end{equation}
and the remainder of the proof is similar to the one dimensional case. The proof for $L = \CC$ is similar $\bb$

 This theorem can be extended  without difficulty to vector spaces defined over locally compact division rings D, since the commutativity of the field has not been used in the proofs.

\subsection { The weal Mellin transform of a second degree character defined on a vector space}

On a locally compact field, we have defined the weak Mellin transform thanks to  the formula $\Mell( \lambda(\phi)f,s) = \Mell(\phi,s)\Mell(f,s)$ which is valid for all Schwartz functions with $\Re(s)>0$. This formula can also be written as  $\lambda(\phi) \abs x^{s-1} = \Mell(\phi,1-s) \abs x^{s-1} $ : the function  $\abs x^{s-1}$ on $\RR^*$ is stable, up to a scalar factor, for the action of the multiplicative group.
In order to generalize this formula to vector spaces defined on locally compact fields, we consider the function $\norm x^s$  on $L^n - \{0\}$ and the natural left action $\lambda$  of $GL_n(L)$ on this  function. If $\phi$ is a general element of $C_c^\infty(GL_n(L))$, then $\lambda(\phi) \norm x^s$ is not equal to $\norm x^s$ up to a scalar factor. We however show that this is  the case if the function $\phi$ is invariant under the action of K, which allows to define the weak radial Mellin transform.

\begin {prop}
Let's consider a function $\nu_s$ defined on $L^n - \{0\}$, invariant under the action of K  and satisfying $\nu_s( \lambda x) = \abs \lambda^s \nu_s(x)$ for all $\lambda \in L$, then $\nu_s$ is equal to the function $\norm x^s$ up to a scalar factor.
\label {radial}
\end {prop}
Proof :  Let's note $C$ the value of $\nu_s(x)$ on $e_1 = (1,0..)$. We then write, using a decomposition $x = kv_1$ given in proposition \ref {transitive} 
\begin{equation} \nu_s(x) = \nu_s(k (x_1,0,0..)) = \abs {x_1}^s \nu_s( 1,0,0..)  = C \norm x^s \bb \end{equation}

Let's note $\cH (GL_n(L))$ the Hecke algebra of $GL_n(L)$, i.e. the algebra of functions in $C_c^\infty(GL_n(L))$ invariant under the left and right action of K.
\begin {prop}
\label {definitionxi}
Let's consider a function $\phi \in \cH (GL_n(L))$ . Then 
For all $s \in \CC$,  there exists a scalar $\xi_s(\phi)$  so that  the equation
\begin{equation} \lambda(\phi) \norm x^{s} = \xi_s(\phi) \norm x^{s} \end{equation}
is valid for all $x \in L^n - \{0\}$
\label {prop:radial2}
\end {prop}
Proof : Let's note $f(x) = \norm x^s$
We have \begin{equation} \lambda(\phi)(f)(kv) = \int_{GL_n(L)}\phi(x) f(x^{-1}kv) d^\times x\end{equation}
writing $x^{-1}k = y^{-1}$, so that $ky = x$, we get 
\begin{equation} = \int_{GL_n(L)}\phi(ky) f(y^{-1}v) d^\times y\end{equation}
\begin{equation}= \lambda(\phi)(f)(v)\end{equation}

it is then  immediate that    $\lambda(\phi) \norm x^{s}$ satisfies the conditions of the previous proposition , so that it is equal to $\norm x^{s}$ up to a scalar factor on $\RR^n - \{0\}$ $\bb$\\

\begin {prop}
The function $\xi_s$ is a character of the Hecke algebra $ \cH (GL_n(L))$ : 
If  $\phi_1$ and $\phi_2$ are in $ \cH (GL_n(L))$,  we have 
\begin{equation} \xi_s(\phi_1 \star \phi_2) =  \xi_s(\phi_1) \xi_s( \phi_2) \end{equation}
\end {prop}
Proof : immediate consequence of the proposition $\lambda(\phi_1 \star \phi_2) = \lambda(\phi_1)\lambda(\phi_2)$
\begin {prop}
\label{transferxi}
Let's consider some function $\varphi$ defined on L so that $M(\varphi,s)$ is well defined. Then $M(\lambda(\phi)\varphi,s)$ is well defined and we have for  $\phi \in {\cH(GL_n(L))}$ the equality 
\begin{equation} M( \lambda(\phi)\varphi,s) = \xi_{s-n}(\phi^*)M(\varphi,s) \end{equation}
whith $\phi^*(g) = \inv {\abs {\det g}} \phi ( g^{-1})$
\end {prop}
Proof : 
We have 
\begin{equation} M( \lambda(\phi)\varphi,s) = \int_{ x \in L^n} \int_{g \in GL_n(L)} \phi(g) \varphi(g^{-1}x) \norm x^s d^\times g \frac {dx}{\norm x^n}   \end{equation}
the double integral is absolutely convergent, so that we can exchange the order of summations, and write $y = g^{-1}x$, so that $dx =  {\abs {\det g}} dy $
\begin{equation} = \int_{g \in GL_n(L)}  \int_{ y \in L^n}\phi(g) \varphi(y) \norm {gy}^{s-n} \abs {\det g} d^\times g dy \end{equation}
replacing $g$ with $g^{-1}$
\begin{equation} = \int_{g \in GL_n(L)}  \int_{ y \in \QQ_v^2}\inv {\abs {\det g}} \phi(g^{-1}) \varphi(y) \norm {g^{-1}y}^{s-n} d^\times g dy \end{equation}
\begin{equation} = \int_{g \in GL_n(L)} \varphi(y) \int_{ y \in L^n}\phi^*(g)  \norm {g^{-1}y}^{s-n} d^\times g dy \end{equation}
\begin{equation} = \int_{y \in L^n} \varphi(y) (\lambda (\phi^*) \norm x^{s-n} )(y)dy \end{equation}
It is immediate that if $\phi$ is in $\cH(GL_n(L))$, then $\phi^*$ is also in  in $\cH(GL_n(L))$ so that we can apply proposition \ref{definitionxi} 
\begin{equation} = \xi_{s-n}(\phi^*)  \int_{g \in GL_n(L)} \varphi(y) \norm y^{s-n} dy \bb \end{equation}
We are now in a position to extend the definition of the  Mellin transform to second degree characters defined on vector spaces. 
\begin {definition}
Let's consider a non degenerate second degree character  f defined  on $L^n$. Choose some function $\phi$ $\in {\cH(GL_n(L))}$  so that $\xi_{s-n}(\phi^*)  \neq 0$ . We define the weak Mellin  transform $M(f,s)$  of f by the formula 
\begin{equation}  M( \lambda(\phi)f,s) = \xi_{s-n}(\phi^*) M(f,s) \end{equation}
this quantity does not depend on the choice of $\phi$ .
\end {definition}
Proof : We have to prove that  the definition does not depend on the choice of $\phi$.
Let's first suppose that L is a local field and consider some function $\phi \in \cH(GL_n(L))$. It is immediate that $\phi \star \cf_{K} = \phi$  ( we assume that the haar measure on $GL_n(L)$ is normalized so that the measure of $K$ is equal to 1), so that we have 
\begin{equation} M( \lambda(\phi) f) = M ( \lambda(\phi) \lambda(\cf_K) f) \end{equation}
applying proposition \ref{transferxi}, we get 
\begin{equation} = \xi_{s-n}(\phi^*) M(\lambda(\cf_K)f,s) \end{equation}
And we observe that $ M(\lambda(\cf_K)f,s)$ does not depend on $\phi$.
Let's now suppose that $L = \RR$ ot $\CC$.
Let's choose some function $\phi_1$ so that $ \xi_{s-n}(\phi_1^*) \neq 0$ and consider the family of functions $f_{b}(x) = f(x)\psi(b.x)$ and assume that the functions $M_1(f_b,s)$ are defined for all b in $L^n$ by the formula
\begin{equation} M(\lambda(\phi_1)f_{b},s) = \xi_{s-n}(\phi_1^*)M(f_{b},s) \end{equation}
The computations described in  proposition \ref{fourier} can be generalized to vector spaces, showing that $\zeta_{f_b}(s)$, considered as a function of $b$ ( or, more precisely, as a distribution on the variable $b$), is the weak Fourier transform of the distribution $f(x) \norm x^{s-n}$, which shows unicity because  $\zeta_{f_b}(s)$ is a continuous function of $b$.$\bb$

\begin {prop}
If $\phi$ is in $C_c^\infty(GL_n(L))$, and f is a Schwartz function, the Fourier transform of the function $\lambda(\phi)f$ is equal to $\lambda( \phi^c)\four(f)$ where $\phi^c$ is defined by the formula $\phi^c(g) = \phi^*(g^t) =  \inv {\abs {\det g}} \phi((g^t)^{-1})$

\end {prop}

Proof : 
We write that 
\begin{equation} \four( \lambda(\phi)f)(y) = \int_{x \in L^n} \int_{g \in GL_n(L)}\phi(g) f(g^{-1}x) \psi(x.y) d^\times g dx \end{equation}
writing $g^{-1}x = z$, we get 
\begin{equation} =  \int_{z \in L^n} \int_{g \in GL_n(L)}\phi(g) f(z) \psi((gz).y) \abs {\det g} d^\times g dz \end{equation}
\begin{equation} =  \int_{z \in L^n} \int_{g \in GL_n(L)}\phi(g) f(z) \psi(z.g^ty) \abs {\det g} d^\times g  dz \end{equation}
\begin{equation} = \int_{ g \in GL_n(L)} \phi(g) \abs {\det g} \four(f)(g^ty) d^\times g \end{equation}
writing $h = (g^t)^{-1}$
\begin{equation} =  \int_{ h \in GL_n(L)} \phi((h^t)^{-1})\inv { \abs {\det h}} \four(f)(h^{-1}y) d^\times h \end{equation}
\begin{equation}\bb\end{equation}
\begin {prop}
The proposition is also valid is f is a non degenerate second degree character
\end {prop}
Proof : considering that $\lambda(\phi)f$ is a Schwartz function, it is enough prove this in the weak sense, using the same method as in proposition \ref{prop fourier lambda phi}. The computation is straightforward $\bb$
\begin {prop}
The weak Mellin  transform of a second degree character satisfies  for $0< \Re(s)< n$ the equation 
\begin{equation} \zeta_{f}(s) = \rho_n(s) \zeta_{\four(f)}(n-s)  = \rho_n(s)\frac {\gamma_f} {\sqrt {\abs \varrho}}\zeta_{\bar f ( \varrho^{-1}x)}(n-s)\end{equation}
where the scalar factor $\rho_n(s)$ has been defined in proposition \ref{LFT vector space}
\end {prop}
Proof : 
Let's first consider a Schwartz function $\varphi$ and any function  $\phi$ $\in {\cH(GL_n(L),\pi)}$  so that $\xi_{s}(\phi^*)  \neq 0$. 
$\lambda(\phi)\varphi$ is again a Schwartz function, so that we have by the local functional equation (proposition \ref{LFT vector space}) the equality 
\begin{equation} M( \lambda(\phi)\varphi,s)  =  \rho_n(s) M(\four( \lambda(\phi)\varphi),n-s) \end{equation}
using the previous proposition,  
\begin{equation} M(\lambda(\phi)\varphi,s)  =  \rho_n(s) M(\lambda(\phi^c) \four(\phi)\varphi,n-s) \end{equation}

If $\phi \in {\cH(GL_n(L))}$, then $\phi^c$ is also in $\cH(GL_n(L))$, so that we can use proposition \ref{transferxi} and we get, introducing the notation $\phi^t(g) = (\phi^c)^*(g) = \phi(g^t)$, 
\begin{equation} \xi_{s-n}(\phi^*) M( \varphi,s) = \rho_n(s) \xi_{-s}(\phi^t) M( \varphi, n-s) \end{equation}

If we compare this with the local functional equation for $\varphi$, we  conclude that $\xi_{s-n}(\phi^*) =  \xi_{-s}(\phi^t)$ for any function $\phi$ in $\cH(GL_n(L))$.
Let's now consider a non degenerate second degree character $f$. We know that $\lambda(\phi)f$ is a Schwartz function, so  that we have 
\begin{equation} M( \lambda(\phi)f,s) = \rho_n(s) M( \four(\lambda(\phi)f),n-s) \end{equation}
Assuming that $\phi$ is in $\cH(GL_n(L))$, using the previous propositions and the definition of the weak Mellin transform, this becomes
\begin{equation} \xi_{s-n}( \phi^*) \zeta_{ f}(s) = \rho_n(s) \xi_{-s}(\phi^t) \zeta_{\four(f)} (n-s) \end{equation}
and finally 
\begin{equation}  \zeta_{ f}(s) = \rho_n(s) \zeta_{\four(f)} (n-s) \bb \end{equation}

This formula shows that  $\zeta_f(s)$ has an analytic continuation, but is not really a functional equation. We note again, however, that if $\varrho$ is scalar, i.e. of the form $\varrho = \alpha Id$ where $\alpha$ is a scalar in $L$, then we get a true functional equation.

\subsection { The zeroes of $\zeta_f$ for second degree characters defined on $\QQ_p^n$ }
\begin {theorem}
Let's consider a non degenerate second degree character f on $\QQ_p^n$ and assume that the associated map $\varrho$ is  a dilation $\varrho = \alpha Id$.  Then the zeroes of the weak Mellin transform of f   are on the axis $\Re(s) = \frac n2$
\end {theorem}
Proof : 
We then suppose that the second degree character can be written as 
\begin{equation} f(x_1,..x_n) = \psi( \inv 2 \alpha \sum_{i = 1}^n x_1^2 + \sum_{i = 1}^n \beta_i x_i ) \end{equation}
with $\alpha \neq 0$ and $\beta_i \in \QQ_p$

Let's compute for some vector $v$ the function
$ \lambda(\cf_{GL_n(\ZZ_p)})f(v) $

We write 
\begin{equation}  \lambda(\cf_{GL_n(\ZZ_p)})f(v) = \int_{GL_n(\ZZ_p)} f(kv) d^\times k \end{equation}
\begin{equation} = \int_{GL_n(\QQ_p)} f(gv) \cf_{GL_n(\ZZ_p)}(g) d^\times g\end{equation}
We know, thanks to the unicity of the Haar measure, that $d^\times g$ is equal to $\frac {dg}{\abs {\det g}^n}$  up to some scalar factor. Let's note $\mu$ this scalar factor.
The integral becomes, noting that the determinant of any element of $GL_n(\ZZ_p)$ has to be a unit :
\begin{equation} = \mu \int_{M_n(\QQ_p)} f(gv_0) \cf_{GL_n(\ZZ_p)}(g) dg\end{equation}

Let's write $ g_{i,j}$ the matrix coefficients of the matrix $g$.
We can assume that $v = re_1 = (r,0,0..0)$  for some $r \in \QQ_p$ because the result is a radial function ( i.e. a function invariant under the action of $GL_n(\ZZ_p)$).
We have  $gv = g(r,0,0..0) = (rg_{1,1}, rg_{2,1},..rg_{n,1})$
and 
\begin{equation}f(gv_0) = \psi(\inv 2 \alpha r^2(\sum_{i = 1}^n g_{i,1}^2) + \sum_{i = 1}^n \beta_i g_{i,1} ) \end{equation} 
This expression is independent of the $g_{i,j}$ with $j \neq 1$, so that the integral can be simplified. We use the following proposition

\begin {prop}
Let's consider some matrix A  in $M_n(\ZZ_p)$ and assume the matrix elements $a_{i,1}$ of the first column of A are not all in $p\ZZ_p$ and are fixed, while the other matrix elements are considered as variables. Then the additive measure of the set of $(a_{i,j})_{j \neq 1}$ satisfying $(a_{i,j}) \in GL_n(\ZZ_p)$ is equal to the measure of $GL_{n-1}(\ZZ_p)$
\end {prop}

Proof : 
We can suppose without loss of generality that the valuation of $a_{1,1}$ is zero. Let's suppose for example that $n = 3$ and write
\begin{equation}  \begin {pmatrix}  a_{1,1} & a_{1,1}x & a_{1,1}y\\ a_{2,1}  & a_{2,1}x+z & a_2y+t \\ a_{3,1} & a_{3,1}x+u
 & a_{3,1}y+v  \end {pmatrix} =  \begin {pmatrix}  a_{1,1} & 0 & 0\\ a_{2,1}  & 1 & 0 \\ a_{3,1} & 0 & 1  \end {pmatrix} \begin {pmatrix}  1 & x & y\\ 0  & z & t \\ 0 & u & v  \end {pmatrix}\end{equation}
 Considering that the matrix $\begin {pmatrix}  a_{1,1} & 0 & 0\\ a_{2,1}  & 1 & 0 \\ a_{3,1} & 0 & 1  \end {pmatrix}$ is in $GL(\ZZ_p)$, 
The left matrix is in $GL_3(\ZZ_p)$ if and only if $x$ and $y$  are in $\ZZ_p$ and if the square matrix $\begin {pmatrix} z&t \\ u&v \end {pmatrix}$ is in $GL_{2}(\ZZ_p)$. Considering that the additive measure of $\ZZ_p$ is equal to 1, the measure of the possible vectors $(x,y,z,t,u,v)$ is then equal to the additive measure of $GL_2(\ZZ_p)$.  We then observe that the determinant of the map $(x,y,z,t,u,v) \mapsto (a_{1,1}x,a_{1,1}y,a_{2,1}x+z,a_{2,1}y+t,a_{3,1}x+u,a_{3,1}y+v)$ is equal to $a_{1,1}^2$ which is a unit. The proof for a general $n$ is the same $\bb$\\

Let's note $\kappa$ the measure of $GL_{n-1}(\ZZ_p)$ : the integral becomes
\begin{equation} \mu \kappa \int_{ g_{1,1}.. g_{n,1} \in D }  \psi(\inv 2 \alpha r^2(\sum_{i = 1}^n g_{i,1}^2) + \sum_{i = 1}^n \beta_i g_{i,1} )dg_{1,1}.. dg_{i,1} \end{equation} 
where the domain of integration D of the integral is the set of all vectors $ g_{i,1}.. g_{i,n}$ in $\ZZ_p^n$ so that at least one of the $g_{i,1}$ is a unit.
This domain can be expressed as the set of vectors which lie  in $ \ZZ_p^n$ but not in  $ (p\ZZ_p)^n$, so that the integral can be written as 
 \begin{multline} \mu \kappa \int_{ g_{1,1}.. g_{n,1} \in \ZZ_p^n }  \psi(\inv 2 \alpha r^2(\sum_{i = 1}^n g_{i,1}^2) + \sum_{i = 1}^n \beta_i g_{i,1} )dg_{1,1}.. dg_{i,1} \\
 - \mu \kappa \int_{ g_{1,1}.. g_{n,1} \in (p\ZZ_p)^n }  \psi(\inv 2 \alpha r^2(\sum_{i = 1}^n g_{i,1}^2) + \sum_{i = 1}^n \beta_i g_{i,1} )dg_{1,1}.. dg_{i,1} \end{multline} 
 It is then natural to introduce the function 
 \begin{equation} \theta_f (r) = \int_{ g_{1,1}.. g_{n,1} \in \ZZ_p^n }  \psi(\inv 2 \alpha r^2(\sum_{i = 1}^n g_{i,1}^2) + \sum_{i = 1}^n \beta_i g_{i,1} )dg_{1,1}.. dg_{i,1} \end{equation} 
 so that the integral can be written as 
 \begin{equation} \mu \kappa ( \theta(r) - \inv {p^n} \theta(pr))\end{equation}

We then get the equality 
\begin{equation}  \lambda(\cf_{GL_2(\ZZ_p)})f(v) = \mu \kappa ( \theta(r) - \inv {p^n}\theta (pr)) \end{equation}
Considering that $\cf_{GL_2(\ZZ_p)} \star \cf_{GL_2(\ZZ_p)} = \cf_{GL_2(\ZZ_p)}$ so that $\xi_s(\cf_{GL_2(\ZZ_p)}) = 1$ for all values of s, the weak Mellin  transform of $f$ is simply the Mellin transorm of $\lambda( \cf_{GL_2(\ZZ_p)})f$ : 
\begin{equation} \zeta_f(s) =  \int_{ v \in \QQ_p^n}  \lambda(\cf_{GL_2(\ZZ_p)})f(v) \norm v^s \frac {dv} { \norm v^n}\end{equation}
this integral is equal, up to a scalar factor, to 
\begin{equation} =  \int_{ r  \in \QQ_p }  \lambda(\cf_{GL_n(\ZZ_p)})f(r e_1) \abs r ^s \frac {dr} { \abs r^n}\end{equation} 
this is equal, up to a scalar factor, to 
\begin{equation} =  \int_{r \in \QQ_p^*}  ( \theta_f(r) - \inv {p^n}\theta_f (pr)) \abs r^s d^\times r \end{equation}
\begin{equation} =    (1 - \inv {p^{n-s}}) \Mell ( \theta_f,s) \end{equation}
Let's now compute the Mellin transform of $\theta_f(r)$. 
We remark that the integral definition of $\theta_f$ naturally splits as the product of one dimensional integrals $\theta_{f_i}$ which we have already computed  in the proof of  theorem \ref{zerounramified}

In order to complete the proof of the theorem, we then have to prove the following proposition 
\begin {prop}
Let's consider n non degenerate second degree characters $f_1$,.. $f_n$ on $\QQ_p$ and note $\theta_{f_1}..\theta_{f_n}$ the associated functions. Assume that the endomorphisms $\varrho_i$ associated to $f_i$ all have the same modulus. Then the zeroes of the Mellin transform of  the product $\theta_{f_1} (r) \theta_{f_2}(r)..\theta_{f_n}(r)$ are on the axis $\Re(s) = \frac n2$
\end {prop}
Proof : let's suppose for example that the valuation of $\varrho$ is even. After rescaling of the functions $f_i$, we can then assume that the valuation of $\varrho$ is zero. Then we have seen in the proof of theorem  \ref{zerounramified} that all the functions $\theta_{f_i}$ are of the form 
\begin{equation} \cf_{\ZZ_p}(x) + (1- \cf_{\ZZ_p}) \inv {\abs x}\end{equation}
or for $k \ge 1$ 
\begin{equation} \cf_{p^k\ZZ_p}(x) + \gamma_f \inv {\abs x} \cf_{p^k\ZZ_p}(\inv x) \end{equation}

We remark that on each formula, the left term cancels for negative valuations of x and the right term cancels for positive or zero valuations of x. If all the terms satisfy $k = 0$, then the proof is immediate. Let's note $n_k$ the number of terms $\theta_{f_i}$ associated to some k. The product of the functions $\theta_{f_i}$ becomes 
\begin{equation} \lbrace (\cf_{\ZZ_p}(x))^{n_0}  \prod_{ i \ge 1} ( \cf_{p^k\ZZ_p}(x))^{n_i}\rbrace  + \gamma \lbrace ((1- \cf_{\ZZ_p}) \inv {\abs x})^{n_0} \prod_{i \ge 1} ( \inv {\abs x} \cf_{p^k\ZZ_p}(\inv x))^{n_i}\rbrace \end{equation} 
where $\gamma$ is the product of all the $\gamma_f$
Let's note m the maximum of the $k's$ appearing with non zero $n_k$.  The first term simplifies to $\cf_{p^m\ZZ_p}$ and, assuming $m\ge 1$ the last term also simplifies as 
\begin{equation} \gamma \inv {\abs x^{n_0+..+n_m} }\cf_{p^m\ZZ_p} ( \inv x) \end{equation}
It is immediate that $n_0+..+n_m = n$ so that we get 
\begin{equation} \cf_{p^m\ZZ_p} + \gamma \inv {\abs x^n} \cf_{p^m\ZZ_p} ( \inv x)\end{equation}
and it is immediate, using the same method as for the case $n = 1$,  that the zeroes of the Mellin transform of this function are on the axis $\Re(s) = \frac n2$.  The proof for odd valuations of $\varrho$ is similar $\bb$

\subsection { The weak Mellin transform of a  second degree character defined on a real  vector vector space }

It is possible to give an explicit description of the weak Mellin transform of the function $\psi(\frac a2 \norm x^2 + b.x)$ on $\RR^n$.

\begin {prop}
\label {trivialrep}
The weak Mellin transform of $f(x) = \psi(\frac a2 \norm x^2 + b.x)$ on $\RR^n$ with $a>0$ is 
\begin{equation} \zeta_{a,b}(s) =   \frac {e^{-\frac {\pi i }4 s}} {\sqrt a^s}\frac {\pi^{\frac {n}2}}{\Gamma( \frac {n}2)}\frac {\Gamma(\frac s2)}{\pi^\frac s2} {_1 F _1}(\frac s2,\frac n2,\frac {\pi i \norm b ^2}a) \end{equation}
\end {prop}
Remark :  this function cancels only for $\Re(s) = \frac n2$ ( cf proposition \ref{zero})

Proof : 
the method is the same as for the case $n = 1$ : the same argument shows that $\zeta_{a,b}(s)$ can be considered also with $a \in \CC$ with $\Im(a) <0$ and $b \in \CC$, and that with this definition $\zeta_{a,b}$ is a continuous in $a$ and analytic in $b$. Let's then  first compute $\zeta_{a,b}(s)$ for $b = 0$, then write its Taylor expansion as a function of $b$.
Let's first suppose $b = 0$ and $a>0$
Let's write $a' = a -i\epsilon $  with $\epsilon >0 $ and compute
\begin{equation} \int_{\RR^n} e^{-2\pi i ( \frac {a'}2 \norm x^2)} \norm x^s \frac {dx}{\norm x^n} \end{equation}
let's write $x = r u$ where $u$ is on the unit sphere $S_{n-1}$.
\begin{equation} = \int_{S_{n-1}} \int_0^\infty e^{-2\pi i \frac {a'}2 r^2} r^{s} r^{n-1}\frac {dr}{r^n}du \end{equation}
Considering that the area of $S_{n-1}$ is $\frac {2 \pi^{\frac {n}2}}{\Gamma( \frac {n}2)}$, we get 
\begin{equation}\frac { \pi^{\frac {n}2}}{\Gamma( \frac {n}2)} \inv {\sqrt {a'}^s} e^{-\frac {\pi i }4 s} \frac {\Gamma(\frac s2)}{\pi^\frac s2} \end{equation}
We then get the result for $a$ real by taking the limit when $\epsilon \rightarrow 0$.\\
Let's now suppose that $b \neq 0$, and consider $\zeta_{a,b}(s)$ as a function of $b$. It is clear that it is even, and that it is also radial. 

It can then we written as a function of $\norm b = r$ as 

\begin{equation} \zeta_{a,b}(s) = \sum_{k \ge 0} a_k r^{2k} \end{equation}
Let's note $\Delta$ the laplacian in $\RR^n$. The values of $a_k$ can be described thanks to the formula $\Delta r^{2k} = 4k(\frac n2+(k-1)) r^{2k-2}$ : 

\begin{equation} \Delta^k r^{2k} = 4^k k! ( \frac n2)_k\end{equation}
so that 
\begin{equation} a_k = \frac {\Delta_b^k \zeta_{a,b}(s)}{ 4^k k! ( \frac n2)_k} \end{equation}
Where the symbol $\Delta_b$ is the laplacian with respect to the variables $b_1,..b_n$
In order to compute $\Delta_b^k \zeta_{a,b}(s)$, we remark that the seconde degree character $f_{a,b}(x) = e^{-2\pi i ( \frac a2 \norm x^2 + b.x)}$ satisfies the partial differential equations
\begin{equation} \frac {\partial }{\partial a} f_{a,b} = - \pi i  \norm x^2 f_{a,b} \end{equation}
and 
\begin{equation} \Delta_b f_{a,b} = (-2\pi i )^2 \norm x ^2f_{a,b} \end{equation}
so that 
\begin{equation} \Delta_b f_{a,b} = (-4\pi i ) \frac {\partial f_{a,b}}{\partial a}\end{equation}
It is not difficult to show, as in proposition \ref{exchange}, (or by taking $a \in \CC$ with $\Im(a)<0$ and taking the limit when $a$ converges to a real value)  that we can exchange the integration and derivation signs 
 
We then get 

\begin{equation}  \Delta_b \zeta_{a,b}= (-4\pi i ) \frac {\partial \zeta_{a,b} }{\partial a}\end{equation}
We have for $b = 0$ the identity
\begin{equation}  \frac {\partial ^k\zeta_{a,b} }{\partial a^k}( b = 0) = \frac { \pi^{\frac {n}2}}{\Gamma( \frac {n}2)} (-1)^k (\frac s2)_k\inv {a^{\frac s2+k}} e^{-\frac {\pi i }4 s} \frac {\Gamma(\frac s2)}{\pi^\frac s2} \end{equation}
so that 
\begin{equation} \zeta_{a,b}(s) =  \frac {e^{-\frac {\pi i }4 s}}  {\sqrt a^s}\frac { \pi^{\frac {n}2}}{\Gamma( \frac {n}2)}\frac {\Gamma(\frac s2)}{\pi^\frac s2} \sum_{k \ge 0}\inv { 4^k k! ( \frac n2)_k}(-1)^k (-4\pi i )^k (\frac s2)_k\inv {a^{k}}     r^{2k}\end{equation} 
We recognize again a confluent hypergeometric function 
\begin{equation} = \frac {e^{-\frac {\pi i }4 s}}  {\sqrt a^s}\frac { \pi^{\frac {n}2}}{\Gamma( \frac {n}2)}\frac {\Gamma(\frac s2)}{\pi^\frac s2} {_1 F _1}(\frac s2,\frac n2,\frac {\pi i r^2}a) \bb\end{equation}

\subsection { weak Mellin transforms associated to non trivial representations of K}

 The  Mellin transform can be used to  decompose a radial function as an integral of functions of the form $\norm x^s$. It is clear that if one wants to get a complete decomposition of a function defined on $\RR^{n}- \{0\}$, one has also to consider functions which are non constant on the unit sphere, and the most natural way to do this is to use the theory of spherical harmonics, i.e. to consider scalar  integrals of the form 
\begin{equation}  \int_{L^n} \varphi(x) Y(\frac {x}{\norm x}) \norm x^s  \frac {dx}{{\norm x}^n} \end{equation}
where Y is some spherical harmonics on the sphere.
Using these scalar integral does not allow to define a weak Mellin transform, but we know that spherical harmonics  are associated to a special class of representation of $O(n)$ called the spherical representations, and we now show how these representations can be used to define the weak Mellin transform as a vector valued integral.

Let's  consider the vector space $L^n$, where $L$ is any locally compact field, note $e_1$ the vector $(1,0,0..)$ and  $K_{e_1}$ the stabilizer of $e_1$ in K, i.e. the subgroup of elements $k$ of $K$ so that $ke_1 = e_1$.
We say that an irreductible representation $(\pi,V_\pi)$ of $K$ is spherical if it has, up to a scalar factor, a unique vector fixed under the action of $K_{e_1}$.

The following proposition can be considered as a generalization of proposition $\ref {radial}$  : 

\begin {prop}
\label {spherical}
Let's consider a spherical representation $(\pi,V_\pi)$ of $K$ and some $s \in \CC$. Then there exists, up to a scalar factor one and only one function $\nu_{s,\pi}$ defined on $L^n - \{0\}$ with values in $V_\pi$ and  satisfying the two following condistions : 
\begin {itemize}
\item  $\nu_{s,\pi}(\mu x) = \abs \mu^s \nu_{s,\pi}(x)$ for all $\mu \in \RR_+^*$ if $L = \RR$ or $L = \CC$ and $\mu = \varpi$ if $L$ is local,  
\item  $\nu_{s,\pi}(kx) = \pi_{s,\pi}(k) \nu_{s,\pi}(x) $ for all k in $K$
\end {itemize}
\end {prop}

Proof : Let's first prove unicity. Note $v_0$ a vector of V invariant under the action of $K_{e_1}$ by the representation $\pi$.
Considering that the action of $K$ on the unit sphere $\norm x = 1$ is transitive, it is immediate that the restriction of  the function $\nu_{\pi,s}$ on the sphere $\norm x = 1$ is uniquely defined by $f(w)$ where $w$ is any vector in the sphere, for example $w = e_1 =  (1,0,0..)$, by the formula 
\begin{equation} f(ke_1) = \pi(k)f(e_1)\end{equation}

 $e_1$ is invariant under the action of $K_{e_1}$. As a consequence,  $f(e_1)$ should then be a vector in $V_\pi$  invariant under the action of $K_{e_1}$  by the representation $\pi$.   We have supposed, however, that there exists, up to a scalar factor, exactly one vector satisfying this property. We then have $\nu_{\pi,s}(e_1) = \mu v_0$ for some scalar $\mu$ and the restriction of  $\nu_{\pi,s}$ to the sphere can be described by the formula
\begin{equation} \nu_{\pi,s} (ke_1) = \pi(k) v_0\end{equation}

inversely,  this equation gives a well defined function $\nu_{\pi,s}$ on the sphere which satisfies the conditions of the proposition. The extension of $\nu_{\pi,s}$ to $L^n - {0}$ using dilations is immediate $\bb$ \\
If $\phi$ is any smooth function with compact support from $GL_n(L)$ to $\End_\CC(V_\pi)$ and $f$ any function defined on $L^n$ with values in $V_\pi$  or $\CC$
we can again define the function $\lambda(\phi)f $ by the formula
\begin{equation} \lambda(\phi)f (x) = \int_{GL(\RR^n)} \phi(g) f(g^{-1}x) dg  \end{equation}
If $f$ has values in $\CC$, $\lambda(\phi)f$ has values in $\End_\CC(V)$. If $f$ has values in $V_\pi$, $\lambda(\phi)f$ has values in $V_\pi$.

For example, if  $f$ is a non degenerate second degree character defined on $L^n$ with values in $\CC$, then the function $\lambda(\phi)f$ has values in $\End_\CC(V)$ and is a Schwartz function ( because each of the matrix coefficient is a Schwartz function).

We also have to replace the  Hecke algebra $\cH(GL_n)$ by the $\pi$-spherical Hecke algebra $\cH(GL_n,\pi)$  associated to the representation $\pi$ , i.e. the set of smooth functions from $GL_n(L)$ to $\End(V_\pi)$ having compact support and satisfying the relation $\phi(k_1gk_2) = \pi(k_1)\phi(g)\pi(k_2)$  . Using these definitions, the generalization of proposition \ref {definitionxi} is immediate :

\begin {prop}
Let's consider a function $\phi \in \cH(GL_n,\pi)$. Then 
For all $s \in \CC$,  there exists a scalar $\xi_{s,\pi} (\phi)$  so that  we have 
\begin{equation} \lambda(\phi) \nu_{s,\pi} = \xi_{s,\pi}(\phi) \nu_{s,\pi} \end{equation}

\end {prop}

Proof : it is immediate that  $\lambda(\phi) \nu_{s,\pi}$ satisfies the conditions of the previous proposition : 

\begin{equation} \lambda(\phi) \nu_{\pi,s}(\mu kx) =  \int_{GL(\RR^n)} \phi(g) \nu_{\pi,s} (g^{-1}\mu k x) dg \end{equation}
\begin{equation} =  \abs {\mu} ^s \int_{GL(\RR^n)} \phi(g) \nu_{\pi,s} (g^{-1} k x) dg \end{equation}
writing $g^{-1}k = g'^{-1}$
\begin{equation} =  \abs \mu ^s \int_{GL(\RR^n)} \phi(kg') \nu_{\pi,s} ((g')^{-1}x) dg \end{equation}
\begin{equation} = \abs \mu ^s \pi(k) \int_{GL(\RR^n)} \phi(g') \nu_{\pi,s} ((g')^{-1}x) dg \end{equation}
 this function is then equal to $\nu_{\pi,s}$ on $L^n - \{0\}$ up to a scalar factor $\bb$

Using this proposition, the method used in the previous sections to define the weak Mellin transform can be easily extended, replacing the function $\norm x^s$ with the functions $\nu_{\pi,s}$  : We simply write 

\begin{equation}M( f,s,\pi) =  \int_{\RR^n} f(x) \nu_{s,\pi}(x) \frac {dx}{\norm x^n}  \in V_\pi \end{equation}

Note that the formula defining $M(f,s,\pi)$ stil makes sense ( the integral is absolutely convergent for $\Re(s)>0$)  if $f$ has values in $GL(V)$, considering that this linear map acts on the vector $\nu_{\pi,s}(x)$ so that if $f$ second degree character and $\phi$ any smooth function with compact support in $GL_n(L)$ with values in $GL(V)$, $M(\lambda(\phi)f,\pi,s)$ is a well defined vector in $V_\pi$

Using this definition of $M(f,s,\pi)$, the definition of the weak Mellin transform of a second degree character as an element of $V_\pi$ can be done using exactly the same method as for the case $n = 1$. 

\begin {prop}
We have for all Schwartz functions defined on $\RR^n$ with values in $\CC$ the formula 
\begin{equation} M( \lambda(\phi)f, \pi,s) = \xi_{s-n,\pi}(\phi^*) M( f,\pi,s) \end{equation}
\end {prop}
Remark : note that $f$ has values in $\CC$, but $\lambda(\phi)f$ has values in $GL(V)$.\\
Proof : We have 
\begin{equation} M( \lambda(\phi)f, \pi,s) = \int_{\RR^n} (\lambda(\phi)f)(x) \nu_{s,\pi}(x) \frac {dx}{\norm x^n} \end{equation}
\begin{equation} = \int_{x \in \RR^n}\int_{g \in GL_n(\RR)} \phi(g)f(g^{-1}x) \nu_{s-n,\pi}(x)  d^\times g  {dx}   \end{equation}
writing $y = g^{-1}x$ 
\begin{equation} = \int_{y \in \RR^n}\int_{g \in GL_n(\RR)}\phi(g)f(y) \nu_{s-n,\pi}(gy) \abs {\det g}d^\times g dy \end{equation}
replacing $g $ with $g^{-1}$
\begin{equation} = \int_{y \in \RR^n}\int_{g \in GL_n(\RR)}f(y)\phi(g^{-1}) \nu_{s-n,\pi}(g^{-1}y) \inv {\abs {\det g}} d^\times g dy  \end{equation}
\begin{equation} = \int_{y \in \RR^n} f(y) \lambda(\phi^*) \nu_{s-n,\pi}(y)  dy \end{equation}
\begin{equation} = \xi_{s-n,\pi}(\phi^*) \int_{y \in \RR^n} f(y) \nu_{s-n,\pi}(y) dy \end{equation}
using this formula, we can define in a reasonable way the weak Mellin transform of a second degree character.

\begin {definition}
Let's consider a non degenerate second degree character  f defined  on $L^n$. Choose some function $\phi$ $\in \cH (GL_n,\pi)$  so that $\xi_{s-n,\pi}(\phi^*)  \neq 0$ . We define the weak Mellin  transform $M(f,s,\pi)$  of f by the formula 
\begin{equation}  M( \lambda(\phi)f,s,\pi) = \xi_{s-n,\pi}(\phi^*) M(f,s,\pi) \end{equation}
this quantity does not depend on the choice of $\phi$ .
\end {definition}
Proof : the proof is the same as the unramified case ( replace $\cf_K$ by $\cf_K(k)\pi(k)$ for the local case)

This weak Mellin transform satisfies the same kind of scaling properties as the usual Mellin transform : 

\begin {prop}
\label {scaling}
Let's suppose that the weak Mellin transform $M(f,s,\pi)$ of a function f on $L^n$  with values in $\CC$ is well defined for some s and $\pi$, and consider some $k \in K$ and $\mu \in \RR_+^*$ if L is $\RR$ or $\CC$ ( or $\mu = \varpi^k$ for some k in $\ZZ$ if L is local) Then we have the formula
\begin{equation} M( f(k\mu x), \pi,s) = \abs \mu^{-s} \pi(k)^{-1} M(f,\pi,s) \end{equation}
\end {prop}
Proof :

 The definition of the weak Mellin transform of $f$  is that for all $\phi$ in $\cH(GL_n,\pi)$, we have 
\begin{equation} M( \lambda(\phi)f,s,\pi) = \xi_{s-n,\pi}(\phi^*) M(f,s,\pi)\end{equation}

Let's note $f_{\mu}$ the function $f(\mu x)$ and $f_k$ the function $f(kx)$. It is then enough to prove that 
\begin{equation} M( \lambda(\phi)f_{\mu}, s,\pi) = \mu^{-s} M(\lambda(\phi) f,s,\pi) \end{equation}
\begin{equation} M( \lambda(\phi)f_{k}, s,\pi) =  \pi(k)^{-1} M(\lambda(\phi) f,s,\pi) \end{equation}
The first identity is immediate.
The second is the consequence of the following computation : 
\begin{equation} M( \lambda(\phi)f_{k}, s,\pi) = \int_{\RR^n}( \int_{G} \phi(g)f( k g^{-1}x)d^\times g) \nu_{s,\pi}(x) \frac {dx}{\norm x^n}  \end{equation}

writing  $h^{-1}  = kg^{-1}$, we get  
\begin{equation} = \int_{\RR^n}( \int_{G} \phi(hk)f(h^{-1}x)  d^\times g ) \nu_{s,\pi}(x ) \frac {dx}{\norm x^n} \end{equation}
using the right equivariance of $\phi$
\begin{equation} = \int_{\RR^n} (\int_{G} \phi(h) f(h^{-1}x)d^\times g  ) \pi(k)  \nu_{s,\pi}(x ) \frac {dx}{\norm x^n}   \end{equation}
using the functional property of $\nu_{s,\pi}$
\begin{equation} = \int_{\RR^n}( \int_{G} \phi(h) f(h^{-1}x)d^\times g ) \nu_{s,\pi}(kx ) \frac {dx}{\norm x^n}   \end{equation}
\begin{equation} = \int_{\RR^n} \lambda(\phi)f(x) \nu_{s,\pi}(kx ) \frac {dx}{\norm x^n}  \end{equation}
writing $y = kx$
\begin{equation} = \int_{\RR^n} \lambda(\phi)f(k^{-1}y) \nu_{s,\pi}(y ) \frac {dy}{\norm y^n}  \end{equation}
\begin{equation} = \int_{\RR^n} (\int_{G} \phi(h) f(h^{-1}k^{-1}y) d^\times h) \nu_{s,\pi}(y ) \frac {dy}{\norm y^n}  \end{equation}
writing $h^{-1}k^{-1} = g^{-1}$
\begin{equation} = \int_{\RR^n} (\int_{G} \phi(k^{-1}g) f(g^{-1}y) d^\times g) \nu_{s,\pi}(y ) \frac {dy}{\norm y^n}  \end{equation}
using the left equivariance of $\phi$

\begin{equation} =  \pi(k)^{-1} M(\lambda(\phi)f,s,\pi) \bb\end{equation} $\bb$

We now specialize to the case $L  = \RR$
\begin {prop}
Let's assume that $L = \RR$  and that $(\pi,V_\pi)$ is a spherical representation of $K_L$. Then the  weak Fourier transform of $\nu_{\pi,s}$ is equal to $\nu_{\pi,-s-n}$ up to a scalar factor
\end {prop}
 Proof : We know ( Cf \cite{Grafakos}, p130) that if $\nu$ is a $C^\infty$ function on $\RR^n-\{0\}$, that is homogeneous, then its Fourier transform $\four(\nu)$ is also a $C^\infty$ function on $\RR^\infty - \{0\}$. It is immediate, using the commutation relation of the Fourier transform, that this fourier transform satisfies $\four(\nu)(\mu x) = \mu^{-n-s} \four(\nu)$ and $\four(\nu)(kg) = \pi( (k^{-1})^t) \nu (g) = \pi( k) \nu (g)$

 Using this proposition, it is not difficult to see that the weak Mellin transform of a second degree character of the form $\psi(\frac a2 x.x + b.x)$ satisfies a functional equation. Let's now consider the location of the zeroes of this function. 
 We first remark that   on $\RR^n$ the function $\zeta_f(s,\pi)$ considered as a function of $s$, is vector valued, but  behaves like a scalar. Let's for example consider the second degree chatacter $f_{a,e_1} =  \psi(\frac a2 x.x + e_1.x)$. We have 
 \begin {prop}
 For any value of s, $\zeta_{f_{a,e_1}}(s)$ is equal, up to a scalar factor, to the unique vector $v_0$ in $V_\pi$ fixed under the action of $K_{e_1}$
 \end {prop}
 Proof : 
 the function $f_{a,e_1}(x)$ remains unchanged if we replace $x$ with $kx$ with $k \in K_{e_1}$. As a consequence of proposition \ref{scaling}, the function $\zeta_{f_{a,e_1}}(s,\pi) $ is a vector of $V_\pi$ invariant under the action of $K_{e_1}$. It is then equal to $v_0$ up to a scalar factor $\bb$\\
 
It is also possible to prove  that the zeroes of $\zeta_f$ lie on the axis $\Re(s) = \frac n2$ if the morphism associated to $f$ is a scalar : 

\begin {theorem}
Let's consider on $\RR^n$  a second degree character of the form $f_{a,b}(x) = \psi(\inv 2 a x.x + b.x)$ with $a \in \RR^*$ and $b \in \RR^n$, $(\pi,V_\pi)$ an irreductible spherical representation of $K = O(n)$,  and note $\zeta_{a,b}(s,\pi)$ the weak Mellin transform of $f_{a,b}$. 
Then 
\begin {itemize}
\item if $b = 0$ and $\pi$ is not trivial,  then $\zeta_{a,b}(s,\pi) = 0$ for all values of s.
\item  if $b \neq 0$,  all the zeroes of $\zeta_{a,b}$ lie on the axis $\Re(s) = \frac n2$
\end {itemize}
\end {theorem} 

Proof : 
The case $b = 0$ is clear : the function $\psi( \inv 2 ax.x)$ is invariant under the action of $K = O(n)$. As a consequence of proposition $\ref{scaling}$, $\zeta_{a,b}(s)$  is a vector in $V_\pi$  invariant under the action of $\pi(k)$ for all $k$ in $K$. Considering that the representation $\pi$  is assumed to be irreductible, the only possible value of  $\zeta_{a,b}(s)$ is zero.

Let's now suppose that $b \neq 0$. We remark  that if $\zeta_{a,b_0}(s) = 0$ for some $b_0$, we have also $\zeta_{a,kb_0}(s) = 0$ for $k \in O(n)$ as a consequence of proposition \ref{scaling}, so that the function cancels on the whole sphere $\norm b = \norm {b_0}$. The idea is then to use Sturm Liouville theory in $\RR^n$ with boundary conditions on the sphere $\norm b = \norm {b_0}$. 
Let's first find the partial differential equation satisfied by $\zeta_{a,b}(s)$ considered as a function of $b$.

We observe that the function $f_{a,b} = \psi( \inv 2 ax.x +b.x)$ satisfies the formula
\begin{equation} \Delta_b f_{a,b} = (-4\pi i ) \frac {\partial } {\partial a }f_{a,b} \end{equation}
where the symbol $\Delta_b$ refers to the laplacian of $f_{a,b}$ considered as a function of the vector variable $b$.

As a consequence, the function $\zeta_{a,b}(s,\pi)$ satisfies the following differential equations : 
\begin{equation}  \Delta_b \zeta_{a,b}(s,\pi)= (-4\pi i ) \frac {\partial \zeta_{a,b}(s,\pi) } {\partial a} \label{PDE}\end{equation}
We also have using the proposition \ref{scaling} for any $\lambda >0$
\begin{equation} \zeta_{ \lambda^2 u ,\lambda v}(s,\pi) =  \lambda^{-s}  \zeta_{ u,v }(s,\pi)\end{equation}
which we can also write as 
\begin{equation} \frac \partial {\partial \lambda}( \lambda^s  \zeta_{ \lambda^2 u ,\lambda v}(s)) = 0 \end{equation}
Let's develop this equation, using the notations $\frac {\partial}{\partial a}$ for the derivation of $\zeta_{a,b}(s)$ with respect to a, and $\nabla_b$ for the gradient of $\zeta_{a,b}(s)$ with respect to the vector variable $b$ : 
\begin{equation} s \lambda^{s-1}  \zeta_{ \lambda^2 u ,\lambda v}(s) + \lambda^s (2 \lambda u) \frac {\partial}{\partial a}\zeta_{ \lambda^2 u ,\lambda v}(s)+ \lambda^s \nabla_b\zeta_{ \lambda^2 u ,\lambda v}(s). v  = 0 \end{equation}
writing $\lambda^2u = a$, $\lambda v = b$,
\begin{equation} s \lambda^{s-1}  \zeta_{ a,b}(s) + \lambda^{s-1} (2 a) \frac {\partial}{\partial a}\zeta_{a ,b}(s)+ \lambda^{s-1} \nabla_b\zeta_{ a,b}(s). b  = 0 \end{equation}
using the partial differential equation $\ref{PDE}$, and removing the term $\lambda^{s-1}$, we get 
\begin{equation} s   \zeta_{ a,b}(s) - \frac {a}{2\pi i } \Delta_b \zeta_{a,b}(s) +  \nabla_b\zeta_{ a,b}(s). b    = 0 \label {PDEV}\end{equation}
In order to apply Sturm Liouville theory, we have to cancel the first order term. We then introduce the vector valued  function $\phi_{a,b}(s)$  defined by the formula
\begin{equation}\zeta_{a,b}(s) = \phi_{a,b}(s) e^{\frac {\pi i b.b}{2a}}\end{equation}

Elementary calculations show that the  equation \ref{PDEV} becomes 
\begin{equation}  \Delta_b \phi_{a,b}(s) +( \frac {\pi i }a (n-2s)  + \norm b^2 ( (\frac {\pi}a)^2) \phi_{a,b}(s) = 0 \end{equation}

Let's multiply this equation with the vector $\bar \phi_{a,b}(s)$ ( the vector whose coordinates are the complex conjugates of the coordinates of $\phi_{a,b}(s)$) and integrate on the ball B defined by  $ \norm b \le \norm {b_0}$. 
\begin{equation}0 =  \int_{b \in B} (\Delta_b \phi_{a,b}(s) +( \frac {\pi i }a (n-2s)  + \norm b^2 ( (\frac {\pi}a)^2) \phi_{a,b}(s)) . \bar \phi_{a,b}(s) db \end{equation}
using the boundary conditions, we get 
\begin{equation} = - \int_{b \in B} \norm { \nabla_b\phi_{a,b}(s)}^2 db +\int_{b \in B}( \frac {\pi i }a (n-2s)  + \norm b^2 ( (\frac {\pi}a)^2)) \norm {\phi_{a,b}(s))}^2  db \end{equation}
so that we have the equality
\begin{equation} \frac {\pi i }a (n-2s) \int_{b \in B} \norm {\phi_{a,b}(s))}^2  db =  \int_{b \in B} \norm { \nabla_b\phi_{a,b}(s)}^2 db - \int_{b \in B}\norm b^2 ( (\frac {\pi}a)^2)) \norm {\phi_{a,b}(s))}^2  db \end{equation}
and this last expression is real. The  function $\phi_{a,b}(s))$ cannot cancel on the whole ball $B$ : Considering that it is a real analytic function as a function of the coordinates of $b$, it would imply that $\phi_{a,b}(s,\pi) = 0$ for all b. We can, however, again consider the function $\zeta_{a,b}(s,\pi)$ as the fourier transform of the distribution $\psi( \frac a2 x.x) D_{\pi,s}$, where the distribution $D_{s,\pi}$ is defined by the formula
\begin{equation} <D_{s,\pi}, \varphi> = \int_{L^n} \varphi(x) \nu_{s,\pi}(x) \frac {dx}{\norm x^n} \end{equation}
 ( the  proof is similar to the one given in proposition \ref{fourier} )  which shows that it cannot be zero. 
 $n - 2s$ has then  to be imaginary, i.e. we have $\Re(s) = \frac n2$
$\bb$


\begin {thebibliography}{20}

\bibitem {Cartier64}P. Cartier \"Uber einige Integralformeln in der Theorie der quadratischen Formen, Math. Zeitschr. 84, p 93-100 (1964) 

\bibitem {Grafakos} L. Grafakos Classical Fourier Analysis, Second Edition, Graduate Texts in Mathematics, Springer GTM 249
\bibitem {Howe} R. Howe On the role of the Heisenberg group in harmonic analysis Bull. AMS, Vol 3, Nb 2 (1980)
\bibitem {Perrin} P. Perrin Repr\'esentations de Schr\"odinger, indice de Maslov et groupe M\'etaplectique, in Lecture Notes in Mathematics, 880 (1981)

\bibitem {Riemann} B. Riemann, \"Uber die Anzahl der Primzahlen unter einer gegebenen Gr\"osse (1859)
\bibitem {Stein}E. Stein, G. Weiss,  Introduction to Fourier analysis on euclidean spaces, Princeton University Press, Princeton, NJ, 1971
\bibitem {Tate} J.Tate, Fourier Analysis in Number Fields and Hecke's Zeta-Functions, Ph.D Thesis, Princeton University, Princeton, NJ,1950
\bibitem {Weil64} A.Weil, Sur certains groupes d'op\'erateurs unitaires, Acta Math, 111, (1964)
\bibitem {Weil66} A.Weil, Fonction z\^eta et distributions, S\'eminaire Bourbaki 312 (1966)

\end {thebibliography}
\end {document}